\documentclass[12pt,oneside,reqno]{amsart}
\usepackage{amssymb}
\usepackage{color}
\usepackage{tikz-cd}
\usetikzlibrary{decorations.markings}
\usepackage{xcolor}
\definecolor{deepgreen}{cmyk}{0.99998,0,1,0}
\usepackage{hyperref}
\hypersetup{hypertex=true,
	colorlinks=true,
	linkcolor=blue,
	anchorcolor=blue,
	citecolor=blue}
\usepackage{comment}
\usepackage{mathrsfs}
\usepackage{cleveref}
\usepackage[%
]{geometry}
\usepackage{cleveref}
\crefname{section}{§}{§§}
\Crefname{section}{§}{§§}
\allowdisplaybreaks[4]

\usepackage{etoolbox}
  
  \makeatletter
  \newcommand*{\rom}[1]{\expandafter\@slowromancap\romannumeral #1@}
  \makeatother
  
  \usepackage{bm}
  \usepackage{bbm}
\theoremstyle{definition}
  \newtheorem{defi}{$\mathbf{Definition}$}[section]
  \newtheorem*{pro}{$\mathbf{Proof}$}
  
  \theoremstyle{plain}
  \newtheorem{theo}[defi]{$\mathbf{Theorem}$}
  \newtheorem{lemma}[defi]{$\mathbf{Lemma}$}
  \newtheorem{coro}[defi]{$\mathbf{Corollary}$}

  \newtheorem{prop}[defi]{$\mathbf{Proposition}$}
  
  \theoremstyle{remark}
  \newtheorem{remark}[defi]{$\mathbf{Remark}$}

\newcommand{\de}{{d}}

\newcommand{\trs}{\mathrm{Tr_s}}
\newcommand{\tro}{\mathrm{Tr}}
\newcommand{\mi}{\mathrm{i}}

\newcommand{\pa}{\partial}

\newcommand{\hti}{\widehat{\otimes}}

\newcommand{\lk}{\left(}
\newcommand{\rk}{\right)}
\newcommand{\li}{\left\langle}
\newcommand{\ri}{\right\rangle}
\newcommand{\lv}{\left\vert}
\newcommand{\rv}{\right\vert}
\newcommand{\lV}{\left\Vert}
\newcommand{\rV}{\right\Vert}

\newcommand{\blk}{\big(}
\newcommand{\brk}{\big)}
\newcommand{\bli}{\big\langle}
\newcommand{\bri}{\big\rangle}

\newcommand{\Blk}{\Big(}
\newcommand{\Brk}{\Big)}
\newcommand{\Bli}{\Big\langle}
\newcommand{\Bri}{\Big\rangle}

\newcommand{\bv}{\big\vert}
\newcommand{\Bv}{\Big\vert}
\newcommand{\bbv}{\bigg\vert}

\newcommand{\bV}{\big\Vert}
\newcommand{\BV}{\Big\Vert}
\newcommand{\bbV}{\bigg\Vert}

\newcommand{\h}{\mathcal{H}}

\definecolor{pink}{RGB}{249,164,186}
\definecolor{grassgreen}{RGB}{128,255,0}

\numberwithin{equation}{section}
\title{Toeplitz operators and the full asymptotic torsion forms} 

\author
{Qiaochu Ma}
\date{}

\begin{document} 
\clearpage\maketitle
\begin{abstract}
	This paper aims to study the asymptotic expansion of analytic torsion forms associated with a certain series of flat bundles $\{F_p\}_{p\in\mathbb{N}^*}$. We prove the existence of the full expansion and give a formula for the sub-leading term, while Bismut-Ma-Zhang \cite{bmz17} have studied the first order expansion and expressed the leading term as the integral of a locally computable differential form.
\end{abstract}
\section*{Introduction}

\subsection{Backgrounds}

In the 1930s, the \emph{Reidemeister torsion} was introduced by Reidemeister \cite{Reidemeister1935HomotopieringeUL} and Franz \cite{Franz35} in their study of lens spaces. The Reidemeister torsion was the first homeomorphic invariant which is \emph{not} homopoty. Let $(F,\nabla^F)$ be a unitary flat vector bundle over a compact manifold $X$. We assume that $H^\bullet(X,F)=0$. The Reidemeister torsion is a real number obtained by a simplicial complex with values in $F$ associated with a triangulation of $X$, which turns out to be independent of the triangulation. 

Ray and Singer asked if there is an analytic interpretation of the Reidemeister torsion. They defined their \emph{analytic torsion} \cite{rs73} as an alternating product of the regularized determinant of Hodge Laplacian and conjectured that it coincides with the Reidemeister torsion for unitary flat bundles. This conjecture was proved by Cheeger \cite{che79} and Müller \cite{mu78} independently. Bismut-Zhang \cite{bz92} and Müller \cite{mu93} simultaneously considered generalizations of this result. Müller \cite{mu93} extended it to the case when $F$ is unimodular, $X$ is oriented and odd-dimensional. Bismut-Zhang \cite{bz92} generalized it to any flat vector bundle with a Hermitian metric.

In this paper, we study the asymptotic property of analytic torsions. Let $\{F_p\}_{{p\in{\mathbb{N}^*}}}$ be a certain family of flat vector bundles over $M$, we obtain the full asymptotic expansion of the analytic torsion of $F_p$ when $p\rightarrow+\infty$. 

Let us now give some background on the results of this paper.

In \cite[Example 1.3]{che79}, Cheeger observed the relation between the Reidemeister torsion of a simplicial complex and the size of its \emph{torsion subgroup} of homology. By using this relation and discussing the asymptotics of analytic torsions of quotients of symmetric spaces by a decreasing sequence of lattices in an underlying Lie group, Bergeron-Venkatesh \cite{bergeron_venkatesh_2013} studied the growth of the size of torsion elements in the homology of an arithmetic group. This was the first application of analytic torsion in arithmetic.

In \cite{mu12}, Müller considered the asymptotics of analytic torsions for a family of flat bundles over a compact $3$-dimensional hyperbolic manifold as a real analogue of Bismut-Vasserot's result on the holomorphic torsion \cite{bisvas90}. Let $\Gamma\backslash\mathbb{H}^3$ be a compact $3$-dimensional hyperbolic manifold with constant sectional curvatures $-4$. For $p\in\mathbb{N}^*$, put $F_p=\text{Sym}^p\mathbb{C}^2$, where $\mathbb{C}^2$ is the flat bundle on $\Gamma\backslash\mathbb{H}^3$ associated with the tautological representation of $\mathrm{SL}_2(\mathbb{C})$ on $\mathbb{C}^2$ and $\text{Sym}^p$ denotes the $p$-th symmetric power. Let $\mathcal{T}({\Gamma\backslash\mathbb{H}^3},\nabla^{F_p})$ be the the analytic torsion of $F_p$. Using Selberg's trace formula, Müller \cite{mu12} obtained
\begin{equation}\label{0.0.2}
	\lim_{p\to+\infty}p^{-2}\mathcal{T}({\Gamma\backslash\mathbb{H}^3},\nabla^{F_p})=\frac{2}{\pi}{\mathrm{Vol}(\Gamma\backslash\mathbb{H}^3)}.
\end{equation}

In \cite{bmz11} and \cite{bmz17}, Bismut-Ma-Zhang gave a general construction of a family of flat vector bundles $\{F_p\}_{{p\in{\mathbb{N}^*}}}$ on any compact manifold and they expressed the asymptotics of analytic torsions as an integral of a local computable differential form on the base manifold. Indeed, they worked in a general setting of the analytic torsion forms of Bismut-Lott \cite{bl95}. 

Our main result is to extend Bismut-Ma-Zhang's work to get a full expansion of torsions and give an explicit formula for the sub-leading term. Now we explain in detail.

\subsection{The existence of full expansion of torsions}
In \cite{bl95}, Bismut-Lott constructed the analytic torsion form as a family extension of the Ray-Singer torsion. Let $\pi\colon M\rightarrow S$ be a fibration of manifolds with compact fiber $X$ of dimension $m$. We set a metric $g^{TX}$ on the relative tangent bundle $TX$ and a horizontal bundle $T^HM\subset TM$ such that 
\begin{equation}\label{0.2.}
TM=T^HM\oplus TX.
\end{equation}
For any  flat vector bundle $(F,\nabla^F)$ over $M$ with a Hermitian metric $g^{F}$. Bismut-Lott's torsion form is a differential form $\mathcal{T}(T^HM,g^{TX},\nabla^F,g^{F})\in \Omega^\text{even}(S)$ and its $0$-degree component $\mathcal{T}_0(T^HM,g^{TX},\nabla^F,g^{F})$ is the function which to $b\in S$ assigns the Ray-Singer torsion of the fibre $(X_b,g^{TX_b})$ over $b$, computed using the flat bundle $(F|_{X_b},\nabla^F|_{X_b},g^{F|_{X_b}})$.

Let $N$ be a compact Kähler manifold with $\dim_\mathbb{C}N=n$, $L$ a \emph{positive} holomorphic line bundle and ${{\xi}}$ a holomorphic vector bundle on $N$. Let $G$ be a Lie group acting \emph{holomorphically} on $N$ and this action can be lifted to $L$ and  ${{\xi}}$. Let $P_G\rightarrow M$ be a flat principal $G$-bundle. Set $\mathcal{N}=P_G\times_GN$, we have a natural projection $q\colon\mathcal{N}\to M$ and the real relative tangent bundle $T_\mathbb{R}{N}=\ker q_*$. For ${p\in{\mathbb{N}^*}}$, let $L^p$ be the $p$-th tensor power of $L$. Let $(F,\nabla^F)$ be another flat vector bundle on $M$ with metric $g^F$. Put
\begin{equation}\label{0.4...}
F_p=\big(P_G\times_GH^{(0,0)}(N,L^p\otimes {{\xi}})\big)\otimes F,
\end{equation}
where $G$ acts naturally on $H^{(0,0)}(N,L^p\otimes {{\xi}})$. We summarize geometric settings as follows:
\begin{figure}[ht]
	\centering
	\begin{tikzcd}[ampersand replacement=\&, column sep=normal,row sep=normal]
		\&T_\mathbb{R}N\arrow[to=2-2]\&{L}^p\otimes{\xi}=P_G\times_G(L^p\otimes {{\xi}})\arrow[to=2-2]\arrow[to=2-3,"R_\bullet q_*"]\&\ \& \\
		N\arrow[r] \&\mathcal{N}=P_G\times_G N\arrow[to=3-2,"q"]\&{P_G\times_G H^{(0,0)}(N,L^p\otimes {{\xi}})}\arrow[to=3-2]\\
		X	\arrow[to=3-2]\&M \arrow[to=4-2,"\pi"] \&\\
		\&	S
	\end{tikzcd}
	\caption{}
	\label{fig1}
\end{figure}
We have a flat connection $\nabla^{F_p}$ induced by the flat connection on $P_G$ and $\nabla^F$. We also denote $P_G\times_GL$ (resp. $P_G\times_G\xi$) by $L$ (resp. $\xi$). Given a metric $g^{T_\mathbb{R}N}$ on $T_\mathbb{R}N$, it induces a volume form $dv_N\in\Omega^{2n}(\mathcal{N})$ along the fibre together with $P_G$. Let $g^L$ (resp. $g^\xi$) be a metic on $L$ (resp. $\xi$) over $\mathcal{N}$. Then the $L^2$-metric on $H^{(0,0)}(N,L^p\otimes {\xi})$ given by $(dv_N,g^L,g^\xi)$ and $g^F$ define a metric $g^{F_p}$ on $F_p$.

In \cite[Definition 9.13]{bmz17}, Bismut-Ma-Zhang introduced a non-degeneracy condition for $g^L$ (see \Cref{De}). Under this condition, Bismut-Ma-Zhang \cite[\S\,4.3, \S\,9.10]{bmz17} proved that, for $p\in\mathbb{N}^*$ large enough, we have $H^\bullet(M,F_p)=0$. In the rest of this section, we \emph{always} assume that $L$ verifies the non-degeneracy condition.

For $a\in\mathbb{R}$, let $\psi_a$ be the automorphism of $\Lambda^\bullet\lk T^*S\rk$ such that, if $\alpha\in \Lambda^k(T^*S)$, then $\psi_a\alpha=a^k\alpha$.  Let $o(TX)$ be the orientation line of $TX$. Our main result is the following theorem (see Theorem \ref{4.21}).
\begin{theo}\label{1.1}
There are locally computable differential forms $W^{L,{{{\xi}}}}_i\in\Omega^\bullet(M,o(TX))$ such that for any $k,\ell\in\mathbb{N}$, there exists $C>0$ such that as $p\rightarrow+\infty$, we have
	\begin{equation}\label{1.0.1}
		\Bv p^{-n-1}\psi_{1/\sqrt{p}}\mathcal{T}\blk T^HM,g^{TX},\nabla^{F_p},g^{F_p}\brk-\sum_{i=0}^kp^{-i}\int_XW^{L,{{{\xi}}}}_i\Bv_{\mathscr{C}^\ell(S)}\leqslant Cp^{-k-1}.
	\end{equation}
\end{theo}
We note that if $\dim X=m$ is odd, $\mathcal{T}\blk T^HM, g^{TX},\nabla^{F_p},g^{F_p}\brk\in\Omega^\text{even}(S)/d\Omega^\text{odd}(S)$ is a topological invariant for $p\in\mathbb{N}^*$ large, by the anomaly formula of Bismut-Lott \cite[Theorem 3.24]{bl95}, it is independent on $(T^HM,g^{TX},g^{F_p})$. Hence each $W_i^{L,\xi}$ is a topological invariant for $i\in\mathbb{N}$. In \cite[Theorem 9.32]{bmz17}, Bismut-Ma-Zhang established Theorem \ref{1.1} for $k=0$ and gave an explicit formula for $W^{L,{{{\xi}}}}_0$.

Let $\widehat{M}\to S$ be a fibre bundle with fibre $\widehat{X}$, and let $\Gamma$ be a discrete group acting fibrewise freely and properly discontinuously on $\widehat{M}$ such that $\Gamma\backslash \widehat{M}=M$. We have the $\Gamma$-torsion form $\mathcal{T}^\Gamma\blk T^HM,g^{TX},\nabla^{F_p},g^{F_p}\brk\in\Omega^\text{even}(S)$ as in \eqref{gc5} (see also \cite{gammat}, \cite{10.4310/jdg/1214448084} and \cite{MATHAI1992369}). Then $\mathcal{T}^\Gamma\blk T^HM,g^{TX},\nabla^{F_p},g^{F_p}\brk$ has the \emph{same} asymptotic expansion as in \eqref{1.0.1}, indeed, we have the following stronger result (see Theorem \ref{Gc5}).
\begin{theo}\label{t03}
The asymptotics of the two torsions differ only by an exponentially decay term: there is $a>0$ such that for $\ell\in\mathbb{N}$, there is $C>0$ that as $p\rightarrow+\infty$, we have
	\begin{equation}\label{gc.15}
		\begin{split}
			\Bv\mathcal{T}^\Gamma\blk T^HM,g^{TX},\nabla^{F_p},g^{F_p}\brk-\mathcal{T}\blk T^HM,g^{TX},\nabla^{F_p},g^{F_p}\brk\Bv_{\mathscr{C}^\ell(S)}\leqslant C{e}^{-a{p}}.
		\end{split}
	\end{equation}
\end{theo}
Theorem \ref{t03} was first established by Bismut-Ma-Zhang \cite[\S\,7.6]{bmz17} when $S=\{\mathrm{pt}\}$, a point. If $S=\{\mathrm{pt}\}$, $X=\Gamma\backslash\bm{G}/\bm{K}$, a compact locally symmetric manifold, and $F_p$ is induced by multiples of the highest weight $\lambda\in\bm{\mathfrak{u}}^*$ of an irreducible $\bm{U}$-representation, where $\bm{U}$ is a compact form of $\bm{G}$ with lie algebra $\bm{\mathfrak{u}}$, Bismut-Ma-Zhang \cite[Proposition 8.12]{bmz17} showed that the nondegeneracy condition is the same to $W_{\bm{U}}\cdot\lambda\cap\bm{\mathfrak{t}}^*=\emptyset$, where $W_{\bm{U}}$ is the Weyl group of $\bm{U}$ and $\bm{\mathfrak{t}}$ is the Lie algebra of a maximal torus in $\bm{K}$, which is exactly the \emph{strong acyclic condition} $\theta \Lambda\neq \Lambda$ \cite[Propostion \rom{2}.6.12]{borel}, Müller-Pfaff \cite[Propositions 1.2, 1.3]{mupf13} gave a new proof of \eqref{gc.15} and showed that $\mathcal{T}^\Gamma\blk g^{TX},\nabla^{F_p},g^{F_p}\brk$ is a \emph{polynomial} of $p\in\mathbb{N}^*$, see also Liu \cite[Theorem 7.4.3]{2005} for another proof of this polynomial property.

\subsection{An explicit formula for $W_1^{L,{\xi}}$ for reductive $G$}

Now we explain another main result, an explicit formula for $W_1^{L,{\xi}}$ when $G$ is a \emph{connected reductive linear} Lie group, $\mathrm{rk}\xi=1$ and $F=\mathbb{C}$ trivial.

Let $G$ be a connected reductive linear Lie group with Lie algebra $\mathfrak{g}$. Let $K\subset G$ be a maximal compact subgroup of $G$ with Lie algebra $\mathfrak{k}$. We have the Cartan decomposition $\mathfrak{g}=\mathfrak{p}\oplus\mathfrak{k}$. Let $U$ be a compact form of $G$ with Lie algebra $\mathfrak{u}=\sqrt{-1}\mathfrak{p}\oplus\mathfrak{k}$.

Recall that $N$ is a compact complex manifold with $\dim_\mathbb{C}N=n$ and $(L,g^{L})$ is a positive line bundle on $N$ with the first Chern form $c_1(L,g^L)$. We assume further that $U$ acts holomorphically on $N$ and this action lifts to a holomorphic unitary action on $L$, then $c_1(L,g^L)$ is a $U$-invariant form. Let $\mu_L\colon N\rightarrow\mathfrak{u}^*$ be the associated moment map obtained from the Kostant formula \eqref{hc1}.

Let $\theta^\mathfrak{g}$ be the $\mathfrak{g}$-valued flat connection form on $P_G$. Let $P_K$ be a reduction of $P_G$ to $K$-principal bundle. Set $\mathfrak{g}_r=P_K\times_K\mathfrak{g}$. By projection of $\theta^\mathfrak{g}$ on $\mathfrak{p}$ and $\mathfrak{k}$ with respect to the Cartan decomposition, we get $\theta^\mathfrak{g}=\theta^\mathfrak{p}+\theta^\mathfrak{k}$.

For $A\in\mathfrak{u}$, let $R_{L}(A)$ be the Duistermaat-Heckman integral
\begin{equation}
	R_{L}(A)=\int_N\exp\big({2\pi i}\langle\mu_L,A\rangle+c_1(L,g^L)\big),
\end{equation}
then we naturally extend $R_L(\cdot)$ to a holomorphic function on $\mathfrak{g}_\mathbb{C}\cong \mathfrak{u}\otimes_\mathbb{R}\mathbb{C}$.

Let $S\mathfrak{g}$ be the symmetric algebra of $\mathfrak{g}$ and we denote by $\overline{S}\mathfrak{g}$ its formal completion. Canonically, $S\mathfrak{g}$ can be identified with the algebra of real differential operators with constant coefficients on $\mathfrak{g}$. Then $\overline{S}\mathfrak{g}$ naturally acts on $R_{L}(\cdot)$ (see also \cite[\S\,1.4]{bmz17}).

Let $\{{e}_i\}_{i=1}^m$ be a local orthonormal frame of ${TX}$ with dual frame $\{{e}^i\}_{i=1}^m$. Let $\widehat{TX}$ be another copy of $TX$ and let $\widehat{\theta}^\mathfrak{p}$ be the restriction of ${\theta}^\mathfrak{p}$ on $\widehat{TX}$.  Set
\begin{equation}
	\begin{split}
\vert\widehat{\theta}^{\mathfrak{p}}\vert^2&=\sum_{i=1}^m\big(\widehat{\theta}^{\mathfrak{p}}(e_i)\big)^2\in\mathscr{C}^\infty(M,S^2\mathfrak{g}_r),\\ \widehat{\theta}^{\mathfrak{p},2}&=\frac{1}{2}\widehat{e}^i\wedge\widehat{e}^j[\widehat{\theta}^{\mathfrak{p}}(\widehat{e}_i),\widehat{\theta}^{\mathfrak{p}}(\widehat{e}_j)]\in\mathscr{C}^\infty(M,\Lambda^2(\widehat{T^*X})\otimes\mathfrak{g}_r).
	\end{split}
\end{equation}

Let $\nabla^{\mathfrak{g}_r,u}$ be the connection on $\mathfrak{g}_r$ induced by $\theta^\mathfrak{k}$. For $t\geqslant 0$, let $\sigma_t$ be a section of $\Lambda^\bullet(T^*M)\hti\Lambda^\bullet(\widehat{T^*X})\otimes S\mathfrak{g}_r$ given by
\begin{equation}\label{hd..14}
	\sigma_t=-\frac{1}{4}\big\langle R^{{TX}}{e}_i,{e}_j\big\rangle\widehat{e}^i\wedge\widehat{e}^j-\theta^{\mathfrak{p},2}+\sqrt{t}\nabla^{\widehat{T^*X}\otimes\mathfrak{g}_r,u}\widehat{\theta}^{\mathfrak{p}}+t\vert\widehat{\theta}^{\mathfrak{p}}\vert^2+t\widehat{\theta}^{\mathfrak{p},2}.
\end{equation}
Denote by $\int^{\hat{B}}\colon\Lambda^\bullet(TM)\hti\Lambda^\bullet(\widehat{TX})\to \Lambda^\bullet(TM)$ the Berezin integral (see \eqref{ab1}). Let $\varphi$ be the endomorphism of $\Lambda^\bullet(T^*M)\otimes\mathbb{C}$ which maps $\alpha\in \Lambda^k(T^*M)\otimes\mathbb{C}$ to $({2\pi i})^{-k/2}\alpha$. 

If ${\xi}$ is a trivial line bundle, we denote the form $W^{L,{\xi}}_i$ in \eqref{1.0.1} by $W^{L}_i$. Bismut-Ma-Zhang \cite[Definition 2.11]{bmz17} gave the following formula for $W^L_0$:
\begin{equation}\label{1.0.13}
	W^L_0=-\sqrt{{2\pi i}}\varphi\int_{0}^{+\infty}\int^{\widehat{B}}\frac{\theta^{\mathfrak{p}}\wedge\widehat{\theta^{\mathfrak{p}}}}{2}\exp(-\sigma_t) R_{L}(0)\frac{d t}{\sqrt{t}}\in \Omega^\bullet(M,o(TX)),
\end{equation}
where the operator $\frac{\theta^{\mathfrak{p}}\wedge\widehat{\theta^{\mathfrak{p}}}}{2}\exp(-\sigma_t)\in \Lambda^\bullet(T^*M)\hti\Lambda^\bullet(\widehat{T^*X})\otimes \overline{S}\mathfrak{g}_r$ acts on the function $R_L(\cdot)$ and we evaluate at $0$.

Now we assume that $(\xi,g^\xi)$ is a holomorphic Hermitian \emph{line} bundle on $N$. Let $(TN,g^{TN})$ be the \emph{holomorphic tangent} bundle of $N$, where $g^{TN}$ is induced by $c_1(L,g^L)$: if $A\in TN$, $g^{TN}=-ic_1(L,g^L)(A,\overline{A})$. We assume that the action of $U$ on $N$ lifts to holomorphic unitary actions on $(\xi,g^\xi)$ and $(TN,g^{TN})$ with moment maps $\mu_{{{\xi}}}$ and $\mu_{\det TN}$.

Formally, for $p\in\mathbb{N}^*$, let ${{{\xi}}^{\frac{1}{p}}}$ be the $p$-th root of $L$ and $(\det TN)^{\frac{1}{2p}}$ be the $2p$-th root of $\det TN$ at the level of cohomology and moment map (see \eqref{hd9}), then $R_{L\otimes {{{\xi}}}^{\frac{1}{p}}\otimes (\det TN)^{\frac{1}{2p}}}$ and $W^{L\otimes {{{\xi}}}^{\frac{1}{p}}\otimes (\det TN)^{\frac{1}{2p}}}_0$ are well defined even if ${{{\xi}}}^{\frac{1}{p}}\otimes (\det TN)^{\frac{1}{2p}}$ may not be.

We have the following formula for $W_1^{L,{{{\xi}}}}$ (see Theorem \ref{Hd4}).
\begin{theo}\label{1.2}
For $k=1$ in \eqref{1.0.1}, as $p\rightarrow+\infty$ we have
\begin{equation}\label{1.0.2}
	p^{-n-1}\psi_{1/\sqrt{p}}\mathcal{T}\blk T^HM,g^{TX},\nabla^{F_p},g^{F_p}\brk=\int_XW_0^{L\otimes {{{\xi}}}^{\frac{1}{p}}\otimes (\det TN)^{\frac{1}{2p}}}+\mathcal{O}\big(p^{-2}\big)\in\Omega^\bullet(S).
\end{equation}
In other words, 
\begin{equation}
	W_0^{L,{{{\xi}}}}+p^{-1}W_1^{L,{{{\xi}}}}=W_0^{L\otimes {{{\xi}}}^{\frac{1}{p}}\otimes (\det TN)^{\frac{1}{2p}}}+\mathcal{O}\big(p^{-2}\big).
\end{equation}
\end{theo}
In \Cref{H.e} we discuss a special case of Theorem \ref{1.2} when $N$ is a coadjoint orbit of $\mathfrak{u}^*$, and we use this to compute asymptotics of torsion for some concrete examples in \Cref{Hf}.

Liu considered the asymptotics of equivariant torsions for compact locally symmetric spaces \cite{MR4269601} and asymptotic torsions for compact locally symmetric orbifolds \cite{2005}. There are analogous results under holomorphic settings, Bismut-Vasserot \cite{bisvas89} studied the asymptotics of holomorphic torsion associated with increasing powers of a positive line bundle over a Kähler manifold and they gave a formula for the leading term. This asymptotic expansion played an important role in the arithmetic Hilbert-Samuel theorem of Gillet–Soulé \cite{gillet}. They also extended their results in \cite{bisvas90}, where the powers of the positive line bundle are replaced by the symmetric powers of a Griffiths positive holomorphic vector bundle. In \cite{fin}, Finski generalized \cite{bisvas89} to obtain a full asymptotic expansion and gave a formula for the second term of the expansion. Puchol \cite{pu} studied asymptotics of holomorphic torsion forms, which extended \cite{bisvas89,bisvas90} to family versions. Savale \cite{sava1,sava} obtained the asymptotics of the eta invariants using semiclassical analysis.

\subsection{Main techniques}

Now we briefly describe the techniques that we will use in the proof of Theorems \ref{1.1}-\ref{1.2} as well as the main points in this paper.

\subsubsection{Toeplitz operators}

The theory of Toeplitz operators is recalled (see \Cref{sC}). Especially, one of the key results is the growth of the exponential of a Toeplitz operator (see Theorem \ref{Cc10}), which is used to get uniform estimations for operators.

\subsubsection{The spectral gap}
Under the nondegeneracy condition (see \cref{De}), Bismut-Ma-Zhang \cite[\S\,9.10]{bmz17} showed that the Hodge-de Rham Laplacian (see \eqref{be5}) has a spectral gap: $D_X^{F_p,2}\geqslant Cp^2$ for $p\in\mathbb{N}^*$ large, this condition is vitally important in analysis, for instance, it ensures that $W^{L,\xi}_i$ in \eqref{1.0.1} is well defined (see Theorem \ref{4.21}).

\subsubsection{Analytic localization method}

The proof of Theorems \ref{1.1} and \ref{1.2} relies on a refinement of Bismut-Ma-Zhang's argument \cite{bmz17}. The strategy to get a full expansion is to study the local asymptotics of certain heat kernels. We will use the analytic localization method of Bismut-Lebeau \cite{bisleb91}. Bismut's family local index theory in the context of families \cite{bgv,bis86} also plays an important role, in particular, we apply the rescaling for two Clifford variables as in \cite[Chapter 4d)]{bz92}.

Following \cite[Chapter 11]{bisleb91} and \cite[\S\, 9.12]{bmz17},  in \cref{Fb}. we prove that the limit of the trace of certain heat kernels used in the definition of analytic torsion forms (see \eqref{fa5}) can be localized, by using the finite propagation speed of the wave operator (see \cite[Appendix D]{mm07}). And as analysis carried through in \cref{Fc}-\cref{Fh}, we use the techniques in \cite[Chapter 4]{mm07} to get estimations for high order expansions of resolvents and heat kernels.

Compared with \cite{bmz17} and \cite{pu}, to get the leading term, one only needs to get a uniform bound and the pointwise convergence for the heat traces, then apply the dominated convergence, while to give the full expansion as in \eqref{1.0.1}, we need to give precise estimations for \emph{each term} in the local expansion of the heat kernel as well as the \emph{remainder}. To do so, the most important trick is to preserve the spectral gap property for the resulting Laplacian in the localization and rescaling procedures. Hence we should carefully choose parameters at each step of the analysis, especially for the weighted norm (see \eqref{fd1}) and the corresponding elliptic regularity (see Theorem \ref{4.13}).

\subsection{The organization of the paper}

This paper is organized as follows. In \Cref{sB}, we review the main results of Bismut-Lott \cite{bl95} of analytic torsion forms. In \cref{sC}, we recall some results on the Toeplitz operator following \cite[\S\,7]{mm07} and analyze the asymptotic behavior of the exponential of a Toeplitz operator. In \cref{sD}, we recall Bismut-Ma-Zhang's definition \cite[\S\,9.7]{bmz17} for a series of flat bundles $\{F_p\}_{p\in\mathbb{N}^*}$ over $M$, which is the main geometric object in the whole paper. In \cref{sF}, we give  the full expansion of the odd superconnection form $h\big( A',g^{\Omega^\bullet(X,F_p)}\big)$ as $p\rightarrow+\infty$. In \cref{sG}, we state and prove the main theorem, which gives the existence of the full expansion of the analytic torsions and the $\Gamma$-torsions associated with $F_p$. In \cref{sH}, we consider the case where $G$ is a reductive Lie group and give an explicit formula for the sub-leading term in the asymptotics of torsion obtained in \cref{sG}.

In the whole paper, we apply the superconnection formalism of Quillen (see \cite{quillen85} and \cite[\S\,1]{bgv}). If $E=E_+\oplus E_-$ is a $\mathbb{Z}_2$-vector space, and $\tau=\pm 1$ defines the $Z_2$-grading, if $A\in\mathrm{End}(E)$, we denote by $\trs[A]$ the supertrace: 
\begin{equation}
\trs[A]=\tro[\tau A].
\end{equation}
For a multi-index $\alpha=(\alpha_1,\cdots,\alpha_k)\in\mathbb{N}^k$ and a multi-variable $Z=(Z_1,\cdots,Z_k)$, set
\begin{equation}
\vert\alpha\vert=\sum_{i=0}^k\alpha_i,\ \  Z^\alpha=Z_1^{\alpha_1}\cdots Z_k^{\alpha_k}, \ \  \frac{\pa^\alpha}{\pa Z^\alpha}=\frac{\pa^{\alpha_1}}{\pa Z^{\alpha_1}}\cdots\frac{\pa^{\alpha_k}}{\pa Z^{\alpha_k}}.
\end{equation}
We also use the Einstein summation convention, if in a term the same index appears twice, that term is assumed to be summed over all possible values of that index.

\subsection*{Acknowledgment}
This work is the main result of our PhD. thesis, which was done at Université Paris Cité. I am deeply grateful to my thesis supervisor Prof. Xiaonan Ma for his patient guidance and numerous suggestions. This project has received funding from the European Union’s Horizon 2020 research and innovation
programme under the Marie Skłodowska-Curie grant agreement No 754362.

\section{The analytic torsion forms}\label{sB}

In this section, we will summarize the main results on the analytic torsion forms following Bismut-Lott \cite{bl95}, which generalize the classical Ray-Singer analytic torsion \cite{rs73}. 

This section is organized as follows. In \cref{Ba}, we introduce the smooth fibration $\pi\colon M\rightarrow S$ and define some associated tensors. In \cref{Bb}, given a flat complex vector bundle $(F,\nabla^F)$ with a Hermitian metric $g^F$, we define a natural unitary connection $\nabla^{F,u}$ and the associated odd characteristic form. In \cref{Be}, we reinterpret the de Rham operator $d^{F}$ as a flat superconnection on the bundle $\Omega^\bullet(X,F)$ over $S$. In \cref{Bi}, we get a transgression formula for the odd closed forms $h(A',g^{\Omega(X,F)}_t)\in\Omega^\bullet(S)$ and define the analytic torsion form. In \cref{Bj}, we recall Bismut-Lott's Lichnerowicz formula associated with the odd closed form.

\subsection{A smooth fibration}\label{Ba}

Let $\pi\colon M \stackrel{}{\rightarrow}  S$ be a fibration of smooth manifolds with compact fibre $X$ of dimension $m$. Let $TX\subset TM$ be the tangent bundle to the fibers $X$. Let $T^HM\subset TM$ be a horizontal subbundle  with \eqref{0.2.} and $P^{TX}\colon TM\rightarrow TX$ the projection map. Since $T^HM\cong\pi^*TS$, for $U\in TS$, let $U^H\in T^HM$ be the horizontal lift of $U$, such that $\pi_*U^H=U$. Let $T^H(\cdot,\cdot)\in \Lambda^2(TS)\otimes TX$ be the curvature of $(\pi,T^HM)$:
\begin{equation}\label{ba1.}
	T^H(U^H,V^H)=-P^{TX}[U^H,V^H],\ \  \text{for\ }U,V\in TS.
\end{equation}
We also have an identification of bundles through the horizontal lift
\begin{equation}\label{ba4}
	\pi^*\Lambda\lk T^*S\rk\widehat{\otimes}\Lambda\lk T^*X\rk\cong \Lambda\lk T^*M\rk.
\end{equation}

Let $g^{TX}$ and $g^{TS}$ be Riemannian metrics on $TX$ and $TS$ respectively. We equip $TM$ with the metric $g^{TM}=g^{TX}\oplus\pi^*g^{TS}$ and the corresponding Levi-Civita connection $\nabla^{TM}$. Let $\nabla^{TX}$ be the connection on $TX$ defined by in \cite[Definition 1.6]{bis86}:
\begin{equation}
\nabla^{TX}=P^{TX}\nabla^{TM}P^{TX},
\end{equation}
and denote its curvature by $R^{TX}$. Let $\nabla^{TS}$ be the Levi-Civita connection of $(TS,g^{TS})$. Let $\nabla^{TM,\oplus}$ be the connection on  $TM$ given by $\nabla^{TM,\oplus}=\pi^*\nabla^{TS}\oplus\nabla^{TX}$. Put
\begin{equation}
S(\cdot)=\nabla^{TM}_\cdot-\nabla^{TM,\oplus}_\cdot\in\Omega^1(M,\mathrm{End}(TM)).
\end{equation}

\subsection{A flat bundle and its odd forms}\label{Bb}

Let $(F,\nabla^F)$ be a complex flat vector bundle on $M$ with a Hermitian metric $g^F$. Following \cite[Definition 4.1]{bz92}, set
\begin{equation}\label{bb1}
	\omega\big(\nabla^F,g^F\big)=\big(g^F\big)^{-1}\nabla^Fg^F\in\Omega^1(M,\mathrm{End}(F)),
\end{equation}
and it takes values in self-adjoint elements of $\mathrm{End}(F)$. 

By \cite[Definitions 4.2, Proposition 4.3]{bz92}, we have the following unitary conneciton $\nabla^{F,u}$ on $F$ with its curvature:
\begin{equation}\label{bb3}
	\nabla^{F,u}=\nabla^F+\frac{1}{2}\omega\big(\nabla^{F},g^F\big),\ \  R^{F,u}=-\frac{1}{4}\omega\big(\nabla^{F},g^F\big)^2.
\end{equation}
For $x\in \mathbb{R}$, set
\begin{equation}\label{bb7}
	h\lk x\rk=x\exp\lk x^2\rk,
\end{equation}
Set $\varphi$ the endomorphism of $\Lambda^\bullet(T^*M)\otimes_\mathbb{R}\mathbb{C}$ sends $\alpha\in \Lambda^k(T^*M)\otimes_\mathbb{R}\mathbb{C}$ to $(2\pi i)^{-k/2}\alpha$. Put
	\begin{equation}\label{bb9}
		h\big(\nabla^F,g^F\big)=\lk2\pi i\rk^{1/2}\varphi\tro\big[h\big(\omega\big(\nabla^{F},g^F\big)/2\big)\big].
	\end{equation}

\begin{theo}[{\cite[Theorems 1.8, 1.9 and 1.11]{bl95}}]
	The odd form $h\lk\nabla^F,g^F\rk$ is real and closed and its cohomology class does not depend on $g^F$.
\end{theo}

\subsection{A superconnection and odd forms}\label{Be}

We make the same assumption as in \cref{Ba}. Let $\Omega^\bullet(X,F)$ be a formal smooth infinite-dimensional $\mathbb{Z}$-graded vector bundle over $S$ whose fibre over $b\in S$ is $C^\infty(X_{b},\Lambda^\bullet(T^*X)\otimes F|_{b})$. By \eqref{ba4} we have an identification of $\mathbb{Z}$-graded vector spaces
\begin{equation}\label{be1}
	\Omega^\bullet\lk M,F\rk\cong \Omega^\bullet(S,\Omega^\bullet(X,F)).
\end{equation}

The exterior differential operator $d^F$ acting on $\Omega^\bullet\lk M,F\rk$, then by \eqref{be1}, it can be considered as a flat superconnection on $\Omega^\bullet\lk X,F\rk$, we also denote it by $A'=d^F$. 

If $U\in TS$, the Lie derivative operator $L_{U^H}$ acts naturally on $\Omega^\bullet(X,F)$.  Let $\nabla^{\Omega^\bullet(X,F)}$ be the connection on $\Omega^\bullet\lk X,F\rk$ such that if $U\in TM$ and $s\in\Omega^\bullet\lk X,F\rk$,
\begin{equation}
	\nabla^{\Omega^\bullet(X,F)}_{U^H}s={L}_{U^H}s.
\end{equation}

We set ${i}_{T^H}=\frac{1}{2}f^\alpha\wedge f^\beta \wedge i_{T^H(f_\alpha,f_\beta)}$  (see \eqref{ba1.}), which acts naturally on $\Lambda^\bullet(T^*S)\hti\Lambda^\bullet(T^*X)\otimes F$. Let $d_X^{F}$ be the fibrewise exterior differential operator on $\Omega^\bullet(X,F)$.
\begin{prop}[{\cite[Proposition 3.4]{bl95}}]
	The superconnection $A'(=d^F)$ satisfies
	\begin{equation}\label{be3}
		A'=d_X^F+\nabla^{\Omega^\bullet(X,F)}+{i}_{T^H}.
	\end{equation}
\end{prop}

Let $g^{\Omega^\bullet(X, F)}$ be $L^2$-metric on $\Omega^\bullet(X, F)$ induced by $(g^{TX},g^F)$. Then let $d_X^{F,*}$ be the fibrewise formal adjoint operator of $d^F_X$. Let $\nabla^{\Omega^\bullet(X,F),*}$ be the adjoint connection of $\nabla^{\Omega^\bullet(X,F)}$. By identifying $TX$ and  $T^*X$ through $g^{TX}$, we can consider $T^H$ as a section of $\Lambda^2(T^*S)\hti T^*X$. Then $T^H\wedge$ acts naturally on $\Lambda^\bullet(T^*S)\hti\Lambda^\bullet(T^*X)\otimes F$.

Let $A''$ be the adjoint of $A'$ with respect to the metric $g^{\Omega^\bullet(X,F)}$ in the sense of \cite[\S\,1.4]{bl95}, which is also a flat superconnection. By \cite[Proposition 3.7]{bl95}, we have
\begin{equation}
	A''=d_X^{F,*}+\nabla^{\Omega^\bullet(X,F),*}-T^H\wedge.
\end{equation}
Set
\begin{equation}\label{be5}
	A=\frac{1}{2}\lk A'+A''\rk,\ \ \ B=\frac{1}{2}\lk A''-A'\rk,\ \ \ D_X^F=d_X^F+d_X^{F,*},
\end{equation}
then $A$ is a superconnection on $\Omega^\bullet(X,F)$, $B$ is a section of $(\pi^*\Lambda^\bullet \lk T^*S\rk\hti \mathrm{End}\lk\Omega^\bullet(X,F)\rk)^{\mathrm{odd}}$ and $D_X^F$ is the \emph{fiberwise Dirac operator}. 

Let $\varphi$ be the action on $\Lambda(T^*S)\otimes\mathbb{C}$ as in \eqref{bb9}.

\begin{defi}
For $B$ given in \eqref{be5}, we set the following odd form similar to \eqref{bb9}:
	\begin{equation}\label{bh1}
		h\blk A',g^{\Omega^\bullet(X,F)}\brk=\lk{2\pi i}\rk^{1/2}\varphi \trs\left[h(B)\right]\in\Omega^\bullet(S).
	\end{equation}
\end{defi}

\subsection{The analytic torsion forms}\label{Bi}

Following \cite[\S\,5.4]{bmz17}, for $t>0$, we set a rescaling metric $g_t^{TX}=g^{TX}/t$ and the associated metric $g^{\Omega^\bullet(X,F)}_t$ on $\Omega^\bullet(X,F)$.

For $a\in\mathbb{R}$, let $\psi_a$ be the automorphism of $\Lambda^\bullet\lk T^*S\rk$ such that, if $\alpha\in \Lambda^k(T^*S)$
then $\psi_a\alpha=a^k\alpha$. For $t>0$, put
\begin{equation}\label{be14}
	C_t=\psi_{\sqrt{t}}^{-1}\sqrt{t}A\psi_{\sqrt{t}}, \ \ \ \ D_t=\psi_{\sqrt{t}}^{-1}\sqrt{t}B\psi_{\sqrt{t}}.
\end{equation}

For $t>0$, let $h\big(A',g^{\Omega^\bullet(X,F)}_t\big)$ be the form as in \eqref{bh1} associated with $g^{\Omega^\bullet(X,F)}_t$. By \cite[(5.28)]{bmz17}, we have
\begin{equation}\label{bh2}
	\begin{split}
		h\big(A',g^{\Omega^\bullet(X,F)}_t\big)=\lk{2\pi i}\rk^{1/2}\varphi \trs\left[h(D_t)\right].
	\end{split}
\end{equation}
By \eqref{be5} and \eqref{be14}, we have $C_t^2=-D_t^2$, together with \eqref{bb7} and \eqref{bh2}, we get
\begin{equation}\label{bh2..}
	\begin{split}
		h\big(A',g^{\Omega^\bullet(X,F)}_t\big)=\lk{2\pi i}\rk^{1/2}\varphi \trs\Big[D_t\exp\big(-C_t^2\big)\Big].
	\end{split}
\end{equation}

Let $H^\bullet\lk X,F\rk=\oplus_{i=0}^{\dim X} H^i(X,F)$ be the $\mathbb{Z}$-graded vector bundle over $S$ whose fibre over $b\in S$ is the cohomology $H^\bullet(X_b,F|_{X_b})$ of the sheaf of locally flat sections of $F|_{X_b}$ on $X_b$. For the fiberwise Dirac operator $D_X^F$ in \eqref{be5}, by the Hodge Theorem, we have
\begin{equation}\label{be6}
H^i\lk X,F\rk\cong\ker D_X^{F,2}\big|_{\Omega^i(X,F)}\ \ \text{for any\ } i\in\mathbb{N}.
\end{equation}
By \cite[\S\,2a]{bl95}, $H^\bullet\lk X,F\rk$ is canonically equipped with a flat connection $\nabla^{H^\bullet(X,F)}$ and a Hermitian metric $g^{H^\bullet(X,F)}$. Let $e(TX,\nabla^{TX})\in\Omega^\bullet(M)$ be the Euler form of $TX$ (see \cite[(1.43)]{bmz17}) and we set $h(\nabla^{H^\bullet(X,F)},g^{H^\bullet(X,F)})=\sum_j(-1)^jh(\nabla^{H^j(X,F)},g^{H^j(X,F)})$.

\begin{theo}[{\cite[Theorem 3.16]{bl95}}]\label{Bh14}
	The forms $h\blk A',g^{\Omega^\bullet(X,F)}_t\brk$ are real, odd and closed, and their cohomology class does not depend on $t>0$. Moreover, as $t\rightarrow0$, we have
\begin{equation}
	h\blk A',g^{\Omega^\bullet(X,F)}_t\brk=\begin{cases}
	\int_Xe(TX,\nabla^{TX})h\big(\nabla^F,g^F\big)+\mathcal{O}(\sqrt{t}),\ \text{as\ } t\rightarrow0,\\
		h\big(\nabla^{H^\bullet(X,F)},g^{H^\bullet(X,F)}\big)+\mathcal{O}(1/\sqrt{t}),\ \text{as\ } t\rightarrow+\infty.
	\end{cases}
\end{equation}
\end{theo}
Now we review the transgression procedure following \cite[\S\,3.9]{bl95}. We enlarge $M, S$ to $\widetilde{M}=M\times \mathbb{R}^*_+, \widetilde{S}=S\times \mathbb{R}^*_+$. Set $\widetilde{\pi}\colon \widetilde{M}\rightarrow\widetilde{S}$ by $\widetilde{\pi}(x,s)=(\pi(x),s)$. Let $\rho_0\colon \widetilde{M}\rightarrow M$ and $\rho_1\colon \widetilde{M}\rightarrow \mathbb{R}^*_+$ be the projection maps. Let $\widetilde{X}$ be the fiber of $\widetilde{\pi}$, then we have $T\widetilde{X}=\rho^*TX$. and we equip $T\widetilde{X}$ with the metric $\rho^*g^{TX}/s$ over $M\times\{s\}$. We have
\begin{equation}\label{bf4}
	\begin{split}
		{\nabla}^{T\widetilde{X}}=\rho^*\nabla^{T{X}}+\de s \Big( \frac{\pa}{\pa s}-\frac{1}{2s}\Big),\ \ R^{T\widetilde{X}}=\rho_0^*R^{TX}.
	\end{split}
\end{equation}
On $\widetilde{S}$, we have the odd form $h\big(\widetilde{A}',\widetilde{g}^{\Omega(X,F)}_t\big)$ analogous to $h\big({A}',{g}^{\Omega(X,F)}_t\big)$ on $S$.
\begin{defi}
	Let $h^\wedge\big({A}',{g}^{\Omega(X,F)}_{t}\big)$ be the even form on $S$ satisfies (see \cite[(3.114)]{bl95})
	\begin{equation}\label{bh13}
		h\big(\widetilde{A}',\widetilde{g}^{\Omega(X,F)}_t\big)\big|_{s=1}=h\big({A}',{g}^{\Omega(X,F)}_{t}\big)+{\de s}\wedge h^\wedge\big({A}',{g}^{\Omega(X,F)}_{t}\big).
	\end{equation}
\end{defi}
Set  
\begin{equation}
	\begin{split}
\chi'(X,F)=\sum_{j=0}^m(-1)^jj\dim H^j(X,F),\ \ \chi(X,F)=\sum_{j=0}^m(-1)^j\dim H^j(X,F).
	\end{split}
\end{equation}
By Theorem \ref{Bh14} and \eqref{bh13} we have:
\begin{theo}[{\cite[Theorems 3.20, 3.21]{bl95}}]\label{Bh17}
	The form $h^\wedge\blk A',g^{\Omega^\bullet(X,F)}_t\brk$ is even, and
\begin{equation}
	\begin{split}
		&\frac{\pa}{\pa t}h\blk A',g^{\Omega^\bullet(X,F)}_t\brk=\frac{1}{t}h^\wedge\big( A',g^{\Omega^\bullet(X,F)}_t\big),\\
&h^\wedge\big( A',g^{\Omega^\bullet(X,F)}_t\big)=\begin{cases}
	\mathcal{O}(\sqrt{t}),\ \text{as\ } t\rightarrow0,\\
	\big(\frac{1}{2}\chi'(X,F)-\frac{m}{2}\chi(X,F)\big) h'(0)+\mathcal{O}(1/\sqrt{t}),\ \text{as\ } t\rightarrow+\infty.
\end{cases}
	\end{split}
\end{equation}
\end{theo}

Now we give the definition of the analytic torsion form $\mathcal{T}\big( T^HM,g^{TX},\nabla^{F},g^{F}\big)$.
\begin{defi}[see {\cite[Definition 3.22]{bl95}}]
	Set
	\begin{equation}\label{bi1}
		\begin{split}
			\mathcal{T}\big( T^HM,g^{TX},\nabla^{F},g^{F}\big)=&-\int_0^{+\infty}\Big[h^\wedge\blk A',g^{\Omega^\bullet(X,F)}_t\brk\\
			&+\Blk\frac{m}{4}\chi(X,F)-\frac{1}{2}\chi'(X,F)\Brk \blk h'(0)-h'(i\sqrt{t}/2)\brk\Big]\frac{d t}{t}.
		\end{split}
	\end{equation}
	By Theorems \ref{Bh14} and \ref{Bh17}, the above integral is well defined.
\end{defi}

\begin{theo}[{\cite[Theorem 3.23]{bl95}}]
	The form $\mathcal{T}\big( T^HM,g^{TX},\nabla^{F},g^{F}\big)$ is even and real, and it verifies the following transgression formula:
	\begin{equation}\label{bi2}
		\begin{split}
			\de \mathcal{T}\big( T^HM,g^{TX},\nabla^{F},g^{F}\big)=\int_Xe(TX,\nabla^{TX})h\big(\nabla^F,g^F\big)-h\big(\nabla^{H^\bullet(X,F)},g^{H^\bullet(X,F)}\big).
		\end{split}
	\end{equation}
\end{theo}

\subsection{Lichnerowicz type formulas}\label{Bj}

Let $z$ be an odd Grassmannian variable anticommutes with odd variables we used before. If $a\in \mathbb{R}[z]\hti\pi^*\Lambda^\bullet\lk T^*S\rk\hti\Lambda^\bullet\lk T^*X\rk$, we set
\begin{equation}\label{bj26}
[a]^z=c,\ \ \text{for\ } a=b+zc\ \text{where\ } b,c\in\pi^*\Lambda^\bullet\lk T^*S\rk\hti\Lambda^\bullet\lk T^*X\rk.
\end{equation}
By \eqref{bh2..}, we get
\begin{equation}\label{bj4}
	\begin{split}
		h\big(A',g^{\Omega^\bullet(X,F)}_t\big)=\lk{2\pi i}\rk^{1/2}\varphi \trs\Big[\exp\big(-\lk C_t^2-zD_t\rk\big)\Big]^z.
	\end{split}
\end{equation}

From now on, we will always use Latin indices $i,j,\cdots$ for vertical variables, and Greek indices $\alpha,\beta,\cdots$ for horizontal variables. Let $e_1,\cdots,e_m$ be a local orthonormal frame of $TX$, and let $f_1,\cdots,f_r$ be a basis of $TS$. The corresponding dual bases are denoted with upper indices. We set the following Clifford actions on $\Lambda(T^*X)$:
	\begin{equation}\label{aa.2}
		c({e}_j)={e}^j\wedge-{i}_{{e}_j},\ \ \ \ \ \widehat{c}({e}_j)={e}^j\wedge+{i}_{{e}_j},
	\end{equation}
then for $1\leqslant i,j\leqslant m$ we have
\begin{equation}\label{aa3}
	[c({e}_i),c({e}_j)]=-2\delta_{ij},\ \ [\widehat{c}({e}_i),\widehat{c}({e}_j)]=2\delta_{ij}, \ \  [{c}({e}_i),\widehat{c}({e}_j)]=0.
\end{equation}
By \eqref{aa3}, the actions in \eqref{aa.2} extend to an isomorphism of algebras
\begin{equation}\label{aa5}
	c\hti\widehat{c}\colon{c}({TX})\hti\widehat{c}({TX})\rightarrow \mathrm{End}(\Lambda({T^*X})).
\end{equation}

Put
\begin{equation}
	\mathcal{R}^F=-\frac{1}{4}\big\langle R^{TX}e_i,e_j \big\rangle\widehat{c}\lk e_i\rk\widehat{c}\lk e_j\rk-\frac{1}{4}\omega^2\big(\nabla^{F},g^{F}\big).
\end{equation}
Let $r^X$ be the scalar curvature of the fibre $\big(X,g^{TX}\big)$.

\begin{defi}
	Let $\Lambda^F$	 be a section of $\mathbb{R}[z]\hti\Lambda^\bullet(T^*S)\hti\mathrm{End}(\Lambda({T^*X}))\otimes \mathrm{End}(F)$ by
	\begin{equation}\label{bj7}
		\begin{split}
			\Lambda^{F}=&\frac{1}{4}r^X+\frac{1}{2}c(e_i)c(e_j)\mathcal{R}\lk e_i,e_j\rk+\frac{1}{2}f_\alpha f_\beta\mathcal{R}\big(f_\alpha^H,f_\beta^H\big)+c(e_i) f_\alpha\mathcal{R}\big(e_i,f_\alpha^H\big)\\
			&+\frac{1}{4}{\omega}\big(\nabla^{F},g^{F}\big)(e_i)^2+\frac{1}{8}\widehat{c}(e_i)\widehat{c}(e_j){\omega}\big(\nabla^{F},g^{F}\big)^2\lk e_i,e_j\rk\\
			&-\frac{1}{2}f^H_\alpha\widehat{c}(e_i)\nabla^{TX\otimes F_p,u}_{f_\alpha^H}\omega\big(\nabla^{F},g^{F}\big)\lk e_i\rk-\frac{1}{2}c(e_i)\widehat{c}(e_j)\nabla^{TX\otimes F_p,u}_{e_i}\omega\big(\nabla^{F},g^{F}\big)(e_j)\\
			&-\frac{1}{2}zc(e_i)\omega\big(\nabla^{F},g^{F}\big)(e_i)-\frac{1}{2}zf_\alpha\omega\big(\nabla^{F},g^{F}\big)\big(f^H_\alpha\big).
		\end{split}
	\end{equation}
\end{defi}

\begin{defi}
	Let ${}^{0}\nabla^{\mathbb{R}[z]\hti\pi^*\Lambda^\bullet\lk T^*S\rk\hti\Lambda^\bullet\lk T^*X\rk}$ be the fibrewise connection along $X$ by:
	\begin{equation}\label{bi.33}
		\begin{split}
			{}^{0}\nabla^{\mathbb{R}[z]\hti\pi^*\Lambda^\bullet\lk T^*S\rk\hti\Lambda^\bullet\lk T^*X\rk}=&\nabla^{\pi^*\Lambda^\bullet\lk T^*S\rk\hti\Lambda^\bullet\lk T^*X\rk}+\frac{1}{2}\big\langle Se_i,f_\alpha^H\big\rangle c\lk e_i\rk f^\alpha\\
			&+\frac{1}{4}\big\langle Sf_\alpha^H,f_\beta^H\big\rangle f^\alpha f^\beta-\frac{z}{2}e^i\wedge\widehat{c}(e_i).
		\end{split}
	\end{equation}
\end{defi}
By \cite[Theorem 4.14]{bis86}, the curvature of ${}^{0}\nabla^{\mathbb{R}[z]\hti\pi^*\Lambda^\bullet\lk T^*S\rk\hti\Lambda^\bullet\lk T^*X\rk}$ is given by
\begin{equation}\label{bj.6}
	\begin{split}
		{}^{0}R^{\mathbb{R}[z]\hti\pi^*\Lambda^\bullet\lk T^*S\rk\hti\Lambda^\bullet\lk T^*X\rk}(e_i,e_j)=\frac{1}{2}\big\langle R^{TX}(e_k,e_\ell)e_i,e_j\big\rangle \big(c(e_k)c(e_\ell)-\widehat{c}(e_k)\widehat{c}(e_\ell)\big)&\\
		+ \big\langle R^{TX}(f_\alpha,e_\ell)e_i,e_j\big\rangle f^\alpha c(e_\ell)+\frac{1}{2}\big\langle R^{TX}(f_\alpha,f_\beta)e_i,e_j\big\rangle f^\alpha\wedge f^\beta.&
	\end{split}
\end{equation}

 On $\mathbb{R}[z]\hti\pi^*\Lambda^\bullet\lk T^*S\rk\hti\Lambda^\bullet\lk T^*X\rk\otimes F$, let ${}^0\nabla^{\mathbb{R}[z]\hti\pi^*\Lambda^\bullet\lk T^*S\rk\hti\Lambda^\bullet\lk T^*X\rk\otimes F,u}$ be the fibrewise connection induced by ${}^0\nabla^{\mathbb{R}[z]\hti\pi^*\Lambda^\bullet\lk T^*S\rk\hti\Lambda^\bullet\lk T^*X\rk}$ and $\nabla^{F,u}$.

Put $\mathcal{L}^F_t=C_t^2-zD_t$, which appears in \eqref{bj4} and  we set $\mathcal{L}^F=\mathcal{L}^F_4$. Let ${N}_1$ be the number operator of the exterior algebra $\mathbb{R}[z]\hti \pi^*\Lambda^\bullet(T^*S)$ that acts by multiplication by the total degree of each term. By \eqref{be14}, we have
\begin{equation}\label{bj13}
	\mathcal{L}^F_{4t}=\theta_{\sqrt{t}}^{-1}t\mathcal{L}^F\theta_{\sqrt{t}},\ \ \text{where\ } {\theta}_{\sqrt{t}}={\sqrt{t}}^{{N}_1}.
\end{equation}

\begin{theo}[{\cite[Theorem 3.11]{bl95}}]\label{2.19}
	The following Lichnerowicz type formula holds:
	\begin{equation}\label{bj14}
		\begin{split}
			\mathcal{L}^F=-\big({}^0\nabla^{\mathbb{R}[z]\hti\pi^*\Lambda^\bullet\lk T^*S\rk\hti\Lambda^\bullet\lk T^*X\rk\otimes F,u}_{e_i}\big)^2+{}^0\nabla^{\mathbb{R}[z]\hti\pi^*\Lambda^\bullet\lk T^*S\rk\hti\Lambda^\bullet\lk T^*X\rk\otimes F,u}_{\nabla^{TX}_{e_i}e_i}+\Lambda^{F}.
		\end{split}
	\end{equation}
\end{theo}

\section{Toeplitz operators}\label{sC}

In this section, we describe the formalism of the Toeplitz operator introduced by Berezin \cite{bere} and Boutet de Monvel-Guillemin \cite{mo81}, and developed by Bordemann-Meinrenken-Schlichenmaier \cite{bord}, Schlichenmaier \cite{sch} and Ma-Marinescu \cite{mm07}, \cite{mm08}, \cite{mm12}. 

This section is organized as follows. In \cref{Ca}, we introduce the positive line bundle $L$ on a compact Kähler manifold $N$ and the Berezin-Toeplitz quantization, then we give some uniform estimations that will be used later. In \cref{Cc}, we investigate the exponential of a Toeplitz operator and give an estimation for each term in the expansion and the remainder.

\subsection{The algebras of Toeplitz operators}\label{Ca}

Let $N$ be a compact complex manifold of $\mathrm{dim}_{\mathbb{C}}N=n$ with the complex structure $J$. Let $g^{T_{\mathbb{R}}N}$ be a $J$-invariant  Riemannian metric on $T_{\mathbb{R}}N$. Denote the induced Riemannian volume form by $d{v}_N$.

Let $(L,g^L)$ (resp. $(\xi,h^\xi)$) be a holomorphic line bundle (resp. holomorphic vector bundle) on $N$. Let $\nabla^L$ (resp. $\nabla^{{{\xi}}}$) be the associated Chern connections on $L$ (resp. ${{{\xi}}}$) with curvature $R^L$ (resp. $R^{{{\xi}}}$). We assume that $(L,g^L,\nabla^L)$ is a positive line bundle. Then $c_1(L,g^L)=\frac{\sqrt{-1}}{2\pi}R^L$ defines a Kähler form on $N$. We note that $g^{T_\mathbb{R}N}$ is \emph{not} necessarily given by the Kähler metric induced by $c_1(L,g^L)$.

For ${p\in{\mathbb{N}^*}}$, put $L^p=L^{\otimes p}$, the $p$-th tensor power of $L$. We have the $L^2$-inner product on $\mathscr{C}^\infty(N,L^p\otimes {{{\xi}}})$ induced by $(g^{T_\mathbb{R}N},g^L,h^{{{\xi}}})$. We denote the corresponding norm by $\lV\cdot\rV_{L^2}$ and by $L^2(N,L^p\otimes {{{\xi}}})$ the completion of $\mathscr{C}^\infty(N,L^p\otimes {{{\xi}}})$ with respect to the $L^2$-norm.

We consider the space $H^{(0,0)}(N,L^p\otimes {{{\xi}}})$ of holomorphic sections of $L^p\otimes\xi$, and let 
\begin{equation}
	P_p\colon L^2(N,L^p\otimes {{{\xi}}})\longrightarrow H^{(0,0)}(N,L^p\otimes {\xi})
\end{equation}
be the orthogonal projection with respect to the $L^2$-product.

\begin{defi}[{\cite[(7.2.6)]{mm07}}]\label{Cb2}
	The Berezin-Toeplitz quantization of a smooth section $\mathcal{H}\in \mathscr{C}^\infty(N,\mathrm{End}({{{\xi}}}))$ is a sequence of linear operators $\{T_{\mathcal{H},p}\}_{{p\in{\mathbb{N}^*}}}$ given by
	\begin{equation}
		T_{\mathcal{H},p}\colon L^2(N,L^p\otimes {{{\xi}}})\rightarrow L^2(N,L^p\otimes {{{\xi}}}),\ \ \ T_{\mathcal{H},p}=P_p\mathcal{H}P_p.
	\end{equation}
\end{defi}
The operator $T_{\mathcal{H},p}$ has a smooth kernel $T_{\mathcal{H},p}(x,x')$ with respect to $dv_N$. On the diagonal, $T_{\mathcal{H},p}(x,x)\in \mathrm{End}(L^p\otimes {{{\xi}}})=\mathrm{End}({{{\xi}}}
)$. For $\ell\in\mathbb{N}$, let $\lv\cdot\rv_{\mathscr{C}^\ell(N)}$ be the $\ell$-th order smooth norm on $\mathscr{C}^\infty(N,\mathrm{End}({{\xi}}))$ induced by $(\nabla^{T_\mathbb{R}N},\nabla^E,g^{T_\mathbb{R}N},h^{{{\xi}}})$. By the proof of \cite[Lemma 7.2.4]{mm07}, we have the following expansion with a uniform estimation for the remainder.
\begin{theo}\label{Cb3}
	There is a series of differential operators $\{A_i\}_{i=0}^{+\infty}$ of order no more than $2i$ such that for any $k,\ell\in\mathbb{N}$, there exists $C>0$ such that for any ${p\in{\mathbb{N}^*}}$ and $\h\in\mathscr{C}^\infty(N,\mathrm{End}({{{\xi}}}))$, 
\begin{equation}\label{cb6}
	\Bv p^{-n}T_{\mathcal{H},p}(x,x)-\sum_{i=0}^kp^{-i}A_{i}(\mathcal{H})(x)\Bv_{\mathscr{C}^\ell(N)}\leqslant C\lv\mathcal{H}\rv_{\mathscr{C}^{2k+2}(N)}p^{-k-1}.
\end{equation}
And the operators $\{A_i\}_{i=0}^{+\infty}$ vary smoothly with respect to $(J,g^{T_\mathbb{R}N},h^{L},h^{{{\xi}}})$. In particular,
\begin{equation}\label{cb4..}
A_{0}(\mathcal{H})dv_N=\frac{c_1(L,g^L)^n}{n!}.
\end{equation}
\end{theo}
	By the Kodaira vanishing Theorem \cite[Theorem 1.5.4]{mm07} and the Riemann-Roch-Hirzebruch Theorem \cite[Theorem 1.4.6]{mm07}, for $p\in\mathbb{N}^*$ large enough, $\dim H^{(0,0)}(N,L^p\otimes {{{\xi}}})$ is a polynomial of $p\in\mathbb{N}^*$ with leading term
	\begin{equation}\label{cb9}
		\begin{split}
			\dim H^{(0,0)}(N,L^p\otimes {{{\xi}}})=\dim_\mathbb{C}\xi\cdot\int_N\frac{c_1(L,g^L)^n}{n!}p^n+\mathcal{O}\big(p^{n-1}\big),
		\end{split}
	\end{equation}
and clearly  \eqref{cb4..} is a local version of \eqref{cb9}.

If $A\in\mathrm{End}\lk L^2\lk N,L^p\otimes {{{\xi}}}\rk\rk$, let $\lV A\rV$ be its operator norm. If $A$ is of trace class, we denote its trace norm by $\lV A\rV_1$. If $A=P_pAP_p$, then $A$ is of trace class, by \eqref{cb9}, there is $C>0$ such that for any ${p\in{\mathbb{N}^*}}$,
\begin{equation}\label{cb12}
	\lV A\rV_1\leqslant \lV A\rV\dim_\mathbb{C} H^{(0,0)}(N,L^p\otimes {{{\xi}}})\leqslant C\lV A\rV p^n.
\end{equation}

Following Ma-Marinescu \cite[Definition 7.2.1]{mm07}, we now define the Toeplitz operator.
\begin{defi}
	A Toeplitz operator is a family of operators  $\big\{T_p\mid T_p\in \mathrm{End}(L^2\lk N,L^p\otimes {{{\xi}}}\rk)\big\}_{{p\in{\mathbb{N}^*}}}$ bounded with respect to the $L^2$-norm such that 
$T_p=P_pT_pP_p$, and that there exists $\{\mathcal{H}_i\mid \mathcal{H}_i\in \mathscr{C}^\infty(N,\mathrm{End}({{{\xi}}}))\}_{i\in \mathbb{N}}$ such that for any $k\in\mathbb{N}$, there is $\bm{c}_k>0$ such that
	\begin{equation}\label{cb16}
		\BV T_p-\sum_{i=0}^k p^{-i}T_{\mathcal{H}_i,p}\BV\leqslant \bm{c}_kp^{-k-1}.
	\end{equation}
	As in \cite[(7.2.4), (7.2.5)]{mm07}, we use the notation of formal expansion to denote \eqref{cb16} as
	\begin{equation}\label{cb17}
		T_p=\sum_{i=0}^{+\infty} p^{-i}T_{\mathcal{H}_i,p}+\mathcal{O}\lk p^{-\infty}\rk,
	\end{equation}
	and we replace $\sum_{i=0}^{+\infty}$ with $\sum_{i=0}^{k}$ if we only refer to the first $k$ terms.
\end{defi}

Now we introduce the asymptotic trace symbol for ease of notation. For the Toeplitz operator $\{T_p\}_{{p\in{\mathbb{N}^*}}}$ in \eqref{cb16}, for $k\in\mathbb{N}$, if we set 
\begin{equation}\label{cb19}
	\begin{split}
		\tro_{[k]}[T_p]=\sum_{i+j=k}\int_N\tro^{\xi}[A_i(\mathcal{H}_j)]\de v_N,
	\end{split}
\end{equation}
then by \eqref{cb6}, \eqref{cb12} and  \eqref{cb16}, there is $C_{\bm{c}_k,\{\lv\mathcal{H}_{i}\rv_{\mathscr{C}^{2i}(N)}\}_{i=1}^{k+1}}>0$ which depends on $\bm{c}_k$ and $\{\lv\mathcal{H}_{i}\rv_{\mathscr{C}^{2i}(N)}\}_{i=1}^{k+1}$ such that 
\begin{equation}\label{cb20}
	\begin{split}
		\Bv p^{-n}\tro^{H^{(0,0)}(N,L^p\otimes {{{\xi}}})}[T_p]-\sum_{i=0}^kp^{-i}\tro_{[i]}[T_p]\Bv\leqslant C_{\bm{c}_k,\{\lv\mathcal{H}_{i}\rv_{\mathscr{C}^{2i}(N)}\}_{i=1}^{k+1}}p^{-k-1}.
	\end{split}
\end{equation}

By the proof of \cite[Theorem 7.4.1]{mm07} and \cite[Theorem 0.3]{mm12}, we have the following product formula of Toeplitz operators with estimation for the remainder.
\begin{theo}\label{Cb7}
	The set of Toeplitz operators is an algebra. In particular, for any $k\in\mathbb{N}$, there is $C>0$ such that for any $\mathcal{H},\mathcal{H}'\in\mathscr{C}^\infty(N,\mathrm{End}({{{\xi}}}))$ and ${p\in{\mathbb{N}^*}}$, we have
\begin{equation}\label{cb24}
	\bV T_{\mathcal{H},p}T_{\mathcal{H}',p}-\sum_{j=0}^{k} p^{-i}T_{C_j\lk\mathcal{H},\mathcal{H}'\rk,p}\bV\leqslant Cp^{-k-1}\lv\mathcal{H}\rv_{\mathscr{C}^{2k+2}(N)}\cdot\lv\mathcal{H}'\rv_{\mathscr{C}^{2k+2}(N)}.
\end{equation}
where $C_j(\cdot,\cdot)$ is a smooth bidifferential operator of total degree no more than $2j$, and
	\begin{equation}\label{cb23}
		\begin{split}
			C_0(\mathcal{H},\mathcal{H}')=\mathcal{H}\mathcal{H}',\ \ \ C_1(\mathcal{H},\mathcal{H}')-C_1(\mathcal{H}',\mathcal{H})=i\{\mathcal{H},\mathcal{H}'\}\  \text{if\ } \mathcal{H},\mathcal{H}'\in\mathscr{C}^\infty(N,\mathbb{C}),
		\end{split}
	\end{equation}
where $\{\mathcal{H},\mathcal{H}'\}$ is the Poisson bracket of $2\pi c_1(L,g^L)$.
\end{theo}

\subsection{The exponential of Toeplitz operators}\label{Cc}

To discuss the exponential of Toeplitz operators, we first study the inverse of a Toeplitz operator. We follow Ma-Zhang \cite[Lemma 4.1]{mz08} where they have proved that if the principal symbol $\h_0$ of a $T_p$ as in \eqref{cb16} is invertible, then $T_p^{-1}$ is also a Toeplitz operator, here we give a uniform estimation for the remainder term in the expansion.
\begin{lemma}\label{Cc9}
Let $T_p$ be a Toeplitz operator as in \eqref{cb16}. If $\h_0(x)$ is invertible for all $x\in N$, then $T_p$ is invertible for $p$ large, and $T_p^{-1}$ is also a Toeplitz operator. Moreover, if we write $T_p^{-1}=\sum_{i=0}^{\infty}p^{-i}T_{\mathcal{G}_i,p}+\mathcal{O}(p^{-\infty})$ as in \eqref{cb17}, then we have $\mathcal{G}_0=\mathcal{H}_0^{-1}$.
\end{lemma}

\begin{pro}
We recall a basic fact: for invertible operators $A$ and $B$, we have
 \begin{equation}\label{cc.31}
 	\begin{split}
 		A^{-1}-B=B((1-(1-AB))^{-1}-1)=B(1-AB)(1-(1-AB))^{-1}.
 	\end{split}
 \end{equation}

For $k\in\mathbb{N}$ and $C_i(\cdot,\cdot)$ given in Theorem \ref{Cb7}, we set
	\begin{equation}\label{cc11}
		\begin{split}
\mathcal{G}_{0}=\mathcal{H}_0^{-1},\ \ \mathcal{G}_{k+1}=-\mathcal{H}_0^{-1}\sum_{\substack{i+j\leqslant k+1\\(i,j)\neq(0,0)}} C_i\blk\mathcal{H}_j,\mathcal{G}_{k+1-i-j}\brk.
\end{split}
	\end{equation}
By \eqref{cb16}, \eqref{cb24} and \eqref{cc11}, we can prove inductively that
	\begin{equation}\label{cb14}
		\bbV T_p\Blk\sum_{i=0}^{k}p^{-i}T_{\mathcal{G}_i,p}\Brk-1\bbV\leqslant \bm{C}p^{-k-1},\ \ \bbV \Blk\sum_{i=0}^{k}p^{-i}T_{\mathcal{G}_i,p}\Brk T_p-1\bbV\leqslant \bm{C}p^{-k-1}.
	\end{equation}
By \eqref{cb14}, $T_p$ and $\sum_{i=0}^{k}p^{-i}T_{\mathcal{G}_i,p}$ have both left and right inverse for $p\in\mathbb{N}^*$ large, hence they are invertible. We take $A=T_p,B=\sum_{i=0}^{k}p^{-i}T_{\mathcal{G}_i,p}$ in \eqref{cc.31}, by \eqref{cb14} and the obvious inequality $\lV(1-(1-AB))^{-1}\rV\leqslant (1-\lV(1-AB)\rV)^{-1}$, we get
\begin{equation}\label{cb17..}
\BV T_p^{-1}-\sum_{i=0}^kp^{-i}T_{\mathcal{G}_i,p}\BV\leqslant \bm{C}p^{-{k-1}},
\end{equation}
which completes the proof.\qed
\end{pro}

\begin{remark}
By the above proof, we can indeed take $\bm{C}$ in \eqref{cb14} and \eqref{cb17..} as the sum of terms of the following form:
\begin{equation}\label{cc..33}
	C\vert\mathcal{H}_0^{-1}\vert_{\mathscr{C}^{0}(N)}^{\alpha}\cdot\prod_{i=0}^{k}\bm{c_i}^{\beta_i}\vert\mathcal{H}_i\vert_{\mathscr{C}^{2(k+1)}(N)}^{\gamma_i}
\end{equation}
where $C>0$, $\alpha\in\mathbb{N}$ and $\beta,\gamma\in \mathbb{N}^{k+1}$.
\end{remark}

\begin{theo}\label{Cc10}
	If $T_p$ is a Toeplitz operator with expansion \eqref{cb16}, $\exp(-tT_p)$ is also a Toeplitz operator for any $t\geqslant 0$ with the following expansion in the sense of \eqref{cb17}:
	\begin{equation}\label{cc18}
		\exp(-tT_p)=\sum_{i=0}^{+\infty}p^{-i}T_{\mathcal{J}_i(t),p}+\mathcal{O}(p^{-\infty}),\ \ \ \mathcal{J}_0(t)=e^{-t\mathcal{H}_0}.
	\end{equation}
Put $\bm{h}=\inf_{x\in N}\mathrm{Re}\big(\mathrm{Spec}\h_0(x)\big)$, then for any $\varepsilon>0,i,k,\ell\in\mathbb{N}$, there is $C>0$ such that for any $p\in\mathbb{N}^*$,
\begin{equation}\label{cc.22}
	\begin{split}
		&\lv\mathcal{J}_i(t)\rv_{\mathscr{C}^k(N)}\leqslant  C\exp(-(\bm{h}-\varepsilon)t),\\
		&\BV\exp(-tT_p)-\sum_{i=0}^kp^{-i}T_{\mathcal{J}_i(t),p}\BV\leqslant C\exp(-(\bm{h}-\varepsilon)t)p^{-k-1}.
	\end{split}
\end{equation}
\end{theo}

\begin{pro}
 By the product formula \eqref{cb24} and \eqref{cb23}, for $k\in\mathbb{N}$, we have $T_{p}^k=T_{\mathcal{H}_0^k,p}+\mathcal{O}(p^{-1})$. If we ignore the divergence of infinite sums of remainders, we get $\exp({T_{p}})=\sum_{k=0}^{+\infty}\frac{1}{k!}T_{\mathcal{H}_0^k,p}+\mathcal{O}(p^{-1})=T_{\exp(\mathcal{H}_0),p}+\mathcal{O}(p^{-1})$, which is just the first order of \eqref{cc18}.
	
	Now we give rigorous proof.  
	By Lemma \ref{Cc9}, if $\lambda\notin \mathrm{Spec}(\h_0)$, then $\lambda-T_p$ is invertible for $p\in\mathbb{N}^*$ large and we have the following expansion in the sense of \eqref{cb17}:
	\begin{equation}\label{cc30}
		\lk \lambda-T_p\rk^{-1}=\sum_{i=0}^{+\infty}p^{-i}T_{\mathcal{I}_i(\lambda),p}
	\end{equation}
	where $\mathcal{I}_i(\lambda)\in \mathscr{C}^\infty(N,\mathrm{End}({{{\xi}}}))$. As in Theorem \ref{Cb7}, $C_i(\cdot,\cdot)$ is a bilinear differential operator for each $i\in \mathbb{N}$, then by \eqref{cc11}, each $\mathcal{I}_i(\lambda)$ for $i\in\mathbb{N}$ is the sum of terms of the following form:
		\begin{equation}\label{cc31}
		\lk\lambda-\mathcal{H}_0\rk^{-n_0}f_1(\lambda)\lk\lambda-\mathcal{H}_0\rk^{-n_1}\cdots f_j(\lambda)\lk\lambda-\mathcal{H}_0\rk^{-n_j}
	\end{equation}
	where $n_j\in \mathbb{N}$ and $f_j(\lambda)\in \mathscr{C}^\infty(N,\mathrm{End}({{{\xi}}}))[\lambda]$, the space of polynomials of $\lambda$ with coefficients in $\mathscr{C}^\infty(N,\mathrm{End}({{{\xi}}}))$. In particular, the leading term is
	\begin{equation}\label{cc32}
		\mathcal{I}_0(\lambda)=\lk \lambda-\mathcal{H}_0\rk^{-1}.
	\end{equation}

For $\varepsilon>0$, we choose a bounded curve $\Gamma_\varepsilon$ surrounds $\mathrm{Spec}(\h_0)$ counterclockwise on the complex plane such that $\inf_{\lambda\in\Gamma_\varepsilon} \mathrm{Re}(\lambda)\geqslant \bm{h}-\epsilon$. By the Cauchy integral formula, we have
\begin{equation}\label{cc29}
	\exp(-tT_p)=\frac{1}{{2\pi i}}\int_{\Gamma_\varepsilon} {e}^{-t\lambda}\lk \lambda-T_p\rk^{-1}\de \lambda.
\end{equation}
We denote the $\bm{C}$ in \eqref{cb17..} associated with Toeplitz operator $(\blk\lambda-T_p\brk)^{-1}$ by $\bm{C}_{\lambda}$. By \eqref{cc..33}, $\bm{C}_{\lambda}$ is uniformly bounded for $\lambda\in \Gamma$. By \eqref{cc30} and \eqref{cc29}, we have
	\begin{equation}\label{cc33}
		\bbV\exp(-tT_p)-\sum_{i=0}^kp^{-i}T_{\frac{1}{{2\pi i}}\int_{\Gamma_\varepsilon} {e}^{-t\lambda}\mathcal{I}_i(\lambda)\de\lambda,p}\bbV \leqslant p^{-k-1}\frac{1}{2\pi}\int_{\Gamma_\varepsilon} {e}^{-t\mathrm{Re}(\lambda)}\bm{C}_{\lambda}\de \lambda,
	\end{equation}
from which we get \eqref{cc18} and the first inequality of \eqref{cc.22} with
	\begin{equation}\label{cc34}
		\mathcal{J}_i(t)=\frac{1}{{2\pi i}}\int_{\Gamma_\varepsilon} {e}^{-t\lambda}\mathcal{I}_i(\lambda)\de \lambda.
	\end{equation}
And by \eqref{cc32}, we get $\mathcal{J}_0(t)$ in \eqref{cc18}. From \eqref{cc31} and \eqref{cc34}, we obtain the second inequality of \eqref{cc.22}.\qed
\end{pro}

\begin{remark}\label{Cc11}
In this study, estimations for each term in \eqref{cc.22} and the remainder term are essential. By \eqref{cb19} and \eqref{cc.22}, for any $t\geqslant0$, we have
	\begin{equation}\label{cc41}
		\Bv\tro_{[k]}\big[\exp(-tT_p)\big]\Bv\leqslant C{e}^{-{(\bm{h}-\epsilon)}t}.
	\end{equation}
	Moreover, by  \eqref{cb20}, \eqref{cc18} and \eqref{cc.22}, we get
	\begin{equation}\label{cc42}
		\Bv p^{-n}\tro^{H^{(0,0)}(N,L^p\otimes {{{\xi}}})}[\exp(-tT_p)]-\sum_{i=0}^kp^{-i}\tro_{[i]}[\exp(-tT_p)]\Bv\leqslant Cp^{-k-1}{e}^{-(\bm{h}-\epsilon)t}.
	\end{equation}
Equations \eqref{cc41} and \eqref{cc42} ensure the exponential decay of each term in the expansion and the remainder as $t\rightarrow +\infty$, and these play important roles in \Cref{sF} and \Cref{sG}. 
\end{remark}

\section{The geometry of bundles $\{F_p\}_{p\in\mathbb{N}^*}$}\label{sD}

We make the same assumption as in \cref{Ba,Ca}. The purpose of this section is to review some properties of bundles $\{F_p\}_{p\in\mathbb{N}^*}$ over $M$ following \cite[\S\,9]{bmz17}.

This section is organized as follows. In \cref{Da}, we present the geometry of the fibration $\mathcal{N}=P_G\times_GN\rightarrow M$ and the line bundle $L\to\mathcal{N}$. In \cref{Db}, we introduce the bundles $\{F_p\}_{p\in\mathbb{N}^*}$ over $M$ and compute the curvature of the unitary connection $\nabla^{F_p,u}$. In \cref{Dc},  for a smooth family of Toeplitz operators $\{T_p\}_{{p\in{\mathbb{N}^*}}}$, we analyze the asymptotics of $\tro^{F_p}[T_p]$ when $p\rightarrow+\infty$. In \cref{Dd},  we obtain the asymptotics of the odd form $h(\nabla^{F_p},g^{F_p})$ for $p\rightarrow+\infty$. In \cref{De},  we recall the nondegenarated condition of $L$.

\subsection{Geometric settings}\label{Da}

Let $N$ be a compact complex manifold of complex dimension $n$. Let $L$ and $\xi$ be holomorphic bundles on $N$ with $\dim_\mathbb{C}(L)=1$. Let $G$ be a Lie group acting holomorphically on $N$, and this action lifts to holomorphic actions on $L$ and $\xi$.

Let $p\colon P_G\rightarrow M$ be a principal \emph{flat} $G$-bundle. Set
\begin{equation}
	\mathcal{N}=P_G\times_GN.
\end{equation}
We denote by $q$ the projection $q\colon\mathcal{N}\rightarrow M$ with fibre $N$. We still denote the bundle $P_G\times_GL$ (resp. $P_G\times_G\xi$) over $\mathcal{N}$ by $L$ (resp. $\xi$).

Let $T_0^H\mathcal{N}\subset T\mathcal{N}$ be the horizontal bundle determined by the flat connection of $P_G$.
Put $T_\mathbb{R}N=\ker q_*$, the real relative tangent bundle, and let $TN$ be the \emph{holomorphic} relative tangent bundle. Let $J^N$ be the complex structure on $T_\mathbb{R}N$. We clearly have $T_0^H\mathcal{N}\cong q^*TM$, then for $U\in TM$, we denote by $U_0^H\in T_0^H\mathcal{N}$ its horizontal lift.

Let $g^L$ (resp. $g^{{{\xi}}}$) be a Hermitian metric on $L$ (resp. ${\xi}$) over $\mathcal{N}$. Let $\nabla^L$ (resp. $\nabla^{{{\xi}}}$) be the fibrewise Chern on $L$ (resp. ${{{\xi}}}$). Using the flat connection on $P_G$, we can extend $\nabla^L$ and $\nabla^{{{\xi}}}$ to  connections on $L$ and $\xi$ respectively: for $U\in TM$, put 
\begin{equation}
\nabla^{L}_{U^H_0}\ (\text{or\ }\nabla^{\xi}_{U^H_0})=L_{U^H_0}
\end{equation}
These connections are in general non-unitary, and we set
\begin{align}\label{da3}
	\omega\big(L,g^L\big)(U)=\big(g^L\big)^{-1}L_{U^H_0}g^L,\ \ \omega\big({{{\xi}}},g^{{{\xi}}}\big)(U)=\big(g^{{{\xi}}}\big)^{-1}L_{U^H_0}g^{{{\xi}}}.
\end{align}
In what follows, we assume that the first Chern form
\begin{equation}
c_1(L,g^L)=\frac{\sqrt{-1}}{2\pi }\big(\nabla^{L}\big)^2\in\Omega^2(\mathcal{N})
\end{equation}
is \emph{fibrewisely positive}. Let $\partial^N$ denote the fibrewise Dolbeault operators. We have the following identities (see \cite[(9.7)]{bmz17}): for $U,V\in TM$ and $Y\in T_\mathbb{R}N$,
\begin{equation}\label{da7}
c_1(L,g^L)(U^H_0,V^H_0)=0,\ \  c_1(L,g^L)(U^H_0,Y)=\frac{\sqrt{-1}}{2\pi}\partial^N_Y	\omega\big(L,g^L\big)(U).
\end{equation}

Let $g^{T_\mathbb{R}N}$ be a smooth Hermitian metric on $T_\mathbb{R}N$, and let $dv_N$ be the corresponding fibrewise volume form. Then we extend $dv_N$ to be a form on $\mathcal{N}$ that vanishes along the horizontal direction, and still denote it by $dv_N\in\Omega^\bullet(\mathcal{N})$. Recall that, if $U\in TM$, the fibrewise Lie derivative operator ${L}_{U^H_0}$ acts naturally on smooth sections of $\Lambda^\bullet(T^*_\mathbb{R}\mathcal{N})$. We define $\mathrm{div}_N(U)$ by
\begin{equation}\label{da5}
	{L}_{U^H_0}dv_N=\mathrm{div}_N(U)dv_N.
\end{equation}

By \eqref{da7}, when $g^{T_\mathbb{R}N}$ is just the Kähler metric $g^{T_\mathbb{R}N}(\cdot,\cdot)=c_1(L,g^L)(\cdot,J^N\cdot)$, then for any $U\in TM$, we can prove the following formula:
	\begin{equation}\label{da9}
		\mathrm{div}_N(U)=\frac{1}{4\pi}{\Delta}_N\omega\big(L,g^L\big)(U),
	\end{equation}
	where ${\Delta}_N$ is the fibrewise (nonnegative) Laplacian operator with respect to the Kähler metric which acts on $\mathscr{C}^\infty(\mathcal{N})$.

\subsection{The curvature of $\nabla^{F_p,u}$}\label{Db}

Recall the definition of $F_p$ in \eqref{0.4...}:
\begin{equation}
F_p=\big(P_G\times_GH^{(0,0)}(N,L^p\otimes {{\xi}})\big)\otimes F,
\end{equation}
where the flat bundle $(F,\nabla^F)$ with metric $g^F$ plays the role of a shifting. Let $\mathbb{F}(F)$ be the frame bundle of $F$, which is a $\mathrm{GL}(\dim_\mathbb{C}(F))$-principal bundle over $M$. Since $H^{(0,0)}(N,L^p\otimes\xi\otimes \mathbb{C}^{\dim_\mathbb{C}(F)})\cong H^{(0,0)}(N,L^p\otimes\xi)\otimes\mathbb{C}^{\dim_\mathbb{C}(F)}$, we could always assume in what follows that $F=\mathbb{C}$ trivial, or we may replace $(G,N,L,\xi)$ by $(G\times\mathrm{GL}(\dim_\mathbb{C}(F)),N,L,\xi\otimes \mathbb{C}^{\dim_\mathbb{C}(F)})$.

Put $\mathscr{F}_p=P_G\times_G\mathscr{C}^\infty(N,L^p\otimes {{{\xi}}})$, which is an infinite dimensional bundle on $M$ and $F_p$ is a subbundle of $\mathscr{F}_p$. The connection $\nabla^{L^p\otimes {\xi}}$ naturally induces a flat connection on $\nabla^{\mathscr{F}_p}$: if $s$ is a smooth section of $\mathscr{F}_p$ and $U\in TM$, set
\begin{equation}
\nabla^{\mathscr{F}_p}_Us=\nabla^{L^p\otimes {{{\xi}}}}_{U^H_0}s.
\end{equation}
This connection preserves $F_p$, and it induces a flat connection $\nabla^{F_p}$ on $F_p$. We equip $\mathscr{F}_p$ with the $L^2$-metric $g^{\mathscr{F}_p}$ induced by $(g^{L^p},h^{\xi},dv_{N})$, which gives a metric $g^{F_p}$ on $F_p$. Let $P_p$ denote the fibrewise orthogonal projection $\mathscr{F}_p\rightarrow F_p$.

Let $\nabla^{\mathcal{N}}$ be a connection on the infinite dimensional bundle $\mathscr{C}^\infty(\mathcal{N})\to M$ given as follows: if $U\in TM$ and $\h\in\mathscr{C}^\infty(\mathcal{N})$, put $\nabla^{\mathcal{N}}_U=U_0^H(\h)$. Following \cite[(9.34)]{bmz17}, we set
\begin{equation}
\vartheta^L=-\frac{1}{2}\omega\big(L,g^{L}\big),\ \ \ \  \vartheta^\xi=-\frac{1}{2}\omega\big(\xi,g^{\xi}\big).
\end{equation}

We denote by $\widehat{TX}$ a copy of $TX$. We denote by $\widehat{\omega}\big(\nabla^F,g^F\big),\widehat{\vartheta}^L,\widehat{\mathrm{div}}_N$ the restrictions of $\omega\big(\nabla^{F},g^{F}\big),{\vartheta}^L,\mathrm{div}_N$ to $\widehat{TX}$. Let $\{e_i\}_{i=1}^m$ and $\{\widehat{e}_i\}_{i=1}^m$ be a local orthonormal frame of $TX$ and $\widehat{TX}$ respectively. We denote the anticommutator by $[\cdot,\cdot]_+$.
\begin{theo}[{\cite[Theorem 9.27]{bmz17}}]\label{Db4}
	We have the following identity:
	\begin{equation}\label{db16}
		\begin{split}
			\frac{1}{2p}\omega\big(\nabla^{F_p},g^{F_p}\big)=T_{-\vartheta^L-\vartheta^\xi/p+\mathrm{div}_N/2p,p}.
		\end{split}
	\end{equation}
	If $\mathcal{H}\in\mathscr{C}^\infty(N,\mathrm{End}({{\xi}}))$ is a smooth section, we have
	\begin{equation}\label{db17}
		\begin{split}
\nabla^{F_p,u}T_{\mathcal{H},p}=T_{\nabla^\mathcal{N}\mathcal{H},p}+p[T_{\vartheta^L+\vartheta^\xi/p-\mathrm{div}_N/2p,p},T_{\mathcal{H},p}]_+-2pT_{(\vartheta^L+\vartheta^\xi/p-\mathrm{div}_N/2p)\mathcal{H},p}.
		\end{split}
	\end{equation}
Moreover, combined Theorem \ref{Cb7} with \eqref{db16} and \eqref{db17}, we see that $\frac{1}{4p}\omega(\nabla^{F_p},g^{F_p})^2$, $\frac{1}{4p^2}\widehat{\omega}(\nabla^{F_p},g^{F_p})(e_i)^2$ and $\nabla^{F_p,u}T_{\mathcal{H},p}$ are Toeplitz operators.
\end{theo}

\begin{remark}
For a general fibration, Ma-Zhang \cite[Theorems 1.18, 1.19]{mz22}, \cite{MR2286586} showed that the curvature operator on $F_p$ is a Toeplitz operator.
\end{remark}

\subsection{The smooth bundle of Toeplitz operators}\label{Dc}

Analogous to Definition \ref{Cb2}, we call $\{T_p\mid T_p\in\mathscr{C}^\infty\big(M,\mathrm{End}({F}_p)\big)\}_{{p\in{\mathbb{N}^*}}}$ a smooth section of Toeplitz operators if there exists $\{\mathcal{H}_i\mid \mathcal{H}_i\in \mathscr{C}^\infty(\mathcal{N},\mathrm{End}({{{\xi}}}))\}_{i\in \mathbb{N}}$ such that for any $k,\ell\in\mathbb{N}$, there is  $C>0$ satisfies
\begin{equation}\label{dc1}
	\Bv T_p-\sum_{i=0}^k p^{-i}T_{\mathcal{H}_i,p}\Bv_{\mathscr{C}^\ell(M,\mathrm{End}({F}_p))}\leqslant Cp^{-k-1},
\end{equation}
where the norm $\lv\cdot\rv_{\mathscr{C}^\ell(M,\mathrm{End}(F_p))}$ is induced by $(\nabla^{F_p},g^{F_p})$.

\begin{remark}\label{dc6}
As all the expansions in \cref{sC} are smoothly dependent on geometric data $(J,g^{T_\mathbb{R}N},h^{L},h^{E})$, all results in \Cref{sC} have a similar version for $F_p$. In particular, by \eqref{cb24}, all the operators in Theorem \ref{Db4} are smooth sections of Toeplitz operators as well as their derivatives. Also, for $\{T_p\}_{p\in\mathbb{N}^*}$ in \eqref{dc1}, $\tro^{F_p}[T_{p}]$ and $\tro_{[k]}[T_{p}]$ are smooth functions on $M$, and we can replace the absolute values $\lv\cdot\rv$ on the left hand side of \eqref{cb20}, \eqref{cc41} and \eqref{cc42} by a smooth norm $\lv\cdot\rv_{\mathscr{C}^{\ell}(M)}$
\end{remark}

\subsection{The asymptotics of $h(\nabla^{F_p},g^{F_p})$}\label{Dd}
Recall the odd form $h(\nabla^{F_p},g^{F_p})$ of $F_p$ as in \eqref{bb9}.
\begin{prop}
	As $p\rightarrow+\infty$, we have smooth odd forms $\gamma_i\in\Omega^\mathrm{odd}(M),i\in\mathbb{N}$ such that for any $k,\ell\in\mathbb{N}$, there is $C>0$ such that
	\begin{equation}\label{dd1}
		\Big\vert p^{-n}\frac{1}{\sqrt{p}}\psi_{1/\sqrt{p}}h\big(\nabla^{F_p},g^{F_p}\big)-\sum_{i=0}^k\gamma_ip^{-i}\Big\vert_{\mathscr{C}^\ell(M)}\leqslant Cp^{-k-1}.
	\end{equation}
\end{prop}

\begin{pro}
	By Theorem \ref{Db4}, Remark \ref{dc6} and \eqref{bb9}, if we take
	\begin{equation}
\gamma_i=({2\pi i})^{1/2}\varphi\tro_{[k]}\left[\exp\Big(\frac{1}{4p}\omega\big(\nabla^{F_p},g^{F_p}\big)^2+\frac{1}{2p}z\omega\big(\nabla^{F_p},g^{F_p}\big)\Big)\right]^z,
	\end{equation}
then we get the expansion \eqref{dd1}.\qed
\end{pro}

\subsection{Spectral gap}\label{De}

Now we recall the non-degeneracy condition \cite[Definition 9.13]{bmz17}.
\begin{defi}
We say that $\widehat{\vartheta}^L$ is nondegenerated if there is $\bm{a}>0$ such that
\begin{equation}\label{de1}
\sum_{i=1}^{m}\widehat{\vartheta}^L(e_i)^2\geqslant \bm{a} \ \ \text{on\ }\mathcal{N},
\end{equation}
or equivalently, $\widehat{\vartheta}^L\in\mathscr{C}^\infty(\mathcal{N},q^*\widehat{T^*X})$ is nowhere vanishing on $\mathcal{N}$.
\end{defi}
\begin{theo}[{\cite[Theorem 4.4]{bmz17}}]\label{De8}
	If $\widehat{\vartheta}^L$ is nondegenerate, for any $\varepsilon>0$ and $p\in\mathbb{N}^*$ large, the fibrewise Dirac operator of $F_p$ (see \eqref{be6}) satisfies
	\begin{equation}\label{de2}
{D}^{F_p,2}_X\geqslant (\bm{a}-\varepsilon)p^2.
	\end{equation}
In particular, for $p\in\mathbb{N}^*$ large enough, $D^{F_p,2}_X$ is invertible, thus $H^\bullet(X,F_p)=0$. by the Hodge Theorem \eqref{be6}.
\end{theo}

\section{The asymptotics of the odd superconnection forms}\label{sF}

The purpose of this section is to obtain the asymptotics of the odd form $h\big( A',g_{4t/p^2}^{\Omega^\bullet(X,F_p)}\big)$ as $p\rightarrow+\infty$. The main technique we use is the analytic localization method by Bismut-Lebeau \cite{bisleb91}, which is further developed by Dai-Liu-Ma \cite{10.4310/jdg/1143593124} and Ma-Marinescu \cite[\S\,4]{mm07}. 

We treat the two cases $1\leqslant t<+\infty$ and $0\leqslant  t\leqslant 1$ separately in \cref{Fa}-\cref{Fg} and \cref{Fi} respectively. The main difficulty is to get an exponential decay term when $t\to+\infty$ as in Theorem \ref{F1}, to do so,
we devote a large part of this section to carefully handling the localization procedure to ensure that at each step, the operator we get is always ``\emph{non-negative}'' as $p^{-2}{D}^{F_p,2}_X$ in Theorem \ref{De8} (see \cref{Fc}-\cref{Fd}).

Now we give more details on the main points of this section. In \cref{Fa}, we state the main theorem of this section and express the odd superconnection forms in terms of heat kernels. In \cref{Fb}, by the spectral gap property of $\mathcal{L}^{F_p}$, we can use the finite propagation speed of solutions of hyperbolic equations and localize the original question to a problem on $\mathbb{R}^{m}$. In \cref{Fc}, we perform the Bismut-Zhang's rescaling on the operator $\theta_{1/\sqrt{p}}^{-1}p^{-2}\mathcal{L}^{F_p}\theta_{1/\sqrt{p}}$, then we introduce a smooth family of differential operators $\widehat{\mathscr{L}}^{F_p}_v|_{0\leqslant v\leqslant1/p}$ as in \eqref{4.2.9} which links the rescaled operator and a limit operator, and we study the Taylor expansion of $\widehat{\mathscr{L}}^{F_p}_v$ with respect to the parameter $v$. In \cref{Fd}, we introduce graded Sobolev norms with weights $\lV\cdot\rV_{\mu,k,p}$ and give some basic elliptic estimations of $\widehat{\mathscr{L}}^{F_p}_v$, an important step is to carefully analyze the structure of $\widehat{\mathscr{L}}^{F_p}_v$ and choose a suitable weight $\mu$ to preserve the ``positivity" of $\mathrm{Spec}(\widehat{\mathscr{L}}^{F_p}_v)$. In \cref{Fe}, we study the Sobolev  estimations of the resolvent $(\lambda-\widehat{\mathscr{L}}^{F_p}_v)^{-1}$. In \cref{Ff}, we establish the corresponding convergence of the heat kernel. In \cref{Fg}, we analyze the asymptotic trace of the limit kernel when $p\to+\infty$ using results in \cref{Cc} and properties of Gaussian integral. In \cref{Fi}, we deal with the case $0<t\leqslant1$. In \cref{Fh}, we prove the main theorem (see Theorem \ref{F1}) of this section.

\subsection{Asymptotic expansion of the odd forms}\label{Fa}

The constants appearing in the sequel depend on the compact subset of $S$ we are working on. To simplify the statements in what follows, we always assume that $S$ is compact. The following theorem is the main result of this section, and the rest of the section is devoted to its proof.

Recall the odd form $h\big( A',g^{\Omega^\bullet(X,F_p)}_{t}\big)$ and $\theta_a$ in in \eqref{bh2} and \eqref{bj13} respectively.
\begin{theo}\label{F1}
If $\widehat{\vartheta}^L$ is nondegenerate as in \eqref{de1}, there exist $\{c_i(t)\in\Omega^\bullet(M)\mid t>0\}_{i\in\mathbb{N}}$ such that, for any ${\gamma}>0,k,\ell \in\mathbb{N}$, we have $C>0$ that for ${p\in{\mathbb{N}^*}}$ and $t>0$,
	\begin{equation}\label{fa1}
		\begin{split}
			\bbv p^{-n}\frac{1}{\sqrt{p}}\psi_{1/\sqrt{p}}h\big( A',g^{\Omega^\bullet(X,F_p)}_{4t/p^2}\big)-\sum_{i=0}^k\int_Xc_i(t)p^{-i}\bbv_{\mathscr{C}^\ell(S)}\leqslant C{e}^{-(\bm{a}-{\gamma})t}p^{-k-1},
		\end{split}
	\end{equation}
where $\bm{a}\geqslant0$ is given in \eqref{de1}. Moreover, there is $C>0$ such that  for $t>0$, we have
\begin{equation}\label{fa3}
\bbv\int_Xc_k(t)\bbv_{\mathscr{C}^\ell(S)}\leqslant C_{}{e}^{-(\bm{a}-{\gamma})t}.
\end{equation}
\end{theo}
For $p\in\mathbb{N}^*$ and $t>0$, set
\begin{equation}\label{fa.4}
	\mathcal{M}^{F_p}=\theta_{{\sqrt{p}}}p^{-2}{\mathcal{L}}^{F_p}\theta_{{1}/{\sqrt{p}}},\ \ \mathcal{M}_t^{F_p}=\theta_{{\sqrt{p}}/\sqrt{t}}tp^{-2}{\mathcal{L}}^{F_p}\theta_{{\sqrt{t}}/{\sqrt{p}}}.
\end{equation}
Let $\exp(-t\mathcal{M}^{F_p})(x,x')$ and $\exp(-t'\mathcal{M}_t^{F_p})(x,x')$ be the smooth heat kernel of $\mathcal{M}^{F_p}$ and $\mathcal{M}_t^{F_p}$ associated with $dv_X(x')$ respectively. By \eqref{bh2} and \eqref{fa.4}, we have
\begin{equation}\label{fa5}
	\begin{split}
		&\frac{1}{\sqrt{p}}\psi_{1/\sqrt{p}}h\big( A',g^{\Omega^\bullet(X,F_p)}_{4t/p^2}\big)\\
		&=\sqrt{{2\pi i}}\varphi\bigg[\theta_{1/\sqrt{p}}\int_{X}\tro_s^{\Lambda(T^*X)\otimes F_p}\big[\exp(-\mathcal{L}^{F_p}_{4t/p^2})(x,x)\big]dv_X(x)\bigg]^z\\
		&=\sqrt{{2\pi i}}\varphi\bigg[\theta_{\sqrt{t}}^{-1}\int_{X}\tro_s^{\Lambda(T^*X)\otimes F_p}\big[\exp(-t\mathcal{M}^{F_p})(x,x)\big]dv_X(x)\bigg]^z\\
		&=\sqrt{{2\pi i}}\varphi\bigg[\int_{X}\tro_s^{\Lambda(T^*X)\otimes F_p}\big[\exp(-\mathcal{M}_t^{F_p})(x,x)\big]dv_X(x)\bigg]^z.
	\end{split}
\end{equation}
Notice that to prove Theorem \ref{F1}, we only need to consider the case $\bm{a}=0$, or we may replace $\mathcal{L}^{F_p}$ with $\mathcal{L}^{F_p}-\bm{a}p^2$ so that we get extra ${e}^{-\bm{a}t}$.

\subsection{Localization of the problem}\label{Fb}

For $x_0\in M$, let $X$ be the fibre containing $x_0$, we mainly work along this fibre. For ${\varepsilon}>0$, let $B^X\lk{x_0},{\varepsilon}\rk$ and $B^{T_{{x_0}}X}\lk 0,{\varepsilon}\rk$  be the open balls in $X$ and $T_{{x_0}}X$ with center ${x_0}$ and $0$ and radius ${\varepsilon}$ respectively. Let $\exp^X_{{x_0}}$ be the exponential map of $(X,g^{TX})$. For ${\varepsilon}$ small, $\exp^X_{{x_0}}\colon B^{T_{{x_0}}X}\lk 0,{\varepsilon}\rk\to B^X\lk{x_0},{\varepsilon}\rk$ is a diffeomorphism, which gives local coordinates by identifying $T_{{x_0}}X$ with $\mathbb{R}^{m}$ via an orthonormal basis $\{e_i\}_{i=1}^m$ of $T_{{x_0}}X$: 
\begin{equation}\label{fb1}
	Z=\lk Z_1,\cdots,Z_{m}\rk\in \mathbb{R}^{m} \longmapsto  Z_ie_i\in T_{x_0}X.
\end{equation}
We will always identify $B^X\lk {x_0},{\varepsilon}\rk$ with $B^{T_{{x_0}}X}\lk 0,{\varepsilon}\rk$ through the isomorphism in \eqref{fb1}.

For $b\in S$, let $\text{inj}(X_b)$ be the injectivity radius of $X_b=\pi^{-1}(b)$. For any ${\varepsilon}>0$ that
\begin{equation}\label{fb2}
0<\varepsilon<  \mathrm{min}_{b\in S}\ \text{inj}(X_b)/8,
\end{equation} 
(we have $\mathrm{min}_{b\in S}\text{inj}(X_b)>0$ since $S$ is compact), by the compactness of $X$, we can choose a \emph{finite set} $\{x_i\}_{i\geqslant 1}\subset X$ such that $\big\{B^X\lk{x_i},{\varepsilon}\rk\big\}_{i\geqslant 1}$ is an open covering of $X$. For now, $\varepsilon$ is \emph{not} fixed, we will assign it a suitable value after Theorem \ref{4.11}.

For $Z\in B^{T_{{x_i}}X}\lk 0,{\varepsilon}\rk$, we identify $F_{p,Z}$ and $\big(\mathbb{R}[z]\hti\pi^*\Lambda\lk T^*S\rk\hti\Lambda\lk T^*X\rk\big)_{Z}$ with $F_{p,x_i}$ and $\big(\mathbb{R}[z]\hti\pi^*\Lambda\lk T^*S\rk\hti\Lambda\lk T^*X\rk\big)_{x_i}$ by parallel transport with respect to the connections $\nabla^{F_p}$ and ${}^0\nabla^{\mathbb{R}[z]\hti\pi^*\Lambda\lk T^*S\rk\hti\Lambda\lk T^*X\rk}$ along the curve $s\to sZ$ for $s\in[0,1]$. Let $\Gamma^{0}$ and $\Gamma^{F_p,u}$ be the corresponding connection forms of ${}^0\nabla^{\mathbb{R}[z]\hti\pi^*\Lambda\lk T^*S\rk\hti\Lambda\lk T^*X\rk}$ and $\nabla^{F_p}$ on $B^{T_{{x_0}}X}\lk 0,{\varepsilon}\rk$ with respect to this trivialization.

\begin{defi}
	Let $\{f_{x_i}\}_{i\geqslant1}$ be a partition of unity with respect to $\{B^X\lk{x_i},{\varepsilon}\rk\}_{i\geqslant1}$. For $k\in\mathbb{N}$, we define the Sobolev norm $\lV\cdot\rV_{H^{k,p}}$ on $\mathscr{C}^{\infty}(X,\mathbb{R}[z]\hti\pi^*\Lambda(T^*S)\hti\Lambda(T^*X)\otimes F_p)$:
\begin{equation}\label{5.2.5}
	\lV s\rV^2_{H^{k,p}}=\sum_i\sum_{\begin{subarray}{c}\alpha\in\mathbb{N}^{m}, \ 0\leqslant \lv\alpha\rv\leqslant k \end{subarray}}\lV\pa^\alpha f_{x_i} s\rV^2_{L^2},
\end{equation}
and we denote by  $H^{k,p}$ its  completion with respect to the norm $\lV \cdot\rV^2_{H^{k,p}}$.
\end{defi}
\begin{remark}
Note that $\nabla^{F_p,u}$ preserves the metric $g^{F_p}$. Hence the two norms $\lV\cdot\rV_{H^{0,p}}$ and $\lV\cdot \rV_{L^2(X,F_p)}$ are equivalent uniformly for $p$ and we will not distinguish between them.
\end{remark}

Similar to \cite[Lemma 1.6.2]{mm07} and \cite[Lemma 3.6]{pu}, we have the following lemma.
\begin{lemma}\label{Fb4}
	For $k\in\mathbb{N}$, there is $C>0$ such that for ${p\in{\mathbb{N}^*}}$ and $s\in H^{2k,p}$, we have
	\begin{equation}\label{fb5}
		\Vert s\Vert_{H^{2k,p}}\leqslant Cp^{2k}\sum_{j=0}^kp^{-2j}\big\Vert \mathcal{L}^{F_p,j}s\big\Vert_{L^2}.
	\end{equation}
\end{lemma}

\begin{pro}
	By \eqref{bj14}, on $B^X\lk{x_i},{\varepsilon}\rk$ we have
	\begin{equation}\label{fb11}
		\mathcal{L}^{F_p}=-\big({e_i}+\Gamma^{F_p,u}(e_i)+\Gamma^{0}(e_i)\big)^2+\Big({\nabla^{TX}_{e_i}e_i}+\Gamma^{F_p,u}(\nabla^{TX}_{e_i}e_i)+\Gamma^{0}(\nabla^{TX}_{e_i}e_i)\Big)+\Lambda^{F_p}.
	\end{equation}
For any $Z\in B^X\lk{x_i},{\varepsilon}\rk$, we have a classical relation (see \cite[Proposition 1.18]{bgv})
\begin{equation}\label{fb7}
	\Gamma^{F_p,u}_{Z}(\pa_j)=\int_{0}^1R^{F_p,u}_{tZ}(\mathcal{R},\pa_j)dt,\ \ \mathrm{for \ }{\mathcal{R}}=Z_i\pa_i.
\end{equation} 
By Theorem \ref{Db4}, \eqref{bb3}, \eqref{bj7} and \eqref{fb7}, we see that $p^{-1}\Gamma^{F_p,u}_{Z}(\pa_j)$ and $p^{-2}\Lambda^{F_p}$ are smooth family of Toeplitz operators. For a smooth family of Toeplitz operators $T_p$, by \eqref{cb23} and \eqref{db17}, $dT_p=\nabla^{F_p,u}T_p-[\Gamma^{F_p,u},T_p]$ is also a smooth family of Toeplitz operators. 

By \eqref{fb11} and the above argument, the operator norms of $p^{-1}\Gamma^{F_p,u}$ and $p^{-2}\Lambda^{F_p}$ are bounded uniformly for $p\in\mathbb{N}^*$ as well as their derivatives, then we get \eqref{fb5} exactly as the proof of \cite[Lemma 1.6.2]{mm07}.\qed
\end{pro}

\begin{remark}\label{Fb5}
Observe that \eqref{fb5} is still true if we replace $\mathcal{L}^{F_p}$ with $\mathcal{L}^{F_p,*}$ or ${D}_X^{F_p,2}$ since they have the same structure as in \eqref{fb11}.
\end{remark}

Choose $f\in\mathscr{C}^\infty(\mathbb{R})$ even, nonincreasing when $t\geqslant 0$ and
\begin{equation}\label{fb20}
f\lk t\rk = 
	\begin{cases}
		1 &\text{if ${\vert} t{\vert}\leqslant \frac{1}{2}$},\\
		0 &\text{if ${\vert} t{\vert}\geqslant 1$}.
	\end{cases}
\end{equation}

\begin{defi}
For $t,h>0$ and $a\in \mathbb{C}$, set
\begin{equation}\label{fb21}
	\begin{split}
		F_{t,h}\lk a\rk&=\int_\mathbb{R}e^{{\sqrt{2}{i}va}}\exp\lk-v^2/2\rk f\big( \sqrt{t}v/h\big)\frac{\de v}{\sqrt{2\pi}},\\
		G_{t,h}\lk a\rk&=\int_\mathbb{R}e^{\sqrt{2}{{i}va}}\exp\lk-v^2/2\rk\big(1-f\big( \sqrt{t}v/h\big)\big)\frac{\de v}{\sqrt{2\pi}},\\
		H_{t,h}\lk a\rk&=\int_\mathbb{R}e^{\sqrt{2}{{i}va}}\exp\lk-v^2/2t\rk\big(1-f\lk v/h\rk\big)\frac{\de v}{\sqrt{2\pi t}}.
	\end{split}
\end{equation}
\end{defi}
The functions $F_{t,h}(a),G_{t,h}(a),H_{t,h}(a)$ are even holomorphic functions, thus there exist holomorphic functions $\widetilde{F}_{t,h}\lk a\rk$, $\widetilde{G}_{t,h}\lk a\rk$ and $\widetilde{H} _{t,h}\lk a\rk$ such that 
\begin{equation}\label{fb22}
	\widetilde{F}_{t,h}\lk a^2\rk=F_{t,h}\lk a\rk,\ \ \ \widetilde{G}_{t,h}\lk a^2\rk=G_{t,h}\lk a\rk,\ \ \ \widetilde{H}_{t,h}\lk a^2\rk=H_{t,h}\lk a\rk.
\end{equation}
Moreover, the restriction of $\widetilde{F}_{t}$ and $\widetilde{G}_{t}$ to $\mathbb{R}$ lies in the Schwartz space $\mathcal{S}\lk\mathbb{R}\rk$, and
\begin{equation}\label{fb23}
	F_{t,h}\lk a\rk+G_{t,h}\lk a\rk=\exp\lk-a^2\rk,\ \ \ \ H_{t,h}\lk a\rk=G_{t,h}\big(\sqrt{t}a\big).
\end{equation}

Now we fix a $c>0$, let $V_c$ be the following subset of the complex plane:
\begin{equation}\label{fb24}
	V_c=\Big\{a\in\mathbb{C}\colon\mathrm{Re}\lk a\rk\geqslant\frac{1}{4c^2}\mathrm{Im}\lk a\rk^2-c^2\Big\}.
\end{equation}
\begin{lemma}\label{Fb7}
For any $h>0$ and ${k}\in\mathbb{N}$, there are $C,c_1,c_2>0$ such that for $t>0$, we have
	\begin{equation}\label{fb25}
		\sup_{a\in V_c}\bv a^k\widetilde{H}_{t,h}\left(a\right)\bv\leqslant C\exp\Big(c_1t-\frac{c_2}{t}\Big).
	\end{equation}
\end{lemma}
\begin{pro}
By proceeding as in \cite[Theorem 4.2.5]{mm07}, when ${\vert}\text{Im}\lk a\rk{\vert}\leqslant c$, as $i^k a^ke^{iva} = \frac{\partial^k}{\partial v^k}e^{{i}va}$, one can integrate by part the expression of $a^kH_{t,h}\lk a\rk$ given in \eqref{fb21} to obtain that
\begin{equation}\label{fb36}
	\sup_{{\vert}\text{Im}(a){\vert}\leqslant c}{\vert}a^kH_{t,h}(a){\vert}\leqslant C\exp\Big(c_1t-\frac{c_2}{t}\Big).
\end{equation}
Note that for $c>0$, $V_c$ is just the image of $\{a\in\mathbb{C}\colon{\vert}\text{Im}(a){\vert}\leqslant c\}$ by the map $\mathbb{C}\to\mathbb{C}\colon a\mapsto a^2$. By \eqref{fb22} and \eqref{fb36}, we get \eqref{fb25}.\qed
\end{pro}

For any ${\delta}>0$, let $\Gamma$ be the contour in $\mathbb{C}$ defined by $\{x\pm{\delta}\mi\mid x\geqslant -{\delta}\}\cup\{-{\delta}+x\mi\mid-{\delta}\leqslant x\leqslant {\delta}\}$. We choose a suitable ${\delta}$ to make $\Gamma\subset V_c$ and $\min_{a\in\Gamma,b\in V_c}\lv a-b\rv=c^2/2$ (see \eqref{fb24}), and we denote this contour $\Gamma$ by $\Gamma_c$.
\begin{figure}[ht]
	\centering
\begin{tikzpicture}[decoration={markings, mark=at position 0.5 with {\arrow{>}}}] 
	\draw[postaction={decorate}] (3,1.56)--(2.9,1.56);
	\draw[postaction={decorate}] (2.9,-1.56)--(3,-1.56);
	\draw[->] (-2,0)--(7,0);
	\draw[->] (0,-2)--(0,2);
	\draw[postaction={decorate}](-0.5,0.5)--(-0.5,-0.5);
	\draw[postaction={decorate}](6,0.6)--(-0.5,0.5);
	\draw[postaction={decorate}](-0.5,-0.5)--(6,-0.6);
	\node at (7.2,0){\begin{footnotesize}$x$\end{footnotesize}};
	\node at (0,2.2){\begin{footnotesize}$y$\end{footnotesize}};
	\node[fill=white] at (2.8,1){$\Gamma_c$};
	\node[fill=white] at (5,1.4){$V_c$};
	\begin{scope}[rotate=-90]
		\draw[domain=-1.9:1.9] plot (\x,1.75*\x*\x-1.3);
	\end{scope}
\end{tikzpicture}
	\caption{}
	\label{fig2}
\end{figure}

\begin{lemma}
	For $k\in\mathbb{N}$, there are $C>0,\ell\in\mathbb{N}^*$ such that for any ${p\in{\mathbb{N}^*}}$, $\lambda\in\Gamma_c$,
	\begin{equation}\label{fb50}
		\big\Vert\big(\lambda-\mathcal{L}^{F_p}\big)^{-1}s\big\Vert_{H^{2k+2,p}} \leqslant Cp^{\ell}|\lambda|^{\ell}\Vert s\Vert_{H^{2k,p}}.
	\end{equation}
\end{lemma}

\begin{pro}
	First, following \cite[(3.65)]{pu} we prove that for some $C>0,i\in\mathbb{N}^*$,
	\begin{equation}\label{4.1.39}
		\big\Vert\big(\lambda-\mathcal{L}^{F_p}\big)^{-1}s\big\Vert_{L^2} \leqslant Cp^{i}|\lambda|^{i}\Vert s\Vert_{L^2}.
	\end{equation}
	By Theorem \ref{2.19}, we write $\mathcal{L}^{F_p}={D}_X^{F_p,2}+\mathscr{G}^{F_p}$ where $\mathscr{G}^{F_p}$ is a $1$-order fibrewise differential operator with positive degree in $\pi^*\Lambda\lk T^*S\rk$, and
\begin{equation}\label{4.1.35}
	\big(\lambda-\mathcal{L}^{F_p}\big)^{-1}=\sum_{i=0}^{\dim S}\big(\big(\lambda-{D}_X^{F_p,2}\big)^{-1}\mathscr{G}^{F_p}\big)^i\big(\lambda-{D}_X^{F_p,2}\big)^{-1}.
\end{equation}
Note the self adjointness of ${D}_X^{F_p,2}$, for any $\lambda\in\Gamma_c$ and $x\geqslant 0$, we have $\max\{\frac{1}{\lv\lambda-x\rv},\frac{x}{|\lambda-x|}\}\leqslant C\frac{|\lambda|}{\vert\lambda-x\vert}+1\leqslant C|\lambda|$, so there is $C>0$ with
\begin{equation}\label{4.1.36}
	\big\Vert\big(\lambda-{D}_X^{F_p,2}\big)^{-1}s\big\Vert_{L^2}\leqslant C\Vert s\Vert_{L^2},\ \ \big\Vert {D}_X^{F_p,2}\big(\lambda-{D}_X^{F_p,2}\big)^{-1}s\big\Vert_{L^2}\leqslant C|\lambda|\cdot\Vert s\Vert_{L^2}.
\end{equation}
By Theorem \ref{Db4}, Lemma \ref{Fb4}, Remark \ref{Fb5}, \eqref{bj7} and \eqref{bj14}, there is $C>0$ such that for any $\lambda\in\Gamma$,
\begin{equation}\label{4.1.38}
	\begin{split}
		&\big\Vert \mathscr{G}^{F_p}\big(\lambda-{D}_X^{F_p,2}\big)^{-1}s\big\Vert_{L^2}\leqslant Cp^{2}\big\Vert\big(\lambda-{D}_X^{F_p,2}\big)^{-1}s\big\Vert_{H^{2,p}}\\
		&\leqslant Cp^{2}\Big(\big\Vert {D}_X^{F_p,2}\big(\lambda-{D}_X^{F_p,2}\big)^{-1}s\big\Vert_{L^2}+p^2\Vert\big(\lambda-{D}_X^{F_p,2}\big)^{-1}s\Vert_{L^2}\Big)\\
		&\leqslant Cp^{4}|\lambda\vert\cdot \Vert s\Vert_{L^2}.
	\end{split}
\end{equation}
By \eqref{4.1.35}, \eqref{4.1.36} and \eqref{4.1.38}, we get \eqref{4.1.39}. From \eqref{4.1.35} and an argument similar to \eqref{4.1.38}, for some $C>0,j\in\mathbb{N}$, we have
\begin{equation}\label{fb52}
	\begin{split}
		\big\Vert\mathcal{L}^{F_p}(\lambda-\mathcal{L}^{F_p}\big)^{-1}s\big\Vert_{L^2} \leqslant Cp^{j}|\lambda|^{j}\Vert s\Vert_{L^2}
	\end{split}
\end{equation}
By \eqref{fb52}, for any $k\in\mathbb{N}$, there are $C>0,\ell\in\mathbb{N}^*$ such that 
\begin{equation}\label{fb53}
	\begin{split}
		\big\Vert\mathcal{L}^{F_p,k+1}(\lambda-&\mathcal{L}^{F_p}\big)^{-1}s\big\Vert_{L^2}=\big\Vert\mathcal{L}^{F_p}(\lambda-\mathcal{L}^{F_p}\big)^{-1}\mathcal{L}^{F_p,k}s\big\Vert_{L^2}\\
		&\leqslant Cp^{j}|\lambda|^{j}\Vert \mathcal{L}^{F_p,k}s\Vert_{L^2}\leqslant Cp^{\ell}|\lambda|^{\ell}\Vert s\Vert_{H^{2k,p}}.
	\end{split}
\end{equation}
By \eqref{fb5} and \eqref{fb53}, we get \eqref{fb50}.\qed

\end{pro}

Set the fibre product $M\times_SM=\{(x_1,x_2)\in M\times M\mid \pi(x_1)=\pi(x_2)\}$ with the projections $\mathrm{pr}_i(x_1,x_2)\to x_i$. For any two bundles $V,V'$ over $M$, Let $V\times_S V'$ be the bundle on $M\times_SM$ given by $\mathrm{pr}_1^*(V)\otimes \mathrm{pr}_2^*(V')$. If we have two connections $\nabla^V,\nabla^{V'}$ on $V$ and $V'$, let $\nabla^{V\times_S V'}$ be the induced connection on $V\times_S V'$. We define a bundle
\begin{equation}
	\mathbb{E}_p=\mathbb{R}[z]\hti\pi^*\Lambda\lk T^*S\rk\hti\Big(\big(\Lambda(T^*X)\otimes F_p\big)\times_S\big((\Lambda(T^*X))^*\otimes F^*_p\big)\Big)
\end{equation}
over $M \times_S M$. Let $\nabla^{\mathbb{E}_p}$ be the connection on $\mathbb{E}_p$ induced by ${}^{0}\nabla^{\mathbb{R}[z]\hti\pi^*\Lambda^\bullet\lk T^*S\rk\hti\Lambda^\bullet\lk T^*X\rk\otimes F_p,u}$ and $\nabla^{F_p,u}$ (see \eqref{bb3} and \eqref{bi.33}). We notice that $\widetilde{G}_{t,{\varepsilon}}(t\mathcal{L}^{F_p})$ is a smooth section of  $\mathbb{E}_p$ on $M \times_S M$.

We also denote the pull back bundle of $\mathbb{E}_p$ through the diagonal embedding $M\rightarrow M \times_S M$ given by $x_0\mapsto (x_0,x_0)$, then
\begin{equation}\label{fb.57}
\mathbb{E}_{p,x_0}=\mathbb{R}[z]\hti\pi^*\Lambda\lk T_{\pi(x_0)}^*S\rk\hti \mathrm{End}(\Lambda(T_{x_0}^*X))\otimes \mathrm{End}(F_{p,x_0}).
\end{equation}

\begin{theo}\label{Fb10}
	For any ${\varepsilon}>0$ and $\ell\in\mathbb{N}$, there exist $c_1,c_2>0$ and $k\in\mathbb{N}$ such that for any $t>0$ and ${p\in{\mathbb{N}^*}}$, we have
	\begin{equation}\label{4.1.30}
		\big\vert\widetilde{G}_{t,{\varepsilon}}\big(t\mathcal{L}^{F_p}\big)\big\vert_{\mathscr{C}^\ell\lk M\times_SM\rk}\leqslant Cp^{k}\exp\Big(c_1t-\frac{c_2}{t}\Big),
	\end{equation} 
	where the $\mathscr{C}^\ell\lk M\times_SM\rk$-norm is induced by $\nabla^{\mathbb{E}_p}$.
\end{theo}

\begin{pro}
By \eqref{fb25}, for $r\in\mathbb{N}^*$, there is a holomorphic function $\widetilde{H}_{t,h,r}(a)$ defined in a neighborhood of $V_c$ such that (see also \cite[Theorem 13.32]{bisleb91})
\begin{equation}\label{fb60}
	\frac{1}{(r-1)!}\frac{d^{r-1}}{da^{{r-1}}}\widetilde{H}_{t,h,r}(a)=\widetilde{H}_{t,h}(a),
\end{equation}
and for any $h>0$, there are $C,c_1,c_2>0$ such that for any $t>0$, we have
\begin{equation}\label{fb61}
	\sup_{a\in V_c}{\vert}a^k\widetilde{H}_{t,h,r}\left(a\right){\vert}\leqslant C\exp\Big(c_1t-\frac{c_2}{t}\Big).
\end{equation}
Let $\{U_i\}_{i=1}^{\dim S}$ be a set of local vector fields on $S$. Let $\{U^H_i\}$ denote the set of corresponding horizontal lift local vector fields. For any multi-index $\alpha=(i_1,i_2,\cdots)$ in which $1\leqslant i_j\leqslant \dim S$, we denote $U^H_\alpha=\nabla_{U_{i_1}^H}^{\mathbb{E}_{p}}\nabla_{U_{i_2}^H}^{\mathbb{E}_{p}}\cdots$. By \eqref{fb23} and \eqref{fb60}, we get
\begin{equation}\label{fb62}
U^H_\alpha\big(\widetilde{G}_{t,{\varepsilon}}\big(t\mathcal{L}^{F_p}\big)\big)=\frac{1}{{2\pi i}}\int_{\Gamma}\widetilde{H}_{t,{\varepsilon},r}(\lambda)U^H_\alpha(\lambda-\mathcal{L}^{F_p}\big)^{-r}d\lambda.
\end{equation}
Observe that $U^H_\alpha(\lambda-\mathcal{L}^{F_p}\big)^{-r}$ is a linear combination of operators in the form of
\begin{equation}\label{fb63}
(\lambda-\mathcal{L}^{F_p}\big)^{-r_0}U^H_{\alpha_1}\big(\mathcal{L}^{F_p}\big)(\lambda-\mathcal{L}^{F_p}\big)^{-r_1}\cdots U^H_{\alpha_\ell}\big(\mathcal{L}^{F_p}\big)(\lambda-\mathcal{L}^{F_p}\big)^{-r_\ell},
\end{equation}
where $\alpha_j=(j_1,j_2,\cdots)$ is multi-index, $\ell\in\mathbb{N}$, $r_i\geqslant 1$, $\sum r_i=r+\ell$. Note that for any multi-index $\alpha$, $U^H_{\alpha}\big(\mathcal{L}^{F_p}\big)$ is a second order fibrewise differential operator satisfies
\begin{equation}\label{fb64}
\bV U^H_{\alpha}\big(\mathcal{L}^{F_p}\big)s\bV_{H^{k,p}}\leqslant Cp^2\bV s\bV_{H^{k+2,p}}.
\end{equation}
By \eqref{fb50} and \eqref{fb64}, for any $k\in\mathbb{N}$, there are $i,i'\in\mathbb{N}^*$ such that
\begin{equation}\label{fb35}
\big\Vert U^H_{\alpha}\big(\mathcal{L}^{F_p}\big)\big(\lambda-\mathcal{L}^{F_p}\big)^{-1}s\big\Vert_{H^{2k,p}}\leqslant p^{i}\vert\lambda\vert^{i'}\Vert s\Vert_{H^{2k,p}}.
\end{equation}

Let $P, Q$ be differential operators of order $2q$ and $2q'$ with scalar principal symbol and
with compact support in $B^X\lk{x},{\varepsilon}\rk$ and $B^X\lk{x'},{\varepsilon}\rk$ respectively. We take $r\geqslant(2q+2q'+2)\vert\alpha\vert$ in \eqref{fb62}, then there is a $r_j\geqslant 2(q+q')$. By \eqref{fb63}, we split the operator $PU^H_\alpha\big((\lambda-\mathcal{L}^{F_p})^{-k}\big)Q$ into the product of two parts
\begin{equation}\label{fb65}
\begin{split}
&P(\lambda-\mathcal{L}^{F_p}\big)^{-r_0}U^H_{\alpha_1}\big(\mathcal{L}^{F_p}\big)(\lambda-\mathcal{L}^{F_p}\big)^{-r_1}\cdots U^H_{\alpha_{j}}\big(\mathcal{L}^{F_p}\big)(\lambda-\mathcal{L}^{F_p}\big)^{-1}(\lambda-\mathcal{L}^{F_p}\big)^{-2q},\\
&(\lambda-\mathcal{L}^{F_p}\big)^{r_j-2q'-2q-2}(\lambda-\mathcal{L}^{F_p}\big)^{-2q'-1}\cdots U^H_{\alpha_\ell}\big(\mathcal{L}^{F_p}\big)(\lambda-\mathcal{L}^{F_p}\big)^{-r_\ell}Q.
\end{split}
\end{equation}
By \eqref{fb50} and \eqref{fb35}, the first part above is a map from $L^2$-space to itself, and the operator norm is dominated by $Cp^{j_0}\vert\lambda\vert^{j_0'}$. By Remark \ref{Fb5}, the adjoint of the second part has the same structure as the first part, and its operator norm is also dominated by $Cp^{j_1}\vert\lambda\vert^{j_1'}$. Therefore, for some $j,j'\in\mathbb{N}^*$, we have
\begin{equation}\label{fb66}
\bV PU^H_\alpha\big((\lambda-\mathcal{L}^{F_p})^{-k}\big)Qs\bV_{L^2}\leqslant Cp^{j}\vert\lambda\vert^{j'}\Vert s\Vert_{L^2}.
\end{equation}
By \eqref{fb61}, \eqref{fb62}, \eqref{fb66} and the Sobolev inequality, we find that
\begin{equation}
\big|U^H_\alpha\big(\widetilde{G}_{t,{\varepsilon}}\big(t\mathcal{L}^{F_p}\big)\big)(\cdot,\cdot)\big|_{\mathscr{C}^{q}(X\times X)}\leqslant Cp^{k}\exp\Big(c_1t-\frac{c_2}{t}\Big),
\end{equation}
which gives \eqref{4.1.30}. Here we emphasize that the constant in Sobolev inequality is \emph{independent} on the dimension of bundle $F_p$. \qed
\end{pro}

\subsection{Rescaling of the operator $\mathcal{L}^{F_p}$}\label{Fc}

Now we use Bismut-Zhang's rescaling \cite[Chapter 4]{bz92} for two kinds of Clifford variables, which is analogous to Getzler's rescaling \cite{get86}.

Using the same notation as in the beginning of \Cref{Fb}, we fix $x_0\in M$ and identify $\mathbb{R}[z]\hti\pi^*\Lambda\lk T^*S\rk\hti\Lambda\lk T^*X\rk$ and $F_{p}$ to $\big(\mathbb{R}[z]\hti\pi^*\Lambda\lk T^*S\rk\hti\Lambda\lk T^*X\rk\big)_{x_0}$ and $F_{p,x_0}$ on $B^{T_{{x_0}}X}\lk 0,{\varepsilon}\rk$. Let $\Gamma^{0}$ and $\Gamma^{F_p,u}$ be the corresponding connection forms .

Put $X_0=T_{x_0}X\cong \mathbb{R}^m$. Let $g^{TX_0}$ be a Riemannian metric on $X_0$ such that
\begin{equation}\label{fb.73}
	g^{TX_0}=\begin{cases}
		g^{TX}, &\text{ on } B^{T_{x_0}X}\lk0,\mathrm{inj}(X)/2\rk,\\
		g^{T_{x_0}X}, &\text{ outside } B^{T_{x_0}X}\lk0,\mathrm{inj}(X)\rk,
	\end{cases}
\end{equation}
and let $dv_{X_0}$ be the associated volume form. Let $d v_{TX}$ be the Riemannian volume form of $\lk T_{x_0}X,g^{T_{x_0}X} \rk$, and $\kappa\lk x\rk$ be the smooth positive function defined by
\begin{equation}\label{4.1.48}
	d v_{X_0}=\kappa\lk Z\rk d v_{TX},\ \ \kappa\lk0\rk = 1.
\end{equation}

Recall the function $f$ in \eqref{fb20}. For any ${\varepsilon}>0$ satisfies \eqref{fb2}, set $f_{\varepsilon}\colon \mathbb{R}^m\to \mathbb{R}$ by 
\begin{equation}\label{fb.76}
f_{\varepsilon}(Z)=f(\vert Z\vert/{\varepsilon}).
\end{equation}
Put
\begin{equation}
\mathbb{F}_{p}=\mathbb{R}[z]\hti\pi^*\Lambda(T^*S)\hti\Lambda(T^*X)\otimes F_p.
\end{equation}
On the trivial bundle $\mathbb{F}_{p,x_0}$ over $X_0$, we define
\begin{equation}\label{4.1.46}
\nabla^{\mathbb{F}_{p,x_0}}=\nabla+f_{\varepsilon}(Z)\big(\Gamma_Z^{0}+\Gamma_Z^{F_p,u}\big).
\end{equation}
Let ${{\Delta}}^{\mathbb{F}_{p,x_0}}$ be the Laplacian associated with $\nabla^{\mathbb{F}_{p,x_0}}$ and $g^{TX_0}$. Let $\nabla^{TX_0}$ be the Levi-Civita connection of $\lk X_0,g^{TX_0}\rk$ and let $(g^{ij})$ be the inverse of $(g_{ij})=\big(g^{TX_0}(\partial_i,\pa_j)\big)$, then
\begin{equation}\label{4.1.49}
{{\Delta}}^{\mathbb{F}_{p,x_0}}=-g^{ij}\big(\nabla^{\mathbb{E}_{p,x_0}}_{\pa_i}\nabla^{\mathbb{E}_{p,x_0}}_{\pa_j}-\nabla^{\mathbb{E}_{p,x_0}}_{\nabla^{TX_0}_{\pa_i}\pa_j}\big).
\end{equation}
Recall $\Lambda^{F_p}$ in \eqref{bj7}, locally we view $\Lambda^{F_p}$ as an element in $\mathscr{C}^\infty(B^{T_{x_0}X}(0,\varepsilon),\mathbb{E}_{p,x_0})$ (see \eqref{fb.57}). Set 
\begin{equation}\label{4.1.50}
	{\Lambda}_{x_0,p}=f_{\varepsilon}(Z)\Lambda^{F_p}\in \mathscr{C}_0^\infty(T_{x_0}X,\mathbb{E}_{p,x_0}),\ \ \ \ \ {\mathscr{L}}^{F_p}_{x_0}={{\Delta}}^{\mathbb{E}_{p,x_0}}+{\Lambda}_{x_0,p}.
\end{equation}
Then $\mathscr{L}^{F_p}_{x_0}$ is a family of differential operators acting on $\mathscr{C}^\infty(T_{x_0}X,\mathbb{F}_{p,x_0})$ for $x_0\in M$.

Let $\exp(-t\mathcal{L}^{F_p})\lk x,x'\rk,(x,x')$ be the smooth kernel of $\exp(-t\mathcal{L}^{F_p})$ with respect to $dv_{X}$ for $\in M\times_SM$, and $\exp(-t{\mathscr{L}}^{F_p}_{x_0})\lk Z,Z'\rk$ the smooth kernel of $\exp(-t{\mathscr{L}}^{F_p}_{x_0})$ with respect to $dv_{X_0}$ for $(Z,Z')\in T_{x_0}X\times T_{x_0}X$. As in \eqref{fb.57}, we still denote the pull back bundle of $\mathbb{E}_p$ through the projection $TX\times_MTX\to M$ by $\mathbb{E}_p$, then
we can view $\exp(-t{\mathscr{L}}^{F_p}_{x_0})\lk Z,Z'\rk$ as a smooth section of $\mathbb{E}_{p}$ on $TX\times_M TX$.

\begin{theo}\label{4.5}
	For any ${\varepsilon}>0$ satisfies \eqref{fb2} and $\ell\in\mathbb{N}$, there exist $C,c_1,c_2>0$ and $k\in\mathbb{N}$ such that for any ${p\in{\mathbb{N}^*}}, t>0$, we have
	\begin{equation}\label{4.1.52}
		\Bv\exp(-t\mathcal{L}^{F_p})(x_0,x_0)-\exp(-t{\mathscr{L}}^{F_p}_{x_0})(0,0)\Bv_{\mathscr{C}^\ell(M)}\leqslant Cp^{k}e^{c_1t-\frac{c_2}{t}},
	\end{equation}
where $\lv\cdot\rv_{\mathscr{C}^\ell(M)}$ is the $\mathscr{C}^\ell$ norm with respect to the parameter $x_0\in M$.
\end{theo}
\begin{pro}
	By \eqref{4.1.46}, \eqref{4.1.49} and \eqref{4.1.50}, ${\mathscr{L}}^{F_p}_{x_0}$ and $\mathcal{L}^{F_p}$ coincide over $B^{T_{x_0}X}\lk0, {\varepsilon}\rk$, then from \eqref{fb20}, \eqref{fb21}, \eqref{fb22} and the finite propagation speed of wave operator (see \cite[Theorem D.2.1]{mm07}), we obtain that
\begin{equation}
	\widetilde{F}_{t,\varepsilon}\big(t\mathcal{L}^{F_p}\big)\lk x_0,\cdot\rk=\widetilde{F}_{t,\varepsilon}\big( t{\mathscr{L}}^{F_p}_{x_0}\big)\lk0,\cdot\rk.
\end{equation}
	By \eqref{4.1.50}, ${\mathscr{L}}^{F_p}_{x_0}$ has the same structure as $\mathcal{L}^{F_p}$, especially, Lemmas \ref{Fb4} and \ref{Fb10} are still true if we replace $\mathcal{L}^{F_p}$ with ${\mathscr{L}}^{F_p}_{x_0}$. Then \eqref{4.1.52} follows from \eqref{fb23}.\qed
\end{pro}

Put
\begin{equation}\label{fc10}
	{\mathscr{M}}_{x_0}^{F_p}=\theta_{{\sqrt{p}}}p^{-2}{\mathscr{L}}_{x_0}^{F_p}\theta_{{1}/{\sqrt{p}}}.
\end{equation}

\begin{coro}
For any ${\varepsilon}>0$ satisfies \eqref{fb2}, $\gamma>0$ and $\ell\in\mathbb{N}$, there exist $C,c_1,c_2>0$ such that for any $t>0$ and ${p\in{\mathbb{N}^*}}$ large enough, we have
\begin{equation}\label{fb.82}
	\begin{split}
\Bv\theta_{\sqrt{t}}^{-1}\tro_s^{\Lambda(T^*X)\otimes F_p}\big[\exp(-t\mathcal{M}^{F_p})(x_0,x_0)-\exp(-t{\mathscr{M}}^{F_p}_{x_0})(0,0)\big]\Bv_{\mathscr{C}^\ell(M)}\leqslant Ce^{\gamma t-\frac{c_2}{4t}-\sqrt{\gamma c_2}p},
	\end{split}
\end{equation}
where $\lv\cdot\rv_{\mathscr{C}^\ell(M)}$ is the $\mathscr{C}^\ell$ norm with respect to $x_0\in M$, and $c_1$ and $c_2$ are given in \eqref{fb25}.
\end{coro}
\begin{pro}
By \eqref{cb12}, \eqref{fa.4}, \eqref{fc10} and \eqref{4.1.52} by replacing $t$ with $t/p^2$, there is $C>0$ with
\begin{equation}\label{fg...1}
	\begin{split}
&\Bv\theta_{\sqrt{t}}^{-1}\tro_s^{\Lambda(T^*X)\otimes F_p}\Big[\exp{(-t{\mathcal{M}}^{F_p})}\lk x_0,x_0\rk-\exp{(-t{\mathscr{M}}_{x_0}^{F_p})}\lk0,0\rk\Big]\Bv_{\mathscr{C}^\ell(M)}\\
=&\Bv\theta_{\sqrt{\frac{p}{t}}}\tro_s^{\Lambda(T^*X)\otimes F_p}\Big[\exp{(-tp^{-2}{\mathcal{L}^{F_p}})}\lk x_0,x_0\rk-\exp{(-tp^{-2}{\mathscr{L}}_{x_0}^{F_p})}\lk0,0\rk\Big]\Bv_{\mathscr{C}^\ell(M)}\\
\leqslant&C(1+t^{-\dim S/2})p^{\dim S/2}\exp\big(c_1tp^{-2}-\frac{c_2p^2}{t}\big).
	\end{split}
\end{equation}
By the mean value inequality, we have
\begin{equation}\label{fg...2}
	\begin{split}
	c_1tp^{-2}-\frac{c_2p^2}{t}\leqslant \gamma t-\frac{c_2}{2t}-\Big(\frac{\gamma t}{2}+\frac{c_2p^2}{2t}\Big)\leqslant \gamma t-\frac{c_2}{2t}-\sqrt{\gamma c_2}p.
	\end{split}
\end{equation}
It is obvious that  $p^{\dim S/2}e^{-\sqrt{\gamma c_2}p}$ and $(1+t^{-\dim S})e^{-{c_2}/{4t}}$ are uniformly bounded for ${p\in{\mathbb{N}^*}}$ and $t>0$. Therefore, we get \eqref{fb.82} from \eqref{fg...1} and \eqref{fg...2}.\qed
\end{pro}

For $s\in \mathscr{C}^{\infty}\big(T_{x_0}X,\mathbb{E}_{p,x_0}\big), Z\in \mathbb{R}^m, u>0$,  we set $\lk K_us\rk\lk Z\rk=s\lk uZ\rk$ and 
\begin{equation}\label{4.2.2}
{\mathscr{N}}^{F_p}_{x_0}=K_{1/p}{\mathscr{M}}_{x_0}^{F_p}K_{p}.
\end{equation}
Let $\exp(-t{\mathscr{N}}^{F_p}_{x_0})\lk Z,Z'\rk$ be the kernel of $\exp(-t{\mathscr{N}}^{F_p}_{x_0})$ with respect to $d v_{TX}$, then
\begin{equation}\label{fc14}
\begin{split}
&\big(\exp(-t{\mathscr{N}}^{F_p}_{x_0})s\big)(Z)=\big(K_{1/p}\exp(-t{\mathscr{M}}_{x_0}^{F_p})K_ps\big)(Z)\\
&=\int_{\mathbb{R}^m}\exp(-t{\mathscr{M}}_{x_0}^{F_p})(p^{-1}Z,Z')s(pZ')\kappa(Z')dv_{TX}(Z')\\
&=p^{-m}\int_{\mathbb{R}^m}\exp(-t{\mathscr{M}}_{x_0}^{F_p})(p^{-1}Z,p^{-1}Z')\kappa(p^{-1}Z')s(Z')dv_{TX}(Z').
\end{split}
\end{equation}
By \eqref{4.1.48} and \eqref{fc14}, we get
\begin{equation}\label{fc15}
\exp(-t{\mathscr{N}}^{F_p}_{x_0})(0,0)=p^{-m}\exp(-t{\mathscr{M}}_{x_0}^{F_p})(0,0).
\end{equation}

We introduce another copy $\widehat{T_{x_0}X}$ of $T_{x_0}X$. We will add an extra hat when we refer to an element in $\widehat{T_{x_0}X}$. For $s>0$, put
\begin{equation}\label{aa2}
	c_s({e}_j)=\frac{1}{\sqrt{s}}{e^j}\wedge-\sqrt{s}{i}_{{e_j}}\ \ \ \ \ \widehat{c}_s({e}_j)=\frac{1}{\sqrt{s}}{e^j}\wedge+{\sqrt{s}}{i}_{{e}_j}.
\end{equation}
We denote by $\widehat{\mathscr{L}}^{F_p}_{x_0}$ the operator obtained from ${\mathscr{N}}^{F_p}_{x_0}$ by replacing $c(e_i), \widehat{c}(e_i)$ with $c_{1/p}(e_i)$ and $\widehat{c}_{1/p}(\widehat{e}_i)$ respectively. Set
\begin{equation}\label{fc16}
\widehat{\mathbb{F}}_{p}=\mathbb{R}[z]\hti\pi^*\Lambda\lk T_{}^*S\rk\hti\Lambda\lk T_{}^*X\rk\hti \Lambda(\widehat{T_{}^*X})\otimes F_{p},
\end{equation}
then $\widehat{\mathscr{L}}^{F_p}_{x_0}$ acts naturally on $\mathscr{C}^\infty(T_{x_0}X,\widehat{\mathbb{F}}_{p,x_0})$. We denote by $\exp(-t\widehat{\mathscr{L}}^{F_p}_{x_0})\lk Z,Z'\rk$ the smooth kernels of $\exp(-t\widehat{\mathscr{L}}^{F_p}_{x_0})$ with respect to $d v_{T_{}X}(Z')$. Put
\begin{equation}\label{fc17}
	\widehat{\mathbb{E}}_{p}=\mathbb{R}[z]\hti \pi^*\Lambda (T_{}^*S)\hti \mathrm{End}\big(\Lambda( T_{}^*X)\big)\hti\mathrm{End}\big(\Lambda( \widehat{T_{}^*X})\big)\otimes\mathrm{End}(F_p),
\end{equation}
and denote by $\pi\colon TX\times_M TX\to M$ the natural projection, then $\exp(-t\widehat{\mathscr{L}}^{F_p}_{x_0})\lk Z,Z'\rk$ is a smooth section of $\pi^*\widehat{\mathbb{E}}_p$ over $TX\times_MTX$.

For any $H\in \mathrm{End}\big(\Lambda(T^*X)\big)\hti\mathrm{End}\big(\Lambda(\widehat{T^*X})\big)$, it can be expanded uniquely in the following form:
\begin{equation}
	\begin{split}
H=\sum_{\begin{subarray}{c}
				1 \leqslant i_1<\cdots< i_p\leqslant m\\
				1 \leqslant i'_1<\cdots< i'_{p'}\leqslant m
		\end{subarray}}\sum_{\begin{subarray}{c}
		1 \leqslant j_1<\cdots< j_q\leqslant m\\
		1 \leqslant j'_1<\cdots< j'_{q'}\leqslant m
	\end{subarray}}Q^{j_1,\cdots,j_q,j'_1,\cdots,j'_{q'}}_{i_1,\cdots,i_p,i'_1,\cdots,i'_{p'}}e^{i_1}\wedge\cdots\wedge e^{i_p}\wedge i_{e_{j_1}}\cdots i_{\widehat{e}_{j_q}}\\
\wedge \widehat{e}^{i'_1}\wedge\cdots\wedge \widehat{e}^{i'_{p'}}\wedge i_{\widehat{e}_{j'_{1}}}\cdots i_{e_{j'_{q'}}},
	\end{split}
\end{equation}
and we set
\begin{equation}\label{fc57}
	[H]^{\text{max}}=Q_{1,\cdots,n,1,\cdots,n}e^{1}\wedge\cdots\wedge e^{n}\wedge\widehat{e}^{1}\wedge\cdots\wedge \widehat{e}^{n}.
\end{equation}
Now we define a linear map $\int^{\widehat{B}}\colon\Lambda(T^*X)\hti \Lambda(\widehat{T^*X})\to \Lambda(T^*X)$ called the \emph{Berezin integral}: first we assume that $\widehat{TX}$ is oriented by the form $\widehat{e}^1\wedge\cdots\wedge\widehat{e}^m$, then for $\alpha\in \Lambda(T^*X)$ and $\widehat{\beta}\in  \Lambda(\widehat{T^*X})$,
	\begin{equation}\label{ab1}
		\begin{split}
			\int^{\widehat{B}}\alpha\widehat{\beta}=0\ \  \text{if }\deg\widehat{\beta}<m,\ \ \ \int^{\widehat{B}}\alpha\widehat{e}^1\wedge\cdots\wedge\widehat{e}^m=\frac{(-1)^{m(m+1)/2}}{\pi^{m/2}}\alpha,
		\end{split}
	\end{equation}
	if $TX$ is not oriented, let $o(\widehat{TX})$ be the orientation line of $\widehat{TX}$, then $\int^{\widehat{B}}$ indeed defines a linear map from $\Lambda(T^*X)\hti \Lambda(\widehat{T^*X})$ to $\Lambda(T^*X)\hti o(\widehat{TX})$.

If $G\in \mathrm{End}\big(\Lambda(T^*X)\big)$, through \eqref{aa5}, we replace each $\widehat{c}(e)$ in $G$ by $\widehat{c}(\widehat{e})$ to obtain $\widehat{G}\in\mathrm{End}\big(\Lambda(T^*X)\big)\hti\mathrm{End}\big(\Lambda(\widehat{T^*X})\big)$.

By \cite[Proposition 4.9]{bz92}, among the monomials in $c({e}_i)$ and $\widehat{c}({e}_j)$, only the term $c({e}_1)\widehat{c}({e}_1)\cdots c({e}_m)\widehat{c}({e}_m)$ has a nonzero supertrace, and
	\begin{equation}\label{aa6}
		\trs^{\Lambda(T^*X)}[c({e}_1)\widehat{c}({e}_1)\cdots c({e}_m)\widehat{c}({e}_m)]=(-2)^m.
	\end{equation}
From \eqref{ab1} and \eqref{aa6}, we get
\begin{equation}\label{ab.4}
	\begin{split}
		\tro_s^{\Lambda(T^*X)}[G]\cdot e^{1}\wedge\cdots\wedge e^{n}=(4\pi)^{\frac{m}{2}} s^{m}\int^{\widehat{B}}[\widehat{G}]^{\mathrm{max}}.
	\end{split}
\end{equation}

\begin{prop}\label{4.6}
	For any ${p\in{\mathbb{N}^*}},t>0$ and $x_0\in M$, we have
	\begin{equation}\label{4.2.8}
		\mathrm{Tr}_s^{\Lambda(T^*X)\otimes F_p}\big[\exp(-t{\mathscr{M}}_{x_0}^{F_p})\lk0,0\rk\big]dv_{X}=(4\pi)^{\frac{m}{2}}\tro^{F_p}\int^{\widehat{B}}\big[\exp(-t\widehat{\mathscr{L}}^{F_p}_{x_0})(0,0)\big]^{\mathrm{max}}.
	\end{equation}
\end{prop}
\begin{pro}
By the definition of $\widehat{\mathscr{L}}^{F_p}_{x_0}$ and \eqref{ab.4}, we get
\begin{equation}\label{fc19}
	\begin{split}
\tro_s^{\Lambda(T^*X)\otimes F_p}\big[\exp(-t{\mathscr{N}}_{x_0}^{F_p})\lk 0,0\rk\big]dv_X(x_0)=(4\pi)^{\frac{m}{2}}p^{-m}\tro^{F_p}\int^{\widehat{B}}\big[\exp(-t\widehat{\mathscr{L}}^{F_p}_{x_0})(0,0)\big]^{\mathrm{max}},
	\end{split}
\end{equation}
then \eqref{4.2.8} follows immediately from \eqref{fc15} and \eqref{fc19}.\qed
\end{pro}

Let $\widehat{N}$ be the number operator of $\mathbb{R}[z]\hti\pi^*\Lambda\lk T^*S\rk\hti\Lambda\lk T^*X\rk\hti \Lambda(\widehat{T^*X})$ that acts by multiplication by the total degree of this exterior algebra. For any $a>0$, set
\begin{equation}\label{fc20}
\widehat{\theta}_a=a^{\widehat{N}}.
\end{equation} 

Recall the connection forms $\Gamma^{0}_Z$ and $\Gamma^{F_p,u}_Z$ at the beginning of \Cref{Fc}, $f_{\varepsilon}(Z)$ in \eqref{fb.76} and $g^{ij}$ as in \eqref{4.1.49}. Then $f_{\varepsilon}(Z)\Gamma^{0}_Z$ is a $1$-form evaluating in $\mathbb{R}[z]\hti\pi^*\Lambda(T_{\pi(x_0)}^*S)\hti\mathrm{End}(\Lambda(T_{x_0}^*X))$. By replacing $\widehat{c}(e_i)$ with $\widehat{c}(\widehat{e}_i)$ in $f_{\varepsilon}(Z)\Gamma^{0}_Z$, we get a $1$-form $\widehat{\Gamma}^{0}_Z$ on $T_{x_0}X$ taking values in $\mathbb{R}[z]\hti\pi^*\Lambda(T_{\pi(x_0)}^*S)\hti\mathrm{End}(\Lambda(T_{x_0}^*X))\hti\mathrm{End}(\Lambda(\widehat{T_{x_0}^*X}))$. Set
\begin{equation}
\widehat{\Gamma}^{F_p,u}_Z=f_{\varepsilon}(Z)\Gamma^{F_p,u}_Z.
\end{equation}

By \eqref{4.2.2} and the construction of $\widehat{\mathscr{L}}^{F_p}_{x_0}$, we have the following identity by evaluating the tensors on the right-hand side of the equation at $Z/p$:
\begin{equation}\label{fc.22}
	\begin{split}
		\widehat{\mathscr{L}}^{F_p}_{x_0,Z}=-g^{ij}&\Big\{\big(\pa_i+\frac{1}{p}\widehat{\theta}_{1/p}^{-1}\widehat{\Gamma}^{0}\lk\pa_i\rk\widehat{\theta}_{1/p}+p^{-1}\widehat{\Gamma}^{F_p,u}\lk\pa_i\rk\big)\\
		&\ \ \ \ \ \ \ \ \ \cdot\big(\pa_j+\frac{1}{p}\widehat{\theta}_{1/p}^{-1}\widehat{\Gamma}^{0}\lk\pa_j\rk\widehat{\theta}_{1/p}+p^{-1}\widehat{\Gamma}^{F_p,u}\lk\pa_j\rk\big)\\
		&\ \ \ \ \ \ \ \ +\frac{1}{p}\big(\nabla^{TX_0}_{\pa_i}\pa_j+\frac{1}{p}\widehat{\theta}_{1/p}^{-1}\widehat{\Gamma}^{0}\big(\nabla^{TX_0}_{\pa_i}\pa_j\big)\widehat{\theta}_{1/p}+p^{-1}\widehat{\Gamma}^{F_p,u}\big(\nabla^{TX_0}_{\pa_i}\pa_j\big)\big)\Big\}\\		+f_{\varepsilon}&\bigg\{\frac{r^X}{4p^2}-\frac{1}{8p^2}\big\langle R^{TX}(e_i,e_j)e_k,e_\ell \big\rangle c_{1/p}(e_i)c_{1/p}(e_j)\widehat{c}_{1/p}\lk \widehat{e}_k\rk\widehat{c}_{1/p}\lk \widehat{e}_\ell\rk\\
		&-\frac{1}{8p}\big\langle R^{TX}(f_\alpha^H,f_\beta^H\big)e_k,e_\ell \big\rangle f^\alpha f^\beta\widehat{c}_{1/p}\lk \widehat{e}_k\rk\widehat{c}_{1/p}\lk \widehat{e}_\ell\rk\\
		&-\frac{1}{4p^{3/2}}\big\langle R^{TX}\big(e_i,f_\alpha^H\big)e_k,e_\ell \big\rangle c_{1/p}(e_i)f^\alpha\widehat{c}_{1/p}\lk \widehat{e}_k\rk\widehat{c}_{1/p}\lk \widehat{e}_\ell\rk \\
		&-\frac{1}{2p}c_{1/p}(e_i)c_{1/p}(e_j)\frac{1}{4p}\omega\big(\nabla^{F_p},g^{F_p}\big)^2(e_i,e_j)\\
		&-\frac{1}{2}f^\alpha f^\beta\frac{1}{4p}\omega\big(\nabla^{F_p},g^{F_p}\big)^2\big(f_\alpha^H,f_\beta^H\big)\\
		&-\frac{1}{p^{1/2}}c_{1/p}(e_i)f^\alpha\frac{1}{4p}\omega\big(\nabla^{F_p},g^{F_p}\big)^2\big(e_i,f_\alpha^H\big)\\
		&+\frac{1}{4p^2}\widehat{\omega}\big(\nabla^{F_p},g^{F_p}\big)(e_i)^2+\frac{1}{2p}\widehat{c}_{1/p}(\widehat{e}_i)\widehat{c}_{1/p}(\widehat{e}_j)\frac{1}{4p}\widehat{\omega}\big(\nabla^{F_p},g^{F_p}\big)^2(e_i,e_j)\\
		&-\frac{1}{p^{1/2}}f^\alpha\widehat{c}_{1/p}(\widehat{e}_i)\frac{1}{2p}\nabla^{\widehat{TX}\otimes F_p,u}_{f^H_\alpha}\widehat{\omega}\big(\nabla^{F_p},g^{F_p}\big)(e_i)\\
		&-\frac{1}{p}{c}_{1/p}(e_i)\widehat{c}_{1/p}(\widehat{e}_j)\frac{1}{2p}\nabla^{\widehat{TX}\otimes F_p,u}_{e_i}\widehat{\omega}\big(\nabla^{F_p},g^{F_p}\big)(e_j)\\
		&-\frac{1}{p^{1/2}}zc_{1/p}(e_i)\frac{1}{2p}{\omega}\big(\nabla^{F_p},g^{F_p}\big)(e_i)-zf^\alpha\frac{1}{2p}{\omega}\big(\nabla^{F_p},g^{F_p}\big)\big(f^H_\alpha\big)\bigg\}.
	\end{split}
\end{equation}

Note that $\dim F_p$ grows as a polynomial of $p\in\mathbb{N}^*$, locally, we \emph{cannot} reduce $F_p$ to a fixed bundle, so it is nonsense to consider the limit operator $\lim_{p\to+\infty}\widehat{\mathscr{L}}_{x_0,Z}^{F_p}$ naively. In \cite[(9.148)]{bmz17}, Bismut-Ma-Zhang viewed $\widehat{\mathscr{L}}_{x_0,Z}^{F_p}$ as a differential operator with Toeplitz operator \eqref{cb17} coefficients and obtained the first order expansion of $\widehat{\mathscr{L}}_{x_0,Z}^{F_p}$ with respect to $p\in\mathbb{N}^*$ in the sense of Toeplitz operator. 

Now we should describe its full ``expansion". We do \emph{not} expand $\widehat{\mathscr{L}}_{x_0}^{F_p}$ directly in the sense of Toeplitz operator. Instead, we introduce a series of operators $\{\widehat{\mathscr{L}}^{F_p}_v\}_{0\leqslant v\leqslant 1/p}$ such that $\widehat{\mathscr{L}}^{F_p}_{1/p}=\widehat{\mathscr{L}}_{x_0}^{F_p}$, and we view the Taylor expansion $\sum_{i=0}^\ell\frac{1}{i!p^{i}}\frac{\pa^i\widehat{\mathscr{L}}^{F_p}_v}{\pa v^i}|_{v=0}$ as the $\ell$-th order ``expansion" of $\widehat{\mathscr{L}}^{F_p}_{x_0}$. The advantage of this approach is that $v$ is a smooth parameter, so we could use techniques from Ma-Marinescu \cite[\S\,4.1]{mm07}, and this is also convenient for the computation of second order expansion in \Cref{sH}.

Now we give the definition of operators $\{\widehat{\mathscr{L}}^{F_p}_v\}_{0\leqslant v\leqslant 1/p}$ acting on $\mathscr{C}^{\infty}_0(T_{x_0}X,\widehat{\mathbb{F}}_{p,x_0})$. By evaluating the tensors on the right-hand side of the equation at $vZ$, we take
\begin{equation}\label{4.2.9}
	\begin{split}
\widehat{\mathscr{L}}^{F_p}_v=-&g^{ij}\Big\{\big(\pa_i+{v}\widehat{\theta}_{v}^{-1}\widehat{\Gamma}^{0}_{{vZ}}\lk\pa_i\rk\widehat{\theta}_{v}+p^{-1}\widehat{\Gamma}^{F_p,u}_{{vZ}}\lk\pa_i\rk\big)
\\
&\ \ \ \ \ \ \ \ \ \ \cdot\big(\pa_j+{v}\widehat{\theta}_{v}^{-1}\widehat{\Gamma}^{0}_{{vZ}}\lk\pa_j\rk\widehat{\theta}_{v}+p^{-1}\widehat{\Gamma}^{F_p,u}_{{vZ}}\lk\pa_j\rk\big)\\
	&\ \ \ \ \ \ \ \ +{v}\big(\nabla^{TX_0}_{\pa_i}\pa_j+{v}\widehat{\theta}_{v}^{-1}\widehat{\Gamma}^{0}_{{vZ}}(\nabla^{TX_0}_{\pa_i}\pa_j)\widehat{\theta}_{v}+p^{-1}\widehat{\Gamma}^{F_p,u}_{{vZ}}(\nabla^{TX_0}_{\pa_i}\pa_j)\big)\Big\}\\
+&f_\varepsilon\bigg\{\frac{v^2}{4}r^X-\frac{v^2}{8}\big\langle R^{TX}(e_i,e_j)e_k,e_\ell \big\rangle c_{v}(e_i)c_{v}(e_j)\widehat{c}_{v}\lk \widehat{e}_k\rk\widehat{c}_{v}\lk \widehat{e}_\ell\rk\\
		&-\frac{v}{8}\big\langle R^{TX}(f_\alpha^H,f_\beta^H\big)e_k,e_\ell \big\rangle f^\alpha f^\beta\widehat{c}_{v}\lk \widehat{e}_k\rk\widehat{c}_{v}\lk \widehat{e}_\ell\rk\\
		&-\frac{v^{3/2}}{4}\big\langle R^{TX}\big(e_i,f_\alpha^H\big)e_k,e_\ell \big\rangle c_{v}(e_i)f^\alpha\widehat{c}_{v}\lk \widehat{e}_k\rk\widehat{c}_{v}\lk \widehat{e}_\ell\rk\\
		&-\frac{v}{2}c_{v}(e_i)c_{v}(e_j)\frac{1}{4p}\omega\big(\nabla^{F_p},g^{F_p}\big)^2(e_i,e_j)\\
		&-\frac{1}{2}f^\alpha f^\beta\frac{1}{4p}\omega\big(\nabla^{F_p},g^{F_p}\big)^2\big(f_\alpha^H,f_\beta^H\big)+\frac{1}{4p^2}\widehat{\omega}\big(\nabla^{F_p},g^{F_p}\big)(e_i)^2\\
		&-{v^{1/2}}c_{v}(e_i)f^\alpha\frac{1}{4p}\omega\big(\nabla^{F_p},g^{F_p}\big)^2\big(e_i,f_\alpha^H\big)\\
		&+\frac{v}{2}\widehat{c}_{v}(\widehat{e}_i)\widehat{c}_{v}(\widehat{e}_j)\frac{1}{4p}\widehat{\omega}\big(\nabla^{F_p},g^{F_p}\big)^2(e_i,e_j)\\
		&-{v^{1/2}}f^\alpha\widehat{c}_{v}(\widehat{e}_i)\frac{1}{2p}\nabla^{\widehat{TX}\otimes F_p,u}_{f^H_\alpha}\widehat{\omega}\big(\nabla^{F_p},g^{F_p}\big)(e_i)\\
		&-{v}{c}_{v}(e_i)\widehat{c}_{v}(\widehat{e}_j)\frac{1}{2p}\nabla^{\widehat{TX}\otimes F_p,u}_{e_i}\widehat{\omega}\big(\nabla^{F_p},g^{F_p}\big)(e_j)\\
		&-v^{1/2}zc_{v}(e_i)\frac{1}{2p}{\omega}\big(\nabla^{F_p},g^{F_p}\big)(e_i)-zf^\alpha\frac{1}{2p}{\omega}\big(\nabla^{F_p},g^{F_p}\big)\big(f^H_\alpha\big)\Big\}.
	\end{split}
\end{equation}
Combing \eqref{fc.22} and \eqref{4.2.9}, we clearly have $\widehat{\mathscr{L}}^{F_p}_{1/p}=\widehat{\mathscr{L}}^{F_p}_{x_0}$. Now we discuss the Taylor expansion of $\widehat{\mathscr{L}}^{F_p}_v$ with respect to $v$. Set
\begin{equation}\label{fc24}
	\begin{split}
		R^{\widehat{TX}}&=\big\langle R^{{TX}}{e}_i,{e}_j\big\rangle\wedge \widehat{e}^j\wedge{i}_{\widehat{e}_i},\\ \widehat{R}^{TX}&=\frac{1}{2}\big\langle{R}^{TX}(e_k,e_\ell){e}_i,{e}_j\big\rangle \widehat{e}^k\wedge \widehat{e}^\ell\wedge {e}^j\wedge{i}_{{e}_i}.
	\end{split}
\end{equation}
\begin{theo}\label{Fc13}
For $r\in\mathbb{N}, p\in\mathbb{N}^*,0\leqslant v\leqslant1/p$, we can write
\begin{equation}\label{fc26}
\frac{\pa^r\widehat{\mathscr{L}}^{F_p}_v}{\pa v^r}=\mathscr{X}_{v,ij}^{F_p,(r)}(Z)\pa_i\pa_j+\mathscr{Y}_{v,i}^{F_p,(r)}(Z)\pa_i+ \mathscr{Z}_v^{F_p,(r)}(Z),\ \ \text{for\ } \mathscr{X}_{v,ij}^{F_p,(r)}(Z)=-\frac{\pa^rg^{ij}_{vZ}}{\pa v^r}
	\end{equation}
and $\mathscr{Y}_{v,i}^{F_p,(r)}(Z),\mathscr{Z}_v^{F_p,(r)}(Z)\in\mathscr{C}^\infty(T_{x_0}X,\widehat{\mathbb{E}}_{p,x_0})$ with the following properties:
	\begin{enumerate}
		\item  The function $\mathscr{X}_{v,ij}^{F_p,(r)}$ is a constant outside a compact set of $T_{x_0}X$.  Also, if $v\neq 0$, $\mathscr{Y}_{v,i}^{F_p,(r)}(Z)$ and $ \mathscr{Z}_v^{F_p,(r)}(Z)$ have compact supports.
	\item	There is $C>0$ such that for $p\in\mathbb{N}^*$ and $0\leqslant v\leqslant 1/p$, the operator norm $\lV\cdot\rV_{\widehat{\mathbb{E}}_{p,x_0}}$ of $\mathscr{X}_{v,ij}^{F_p,(r)}(Z),\mathscr{Y}_{v,i}^{F_p,(r)}(Z)$ and $\mathscr{Z}_v^{F_p,(r)}(Z)$ are dominated by $C(1+\vert Z\vert^r)$.
	
\item The operator $\mathscr{X}_{0,ij}^{F_p,(r)}(Z)$ (resp. $\mathscr{Y}_{0,i}^{F_p,(r)}(Z),\mathscr{Z}_0^{F_p,(r)}(Z)$) is polynomial in $Z$ (resp. with coefficients in $\widehat{\mathbb{E}}_{p,x_0}$). Moreover, $\mathscr{X}_{0,ij}^{F_p,(r)}(Z)$ is a homogeneous polynomial of $Z$ of degree $r$, the degree in $Z$ of $\mathscr{Y}_{0,i}^{F_p,(r)}(Z)$ ( $\mathscr{Z}_0^{F_p,(r)}(Z)$) is no more than $r$.

\item	In particular,
		\begin{equation}\label{fc28}
	\begin{split}
\widehat{\mathscr{L}}^{F_p}_0=&-{\Delta}^{T_{x_0}X}-\frac{1}{4}\big\langle R^{TX}e_i,e_j \big\rangle_{x_0}\widehat{e}^i\wedge\widehat{e}^j-\frac{1}{4p}\omega\big(\nabla^{F_p},g^{F_p}\big)^2_{x_0}\\
		&+\frac{1}{4p^2}\widehat{\omega}\big(\nabla^{F_p},g^{F_p}\big)_{x_0}(e_i)^2+\frac{1}{4p}\widehat{\omega}\big(\nabla^{F_p},g^{F_p}\big)^2_{x_0}\\
		&-\frac{1}{2p}\nabla^{\widehat{TX}\otimes F_p,u}\widehat{\omega}\big(\nabla^{F_p},g^{F_p}\big)_{x_0}-\frac{1}{2p}z\omega\big(\nabla^{F_p},g^{F_p}\big)_{x_0},\\
		\frac{\pa\widehat{\mathscr{L}}^{F_p}_v}{\pa v}\Big|_{v=0}=&-\frac{Z_j}{4}\Big\{\big\langle R^{TX}_{x_0}\pa_i,\pa_j\big\rangle-\big\langle \widehat{R}^{TX}_{x_0}\pa_i,\pa_j\big\rangle-\frac{1}{2p}\omega\big(\nabla^{F_p},g^{F_p}\big)^2_{x_0}\big(\pa_i,\pa_j\big)\Big\}\pa_i\\
		&+\mathscr{Q}^{F_p}(Z)+\mathscr{Q}^{F_p}{'},
	\end{split}
		\end{equation}
	where $\mathscr{Q}^{F_p}(Z)$ (resp. $\mathscr{Q}^{F_p}{'}$) is a homogeneous polynomials in $Z$ of degree $1$ (resp. $0$) with coefficients in $\widehat{\mathbb{E}}_{p,x_0}$, and
	\begin{equation}\label{fc29}
		\begin{split}
\mathscr{Q}^{F_p}{'}=&\frac{1}{2}R^{\widehat{TX}}-\frac{1}{2}\widehat{R}^{TX}-\frac{1}{4}\big\langle {R}^{TX}(f_\alpha^H,e_i){e}_k,{e}_l\big\rangle \widehat{e}^k\wedge\widehat{e}^l\wedge f^\alpha\wedge{i}_{e_i}\\
&+e^i\wedge{i}_{e_j}\frac{1}{4p}\omega_p^2(e_i,e_j)+f^\alpha\wedge{i}_{e_i}\frac{1}{4p}\omega_p^2\big(f_\alpha^H,e_i\big)+\widehat{e}^i\wedge{i}_{\widehat{e}_j}\frac{1}{4p}\widehat{\omega}_p^2(\widehat{e}_i,\widehat{e}_j)\\
&-f^\alpha\wedge{i}_{\widehat{e}_i}\frac{1}{2p}\nabla^{\widehat{TX}\otimes F_p,u}_{f^H_\alpha}\widehat{\omega}_p(e_i)-{e}^i\wedge{i}_{\widehat{e}_j}\frac{1}{2p}\nabla^{\widehat{TX}\otimes F_p,u}_{e_i}\widehat{\omega}_p(e_j)\\
&-\widehat{e}^j\wedge{i}_{e_i}\frac{1}{2p}\nabla^{\widehat{TX}\otimes F_p,u}_{e_i}\widehat{\omega}_p(e_j)+z{i}_{e_i}\frac{1}{2p}\omega_p(e_i).
		\end{split}
	\end{equation}
	\end{enumerate}
\end{theo}
\begin{pro}
As in the proof of Lemma \ref{Fb4}, by Theorem \ref{Db4}, \eqref{bb3} and \eqref{fb7}, we get
\begin{equation}\label{fc31}
\begin{split}
	p^{-1}\widehat{\Gamma}^{F_p,u}_{vZ}\big(\pa_i\big)\sim \frac{vZ_j}{2}\cdot\frac{1}{4p}\omega\big(\nabla^{F_p},g^{F_p}\big)^2_{x_0}\big(\pa_i,\pa_j\big)+\mathcal{O}(v^2).
\end{split}
\end{equation}
By \eqref{aa2}, we have
\begin{equation}\label{4.2.27}
{v^{1/2}}{c}_{v}(e_i)=e^i\wedge-{v}{i}_{e_i},\ \ \  v^{1/2}\widehat{c}_{v}(\widehat{e}_i)=\widehat{e}^i\wedge+{v}{i}_{\widehat{e}_i}.
\end{equation}
By \eqref{4.2.27}, similar to \eqref{fc31},we have
\begin{equation}\label{fc32}
	\begin{split}
		\widehat{\theta}_{v}^{-1}(\widehat{\Gamma}^0_{vZ})\widehat{\theta}_{v}\big(\pa_i\big)\sim-\frac{1}{4}Z_j
		\Big(\big\langle R^{TX}_{x_0}(e_k,e_\ell)\pa_i,\pa_j\big\rangle e^k\wedge e^\ell-\langle R^{TX}_{x_0}(e_k,e_\ell)\pa_i,\pa_j\big\rangle\widehat{e}^k\wedge\widehat{e}^\ell\\ +\big\langle R^{TX}_{x_0}(f_\alpha,e_\ell)\pa_i,\pa_j\big\rangle f^\alpha\wedge e^\ell+\frac{1}{2}\big\langle R_{x_0}^{TX}(f_\alpha,f_\beta)\pa_i,\pa_j\big\rangle f^\alpha\wedge f^\beta\Big)+\mathcal{O}(v).
	\end{split}
\end{equation}
By \eqref{4.2.9}, \eqref{fc31}, \eqref{4.2.27} and \eqref{fc32}, we can write
\begin{equation}\label{fc35}
\widehat{\mathscr{L}}^{F_p}_v=\mathscr{X}^{F_p}_{v,ij}\pa_i\pa_j+\mathscr{Y}^{F_p}_{v,i}\pa_i+\mathscr{Z}^{F_p}_{v},\ \ \text{where\ }\mathscr{X}^{F_p}_{v,ij}(Z)=-g^{ij}_{vZ},
\end{equation}
and each of $\{\mathscr{Y}^{F_p}_{v,i},\mathscr{Z}^{F_p}_v\}$ is the sum of terms in the form of ${v^k}{\Gamma}_{{vZ}}$, where $k\in\mathbb{N}$, ${\Gamma}\in\mathscr{C}_0^\infty(T_{x_0}X,\widehat{\mathbb{E}}_{p,x_0})$ is a smooth family of Toeplitz operators (see \cref{Dc}) and has support in $B^{T_{x_0}X}(0,{\varepsilon})$, so its norm $\lv\cdot\rv_{\mathscr{C}(T_{x_0}X,\widehat{\mathbb{E}}_{p,x_0})}$ is bounded uniformly for $p\in\mathbb{N}^*$ as well as all its derivatives. For $r\in\mathbb{N}$, we have
\begin{equation}\label{fc36}
\frac{d^r}{dv^r}{v^k}{\Gamma}_{{vZ}}=\sum_{\substack{r_1+r_2=r\\ r_1\leqslant k}}\sum_{\vert\alpha\vert=r_2}C_\alpha{Z^\alpha}{v^{k-r_1}}\Big(\frac{\pa^\alpha {\Gamma}}{\pa Z^\alpha}\Big)_{{vZ}},
\end{equation} 
where $\alpha$ is a multi-index. By \eqref{fc36}, $\frac{\pa^{r}\widehat{\mathscr{L}}^{F_p}_v}{\pa v^{r}}$ is the sum of terms in the form of 
\begin{equation}\label{fc37}
\sum_{\vert\alpha\vert\leqslant r}C_\alpha v^kZ^{\alpha}\Gamma_{\alpha,{vZ}}\pa_i\pa_j+C'_\alpha v^{k'}Z^{\alpha}\Gamma'_{\alpha,{vZ}}\pa_i+C''_\alpha v^{k''}Z^{\alpha}\Gamma''_{\alpha,{vZ}},
\end{equation}
where $\Gamma_\alpha,\Gamma_\alpha',\Gamma_\alpha''\in\mathscr{C}_0^\infty(T_{x_0}X,\widehat{\mathbb{E}}_{p,x_0})$ verify similar property as $\Gamma$. By \eqref{fc37}, we get the first three statements.

By \eqref{4.2.9}, \eqref{fc31} and \eqref{fc32}, we get the first identity in \eqref{fc28}, also, together with the property of geodesic coordinates $(\pa_ig^{jk}\big)(0)=0$, we obtain $\mathscr{X}^{F_p}_{ij,1}(0)$ and $\mathscr{Y}^{F_p}_{i,1}(0)$ in the second identity of \eqref{fc28}. By \eqref{4.2.27}, the contribution of the first four terms in \eqref{4.2.9} to $\mathscr{Q}^{F_p}{'}$ is given by
\begin{equation}\label{4.2.28}
	\begin{split}
	&-\frac{1}{4}\big\langle R^{TX}(e_i,e_j)e_k,e_\ell \big\rangle e^i\wedge e^j\wedge\widehat{e}^k\wedge{i}_{\widehat{e}_\ell}+\frac{1}{4}\big\langle R^{TX}(e_i,e_j)e_k,e_\ell \big\rangle \widehat{e}^k\wedge\widehat{e}^\ell\wedge e_i\wedge{i}_{e_j}\\
		&-\frac{1}{4}\big\langle R^{TX}(f_\alpha^H,f_\beta^H\big)e_k,e_\ell \big\rangle f^\alpha f^\beta\widehat{e}^k\wedge{i}_{\widehat{e}_\ell}\\
		&-\frac{1}{4}\big\langle R^{TX}\big(e_i,f_\alpha^H\big)e_k,e_\ell \big\rangle e^i\wedge f^\alpha\wedge\widehat{e}^k\wedge{i}_{\widehat{e}_\ell}+\frac{1}{4}\big\langle R^{TX}\big(e_i,f_\alpha^H\big)e_k,e_\ell \big\rangle \widehat{e}^k\wedge f^\alpha{i}_{e_i},\\
	\end{split}
\end{equation}
which yields the first three terms in \eqref{fc29}. The other terms in \eqref{fc29} follows from \eqref{4.2.9} and \eqref{4.2.27}. We finish the proof.\qed
\end{pro}

\subsection{Sobolev spaces with weights and estimates on the resolvents}\label{Fd}

In this subsection, we follow the strategy of Bismut-Lebeau \cite[\S\,11 k)-l)]{bisleb91} with some modifications. In particular, we should choose the weight in Theorem \ref{4.11} and the path of integral in Theorem \ref{4.13} carefully to ensure the ``\emph{positivity}" of  $\widehat{\mathscr{L}}_v^{F_p}$: by \eqref{fc28}, we see that $\widehat{\mathscr{L}}_v^{F_p}=\widehat{\mathscr{L}}^{F_p}_0+\mathcal{O}(v)$, and $\widehat{\mathscr{L}}^{F_p}_0$ is the sum of a \emph{non-negative} part $-{\Delta}^{T_{x_0}X}+\frac{1}{4p}\vert\widehat{\omega}(\nabla^{F_p},g^{F_p})\vert^2_{x_0}$ and a \emph{nilpotent} part. So, we want to prove that $\widehat{\mathscr{L}}^{F_p}_v$  is ``nearly" a non-negative operator when $v$ is small: for any ${\gamma}>0$, there is $C>0$ such that for $t>0$ large enough, we have $\Vert\exp(-t\widehat{\mathscr{L}}_v^{F_p})(0,0)\Vert\leqslant C\exp({\gamma} t)$. To eliminate the effect of the nilpotent part in $\widehat{\mathscr{L}}^{F_p}_v$, we introduce Sobolev norms with weights. 

From now on, we will fix a $\gamma>0$. Note that $\widehat{\mathbb{F}}_p$ and $\widehat{\theta}_\mu$ are defined in \eqref{fc16} and \eqref{fc20}.
\begin{defi}\label{Fd1}
For $0<\mu<1$, $k\in\mathbb{N}$ and $s,s'\in\mathscr{C}_0^\infty(T_{x_0}X,\widehat{\mathbb{F}}_{p,x_0})$, set
\begin{equation}\label{fd1}
	\begin{split}
\li s,s'\ri_{0,p,\mu}=\int_{\mathbb{R}^{n}}\big\langle\widehat{\theta}_\mu s,\widehat{\theta}_\mu s'\big\rangle dv_{TX},\ \ \ 
\li s,s'\ri_{k,p,\mu}=\sum_{\lv\alpha\rv\leqslant k}\big\langle{\pa^\alpha s},{\pa^\alpha s'}\big\rangle_{0,p,\mu},
	\end{split}
\end{equation}
and let $\lV\cdot\rV_{0,p,\mu}$ and $\lV\cdot\rV_{k,p,\mu}$ be the corresponding norms. Let $H_{k,p,\mu}$ be the Sobolev space which is the completion of $\mathscr{C}_0^\infty(T_{x_0}X,\widehat{\mathbb{F}}_{p,x_0})$ with respect to $\lV\cdot\rV_{k,p,\mu}$.
	Let $H_{-k,p,\mu}$ be the Sobolev space of negative order with the norm: for $s\in\mathscr{C}_0^\infty(T_{x_0}X,\widehat{\mathbb{F}}_{p,x_0})$,
	\begin{equation}
		\lV s\rV_{-k,p,,\mu}=\sup_{0\neq s'\in H_{k,p,\mu}}\frac{\bv\li s,s'\ri_{0,p,\mu}\bv}{\lV s'\rV_{k,p,\mu}}.
	\end{equation} 
For $k_1,k_2\in\mathbb{Z}$, on the space of bounded linear map from $H_{k_1,p,\mu}$ to $H_{k_2,p,\mu}$, we denote by $\lV \cdot\rV_{p,\mu}^{k_1,k_2}$ the operator norm.
\end{defi}

By \eqref{fd1}, it is clear that
\begin{equation}\label{fd.0}
\lV e^i\wedge\rV_{p,\mu}^{0,0}\leqslant \mu,\ \ \ \lV i_{e_i}\rV_{p,\mu}^{0,0}\leqslant \frac{1}{\mu},
\end{equation}
and \eqref{fd.0} is tailor-made for the next theorem.
\begin{theo}[Elliptic estimations]\label{4.11}
There are $\mu,{\varepsilon}>0$ such that we have $C_i>0$ for $1\leqslant i\leqslant 4$ in which $C_2<{\gamma}/2$ and $p_0\in\mathbb{N}^*$ such that for any $p\geqslant p_0$, $0\leqslant v\leqslant 1/p$ and $s,s'\in\mathscr{C}_0^\infty(T_{x_0}X,\widehat{\mathbb{F}}_{p,x_0})$, we have
	\begin{equation}\label{4.3.1.4}
		\begin{split}
			&\mathrm{Re}\big\langle \widehat{\mathscr{L}}^{F_p}_vs,s \big\rangle_{0,p,\mu}\geqslant C_1\lV\nabla s\rV_{0,p,\mu}^2-C_2\Vert s\Vert_{0,p,\mu}^2,\\
			&\lv\mathrm{Im}\big\langle \widehat{\mathscr{L}}^{F_p}_vs,s \big\rangle_{0,p,\mu}\rv\leqslant C_3\Vert s\Vert_{1,p,\mu}\Vert s\Vert_{0,p,\mu},\\
			&\left\vert\big\langle \widehat{\mathscr{L}}^{F_p}_vs,s' \big\rangle_{0,p,\mu}\right\vert\leqslant C_4\Vert s\Vert_{1,p,\mu}\cdot\Vert s'\Vert_{1,p,\mu}.
		\end{split}
	\end{equation}
\end{theo}

\begin{pro}
We now carefully discuss the structure of $\widehat{\mathscr{L}}^{F_p}_v$. First, we focus on a single term
\begin{equation}
\big(\pa_i+{v}\widehat{\theta}_{v}^{-1}\widehat{\Gamma}_{vZ}^{0}\lk\pa_i\rk\widehat{\theta}_{v}+p^{-1}\widehat{\Gamma}^{F_p,u}_{vZ}\lk\pa_i\rk\big)
\end{equation}
as in \eqref{4.2.9}. Let $\widehat{R}^0$ be the operator obtained by replacing $\widehat{c}(e_j)$ with $\widehat{c}(\widehat{e}_j)$ on the left hand side of \eqref{bj.6}. By \eqref{bj.6} and \eqref{fb7}, we get
\begin{equation}\label{fd5}
	\begin{split}
		{v}\widehat{\theta}_{v}^{-1}\widehat{\Gamma}_{vZ}^{0}\lk\pa_i\rk\widehat{\theta}_{v}&=-\frac{v}{2}\int_0^1f_\varepsilon({vZ}){vZ_j}\big(\widehat{\theta}_{v}^{-1}\widehat{R}^0_{{svZ}}\widehat{\theta}_{v}\big)(\pa_i,\pa_j)ds,\\
		p^{-1}\widehat{\Gamma}^{F_p,u}_{{vZ}}(\pa_i)&=-\int_{0}^1f_\varepsilon({vZ}){vZ_j}\frac{1}{4p}\omega\big(\nabla^{F_p},g^{F_p}\big)^2_{svZ}(\pa_i,\pa_j)ds.
	\end{split}
\end{equation}
Here we use the integral form rather than the Taylor expansion to get an estimation uniformly for ${p\in{\mathbb{N}^*}}$. By Theorem \ref{Db4}, the operator norm of $\frac{1}{4p}\omega(\nabla^{F_p},g^{F_p})^2$ is uniformly bounded for $p\in\mathbb{N^*}$. Since the support of $f_{\varepsilon}(Z)$ is contained in $B^{T_{x_0}X}(0,2{\varepsilon})$, by Theorem \ref{Db4} and \eqref{fd5}, there is $C>0$ such that for any $\varepsilon$ satisfies \eqref{fb2}, ${p\in{\mathbb{N}^*}}, 0\leqslant v\leqslant 1/p$ and $Z\in T_{x_0}X$, we have
\begin{equation}\label{fd.6}
	{v}\Vert\widehat{\theta}_{v}^{-1}\widehat{\Gamma}_{vZ}^{0}\lk\pa_i\rk\widehat{\theta}_{v}\Vert_{\widehat{\mathbb{E}}_{p,x_0}}\leqslant Cv,\ \ \ \big\Vert p^{-1}\widehat{\Gamma}^{F_p,u}_{{vZ}}(\pa_i)\big\Vert_{\widehat{\mathbb{E}}_{p,x_0}}\leqslant C{\varepsilon}.
\end{equation}
Consider the expansion in \eqref{fc35}
\begin{equation}\label{fd.11}
	\begin{split}
\widehat{\mathscr{L}}^{F_p}_v&=\mathscr{X}^{F_p}_{v,ij}(Z)\pa_i\pa_j+\mathscr{Y}^{F_p}_{v,i}(Z)\pa_i+ \mathscr{Z}^{F_p}_v(Z)\\
&=\pa_i\mathscr{X}^{F_p}_{v,ij}(Z)\pa_j+\big(\mathscr{Y}^{F_p}_{v,i}(Z)-\pa_j(\mathscr{X}^{F_p}_{v,ji})(Z)\big)\pa_i+\mathscr{Z}^{F_p}_v(Z).
	\end{split}
\end{equation}
By \eqref{4.2.9}, \eqref{4.2.27} \eqref{fd5} and \eqref{fd.6}, we see that $\mathscr{Z}^{F_p}_v(Z)$ can be separated into three parts $\mathscr{Z}^{F_p}_v(Z)=\mathscr{Z}^{F_p}_{v}{'}(Z)+\mathscr{Z}_{v}^{F_p}{''}(Z)+\mathscr{Z}_{v}^{F_p}{'''}(Z)$, where $\mathscr{Z}_{v}^{F_p}{'}(Z)$ is a non-negative operator, $\mathscr{Z}_{v}^{F_p}{''}(Z)$ adds the degree of exterior algebra in $\widehat{\mathbb{F}}_{p,x_0}$ and there is $C>0$ such that 
\begin{equation}\label{fd.10}
	\Vert\mathscr{Z}_{v}^{F_p}{''}(Z)\Vert_{\widehat{\mathbb{E}}_{p,x_0}}\leqslant C,\ \ \Vert\mathscr{Z}_{v}^{F_p}{'''}(Z)\Vert_{\widehat{\mathbb{E}}_{p,x_0}}\leqslant Cv.
\end{equation}
By \eqref{fb.73} and the positivity of the matrix $\{g^{ij}_{vZ}\}$, there is $C>0$ such that
\begin{equation}\label{fd.13}
 \big\langle \big(\pa_i\mathscr{X}^{F_p}_{v,ij}(Z)\pa_j\big)s,s \big\rangle_{0,p,\mu}\geqslant C\lV\nabla s\rV_{0,p,\mu}^2.
\end{equation}
For some $C>0$, we denote $C_{\varepsilon,\mu,v}=C(\varepsilon+\mu+\frac{v}{\varepsilon}+\frac{v}{\mu})$. Note that if we let $\pa_i$ act on the left hand side of \eqref{fd5}, we get an extra $\frac{v}{\varepsilon}$. This together with \eqref{fd.0} and \eqref{fd.6} implies
\begin{equation}\label{fd.15}
	\begin{split}
		\bv\big\langle \big(\mathscr{Y}^{F_p}_{v,i}-\pa_j(\mathscr{X}^{F_p}_{v,ji})\big)\pa_is,s \big\rangle_{0,p,\mu}\bv\leqslant C_{\varepsilon,\mu,v}\Vert s\Vert_{1,p,\mu}\Vert s\Vert_{0,p,\mu}.
	\end{split}
\end{equation}
Moreover, by \eqref{fd.0}, \eqref{fd.6} and \eqref{fd.10}, we have
\begin{equation}\label{fd.16}
	\begin{split}
		\Bv\big\langle\big(\mathscr{Z}_{v}^{F_p}{''}+\mathscr{Z}_{v}^{F_p}{'''}\big)s,s \big\rangle_{0,p,\mu}\Bv\leqslant C_{\varepsilon,\mu,v} \Vert s\Vert_{0,p,\mu}^2.
	\end{split}
\end{equation}
By \eqref{fd.13}, \eqref{fd.15} and \eqref{fd.16}, we have
\begin{equation}\label{fd.17}
	\begin{split}
		\mathrm{Re}\big\langle\widehat{\mathscr{L}}^{F_p}_vs,s \big\rangle_{0,p,\mu}&\geqslant C\lV\nabla s\rV_{0,p,\mu}^2-C_{\varepsilon,\mu,v}\Big(\Vert \nabla s\Vert_{0,p,\mu}\Vert s\Vert_{0,p,\mu}+\lV s\rV_{0,p,\mu}^2\Big)\\
		&\geqslant \Big(C-C_{\varepsilon,\mu,v}\Big)\lV\nabla s\rV_{0,p,\mu}^2-C_{\varepsilon,\mu,v}\lV s\rV_{0,p,\mu}^2,\\
		\lv\mathrm{Im}\big\langle\widehat{\mathscr{L}}^{F_p}_vs,s \big\rangle_{0,p,\mu}\rv&\leqslant C_{\varepsilon,\mu,v}\Big(\Vert \nabla s\Vert_{0,p,\mu}\Vert s\Vert_{0,p,\mu}+\lV s\rV_{0,p,\mu}^2\Big).
	\end{split}
\end{equation}
By \eqref{fd.17}, we get the first two inequalities in \eqref{4.3.1.4}, and the third one is obvious.\qed
\end{pro}
\begin{remark}
By \eqref{fd5}, there is $C>0$ such that for any $h$ satisfies \eqref{fb2}, ${p\in{\mathbb{N}^*}}, 0\leqslant v\leqslant 1/p$ and $Z\in T_{x_0}X$, we have $p^{-1}\bV\widehat{\Gamma}^{F_p,u}_{{vZ}}(\pa_i)\bV_{\mathrm{End}(F_{p,x_0})}\leqslant Cv\vert Z\vert$. Even if we get an $\mathcal{O}(v)$, it is not uniformly bounded for $Z\in T_{x_0}X$. This explains the necessity of \eqref{fd.6}.
\end{remark}

In the rest of this section, we will always fix a couple of $(\mu,{\varepsilon})$ satisfying Theorem \ref{4.11}.

\begin{prop}\label{4.12}
	For $k\in\mathbb{N}$, there exist $C>0$ and $p_0\in\mathbb{N}^*$ such that for $p\geqslant p_0, 0\leqslant v\leqslant 1/p$ and $Q_1,\cdots Q_k\in\{{\pa_i},Z_i\}_{i=1}^{m}$, we have
	\begin{equation}\label{4.3.1.7}
	\Bv\big\langle\big[Q_{1},[Q_{2},\cdots[Q_{k},\widehat{\mathscr{L}}^{F_p}_v]\cdots]\big]s,s'\big\rangle_{0,p,\mu}\Bv\leqslant C\Vert s\Vert_{1,p,\mu}\Vert s'\Vert_{1,p,\mu}.
	\end{equation}
\end{prop}

\begin{pro}
	By \eqref{fc35} and \eqref{fc36}, for $i,j\in\mathbb{N}$, we have
	\begin{equation}\label{fd9}
	[Z_j,\pa_i]={\delta}_{ij},\ \ \ \pa_i\big({\Gamma}_{{vZ}}\big)={v}(\pa_i{\Gamma})_{{vZ}},
	\end{equation}
where ${\Gamma}_{Z}\in \mathscr{C}_0^\infty(T_{x_0}X,\widehat{\mathbb{E}}_{p,x_0})$. By \eqref{fc35} and \eqref{fd9}, $\big[Q_{1},[Q_{2},\cdots[Q_{k},\widehat{\mathscr{L}}^{F_p}_v]\cdots]\big]$ has the same structure as $\widehat{\mathscr{L}}^{F_p}_v$ in \eqref{fc35}. Hence we easily get \eqref{fd9} by \eqref{4.3.1.4}.\qed
\end{pro}

For $d_1,d_2>0$, set
\begin{equation}\label{fe1}
	V_{d_1,d_2}=\left\{\lambda\in\mathbb{C}\mid \text{Re}(\lambda)\leqslant d_1\text{Im}(\lambda)^2-d_2\right\},
\end{equation}
and let $\pa V_{d_1,d_2}$ be its boundary with counterclockwise orientation. Recall the norm $\lV\cdot\rV_{p,\mu}^{k_1,k_2}$ given in Definition \ref{Fd1}. 
\begin{theo}[Elliptic regularity]\label{4.13}
There exist $d,C>0$ and $p_0\in\mathbb{N}^*$ such that for any $p\geqslant p_0$, $0\leqslant v\leqslant p^{-1}$ and $\lambda\in V_{d,{\gamma}/{2}}$, the resolvent  $\blk\lambda-\widehat{\mathscr{L}}^{F_p}_v\brk^{-1}$ exists, and 
	\begin{align}\label{fe2}
		\begin{split}
	\bV\blk\lambda-\widehat{\mathscr{L}}^{F_p}_v\brk^{-1}\bV^{0,0}_{p,\mu}\leqslant C,\ \ \ \bV\blk\lambda-\widehat{\mathscr{L}}^{F_p}_v\brk^{-1}\bV^{-1,1}_{p,\mu}\leqslant C\blk1+\lv\lambda\rv^2\brk.		\end{split}
	\end{align}
\end{theo}

\begin{pro}
Since $\gamma$ could be arbitrarily small, \eqref{fe2} almost indeed ensures the ``non-negativity'' of $\mathrm{Spec}(\widehat{\mathscr{L}}^{F_p}_v)$.
	
	For $a,b\in\mathbb{R}$ and $\lambda=a+ib\in \mathbb{C}$, we have
	\begin{equation}\label{4.4.4}
		\begin{split}
			&\bv\big\langle\blk\lambda-\widehat{\mathscr{L}}^{F_p}_v\brk s,s\big\rangle_{0,p,\mu}\bv\geqslant\sup\Big\{\text{Re}\bli \widehat{\mathscr{L}}^{F_p}_vs,s\bri_{0,p,\mu}-a\lV s\rV_{0,p,\mu}^2,\Bv\text{Im}\bli \widehat{\mathscr{L}}^{F_p}_v s,s\bri_{0,p,\mu}-b\lV s\rV_{0,p,\mu}^2\Bv\Big\}\\
			&\geqslant \sup\Big\{C_1\lV s\rV_{1,p,\mu}^2-(a+C_1+C_2)\lV s\rV_{0,p,\mu}^2,\lv b\rv\lV s\rV_{0,p,\mu}^2-C_3\lV s\rV_{0,p,\mu}\lV s\rV_{1,p,\mu}\Big\}.
		\end{split}
	\end{equation}
It is obvious that $\lV s\rV_{0,p,\mu}\leqslant \lV s\rV_{1,p,\mu}$, then we deduce from \eqref{4.4.4} that
	\begin{equation}\label{4.4.6}
		\bv\big\langle\blk\lambda-\widehat{\mathscr{L}}^{F_p}_v\brk s,s\big\rangle_{0,p,\mu}\bv\geqslant C(\lambda) \lV s\rV_{0,p,\mu}^2,
	\end{equation}
where $C(\lambda)$ is given by
\begin{equation}\label{fe6}
	C(\lambda)=\inf_{w\geqslant 1}\sup\big\{C_1w^2-(a+C_1+C_2), \lv b\rv-C_3w\big\}.
\end{equation}
For any ${\delta}>0$, by \eqref{fe6}, $C(\lambda)\geqslant {\delta}$ if and only if $\{w\geqslant1 \}$ is contained in
\begin{equation}
\{w\in\mathbb{R}\mid C_1w^2-(a+C_1+C_2)\geqslant{\delta}\}\cup\{v\in\mathbb{R}\mid\lv b\rv-C_3 w\geqslant{\delta}\},
\end{equation}
which means one of the following conditions holds
\begin{equation}\label{fe8}
	\begin{cases}
	C_1-(a+C_1+C_2)\geqslant{\delta},\\
	\sqrt{C_1^{-1}({a+C_1+C_2+{\delta}})}\leqslant C_3^{-1}({\lv b\rv-{\delta}}).
	\end{cases}
\end{equation}
By \eqref{fe8}, $C(\lambda)\geqslant {\delta}$ if and only if one of the following inequalities holds 
\begin{equation}\label{fe.9}
a\leqslant\begin{cases}
-C_2-{\delta},\\
{C_1}C_3^{-2}\lk\lv b\rv-{\delta}\rk^2-C_1-C_2-{\delta},
\end{cases}
\end{equation}
which is the red part of the figure below. By Theorem \ref{4.11}, let ${\delta}$ small enough that 
\begin{equation}
{C_1C_3^{-2}{\delta}^2}-C_1-C_2-{\delta}<\gamma/2< -C_2-{\delta},
\end{equation}
then by the figure below, we can choose a $d$ that satisfies $0<d<{C_1}{C_3^{-2}}$ to make sure that $V_{d,\gamma/2}$ is a \emph{subset} of the red part, in other words, for any $\lambda=a+bi\in V_{d,\gamma/2}$, it satisfies one of the two inequalities \eqref{fe.9}. Therefore, we have $\inf_{\lambda\in V_{d,\gamma/2}}C(\lambda)\geqslant {\delta}>0$.
\begin{figure}[ht]
	\centering
	\begin{tikzpicture}[decoration={markings, mark=at position 0.5 with {\arrow{>}}}] 
		\draw[->] (-5,0)--(9,0);
		\draw[->] (1.5,-3)--(1.5,3);
		\draw (0.5,-3)--(0.5,3);
		\draw[red] (-4.3,0)--(0.5,0);
		\draw[red] (-4.3,-2.1)--(8,-2.1);
		\draw[red] (-4.3,-1.7)--(4.8,-1.7);
		\draw[red] (-4.3,-1.3)--(1.2,-1.3);
		\draw[red] (-4.3,-0.9)--(0.5,-0.9);
		\draw[red] (-4.3,-0.5)--(0.5,-0.5);
		\draw[red] (-4.3,2.1)--(8,2.1);
		\draw[red] (-4.3,1.7)--(4.8,1.7);
		\draw[red] (-4.3,1.3)--(1.2,1.3);
		\draw[red] (-4.3,0.9)--(0.5,0.9);
		\draw[red] (-4.3,0.5)--(0.5,0.5);
		\node[fill=white]at (-0.8,2.9){\begin{footnotesize}$a=-C_2-{\delta}$\end{footnotesize}};
		\node[fill=white]at (4.2,2.9){\begin{footnotesize}$a=db^2-\gamma/2$\end{footnotesize}};
		\node[fill=white]at (5.5,1.2){\begin{footnotesize}$a=C_1C_3^{-2}\lk\lv b\rv-{\delta}\rk^2-C_1-C_2-{\delta}$\end{footnotesize}};
		\node at (9.2,0){\begin{footnotesize}$a$\end{footnotesize}};
		\node at (-1.9,0.2){\begin{footnotesize}$-\gamma/2$\end{footnotesize}};
		\node at (1.5,3.2){\begin{footnotesize}$b$\end{footnotesize}};
		\draw[postaction={decorate}] (3,2.2)--(2.9,2.2);
		\draw[postaction={decorate}] (2.9,-2.2)--(3,-2.2);
		\begin{scope}[rotate=-90]
			\draw[domain=-2:2] plot (\x,3*\x*\x-4);
		\end{scope}
		\begin{scope}[rotate=-90]
			\draw[domain=-2.5:2.5] plot (\x,0.9*\x*\x-1.4);
		\end{scope}
	\end{tikzpicture}
	\caption{}
	\label{fig3}
\end{figure}
Then by \eqref{4.4.6}, if $\lambda\in V_{d,\gamma/2}$, the 
	the resolvent  $\blk\lambda-\widehat{\mathscr{L}}^{F_p}_v\brk^{-1}$ exists and 
	\begin{equation}\label{4.4.8}
		\bV\blk\lambda-\widehat{\mathscr{L}}^{F_p}_v\brk^{-1}\bV^{0,0}_{p,\mu}\leqslant C.
	\end{equation}
	By \eqref{4.3.1.4}, it is clear that $\widehat{\mathscr{L}}^{F_p}_v$ can be extended to a continuous linear map from $H_{1,p,\mu}$ into $H_{-1,p,\mu}$. Using \eqref{4.4.4}, if we fix $\lambda_0=-1-C_1-C_2$, then
	\begin{equation}\label{4.4.9}
		\bV\blk\lambda_0-\widehat{\mathscr{L}}^{F_p}_v\brk^{-1}\bV^{-1,1}_{p,\mu}\leqslant 1.
	\end{equation}
	If $\lambda\in V_{d,\gamma/2}$, then
	\begin{equation}\label{4.4.10}
		\blk\lambda-\widehat{\mathscr{L}}^{F_p}_v\brk^{-1}=\blk\lambda_0-\widehat{\mathscr{L}}^{F_p}_v\brk^{-1}+\lk\lambda_0-\lambda\rk\blk\lambda-\widehat{\mathscr{L}}^{F_p}_v\brk^{-1}\blk\lambda_0-\widehat{\mathscr{L}}^{F_p}_v\brk^{-1}.
	\end{equation}
	By \eqref{4.4.8}, \eqref{4.4.9} and \eqref{4.4.10}, we see that
	\begin{equation}\label{4.4.11}
		\begin{split}
			\bV\blk\lambda-\widehat{\mathscr{L}}^{F_p}_v\brk^{-1}\bV_{p,\mu}^{-1,0}&\leqslant 1+C\lv\lambda-\lambda_0\rv,
		\end{split}
	\end{equation}
	on the other hand, we also have
	\begin{equation}\label{4.4.12}
		\blk\lambda-\widehat{\mathscr{L}}^{F_p}_v\brk^{-1}=\blk\lambda_0-\widehat{\mathscr{L}}^{F_p}_v\brk^{-1}+\lk\lambda_0-\lambda\rk\blk\lambda_0-\widehat{\mathscr{L}}^{F_p}_v\brk^{-1}\blk\lambda-\widehat{\mathscr{L}}^{F_p}_v\brk^{-1}.
	\end{equation}
	From \eqref{4.4.11} and \eqref{4.4.12}, we obtain 
	\begin{equation}
		\begin{split}
			\bV\blk\lambda-\widehat{\mathscr{L}}^{F_p}_v\brk^{-1}\bV_{p,\mu}^{-1,1}\leqslant 1+C\lv\lambda_0-\lambda\rv\lk1+C\lv\lambda-\lambda_0\rv\rk\leqslant C\blk1+\lv\lambda\rv^2\brk.
		\end{split}
	\end{equation}
	This completes the proof.\qed
\end{pro}

In \Cref{Fe}-\Cref{Fg}, we will always work under the assumption that $p\geqslant p_0$, $0\leqslant v\leqslant p^{-1}$ and $\lambda\in V_{d,{\gamma}/{2}}$ as in Theorem \ref{4.13}.
\subsection{Regularizing properties of the resolvents}\label{Fe}

In this subsection, we follow closely Ma-Marinescu \cite[Theorems 4.1.12-4.1.14]{mm07} with some necessary modifications.

\begin{theo}[Higher regularity]\label{4.14}
	For any $k\in\mathbb{N}$, the resolvent $\blk\lambda-\widehat{\mathscr{L}}^{F_p}_v\brk^{-1}$ in \eqref{fe2} maps $H_{k,p,\mu}$ to $H_{k+1,p,\mu}$. For any multi-index $\alpha\in\mathbb{N}^m$, there is $C>0$ such that
	\begin{equation}\label{fe18}
		\bV Z^\alpha\blk\lambda-\widehat{\mathscr{L}}^{F_p}_v\brk^{-1}s\bV_{k+1,p,\mu}\leqslant C\blk1+\lv\lambda\rv^2\brk^{k+\lv\alpha\rv+1}\sum_{\alpha'\leqslant\alpha}\bV Z^{\alpha'}s\bV_{k,p,\mu}
	\end{equation}
\end{theo}
\begin{pro}
	Using Proposition \ref{4.12} and Theorem \ref{4.13} in the same way as \cite[Theorem 1.6.10]{mm07} follows from \cite[Theorem 1.6.8, Proposition 1.6.9]{mm07}, we get Theorem \ref{4.14}.\qed
\end{pro}
For $k\in\mathbb{N}^*,r\in\mathbb{N}$, set
\begin{equation}
	I_{k,r}=\Big\{\lk\bm{k},\bm{r}\rk=\lk k_i,r_i\rk\in \lk\mathbb{N^*}\rk^{j+1}\times\mathbb{N}^{j}\ \Big| \sum_{i=0}^jk_i=k+j,\sum_{i=1}^jr_i=r\Big\}.
\end{equation}
For $(\bm{k},\bm{r})\in I_{k,r}$, we set
\begin{equation}\label{fe24}
	A^{\bm{k}}_{\bm{r}}\lk\lambda,v,p\rk=\blk\lambda-\widehat{\mathscr{L}}^{F_p}_v\brk^{-k_0}\frac{\pa^{r_1}\widehat{\mathscr{L}}^{F_p}_v}{\pa v^{r_1}}\blk\lambda-\widehat{\mathscr{L}}^{F_p}_v\brk^{-k_1}\cdots \frac{\pa^{r_j}\widehat{\mathscr{L}}^{F_p}_v}{\pa v^{r_j}}\blk\lambda-\widehat{\mathscr{L}}^{F_p}_v\brk^{-k_j}.
\end{equation}
Then there exist $a^{\bm{k}}_{\bm{r}}\in\mathbb{R}$ such that
\begin{equation}\label{4.4.21}
	\frac{\pa^{r}}{\pa v^{r}}\blk\lambda-\widehat{\mathscr{L}}^{F_p}_v\brk^{-k}=\sum_{({\bm{k}},{\bm{r}})\in I_{k,r}}a^{\bm{k}}_{\bm{r}}A^{\bm{k}}_{\bm{r}}\lk\lambda,v,p\rk.
\end{equation}

\begin{theo}\label{4.15}
	For any $\ell\in\mathbb{N}$, $k>2(\ell+ r+1)$ and $({\bm{k}},{\bm{r}})\in I_{k,r}$, there are $C> 0$
	and $j\in\mathbb{N}$ such that for any $\lambda\in  V_{d,\gamma/2}$ and $\alpha,\alpha'\in\mathbb{N}^m$ that $\lv\alpha\rv,\lv\alpha'\rv\leqslant \ell$, we have
	\begin{equation}\label{fe26}
		\bV\pa^{\alpha} A^{\bm{k}}_{\bm{r}}\lk\lambda,v,p\rk\pa^{\alpha'}s\bV_{0,p,\mu}\leqslant C(1+{\vert}\lambda{\vert})^j \sum_{{\vert}\alpha{\vert}\leqslant r}\lV Z^\alpha s\rV_{0,p,\mu}.
	\end{equation}
\end{theo}

\begin{pro}
	By Theorem \ref{4.14}, if $\lv\alpha\rv\leqslant \ell$, there is $C>0$ such that for any $\lambda\in V_{d,\gamma/2}$
	\begin{equation}\label{4.4.23}
		\bV\pa^\alpha\blk\lambda-\widehat{\mathscr{L}}^{F_p}_v\brk^{-\ell}\bV_{p,\mu}^{0,0}\leqslant C\blk1+\lv\lambda\rv^2\brk^{\ell+1}.
	\end{equation}

Note that if $\mu=1$ in \eqref{fd1}, the product $\langle\cdot,\cdot\rangle_{0,p,1}$ is just the usual $L^2$ inner product. Let $\widehat{\mathscr{L}}_{v}^{F_p,*}$ be the formal adjoint operator of $\widehat{\mathscr{L}}^{F_p}_v$ with respect to $\langle\cdot,\cdot\rangle_{0,p,1}$. Then $\widehat{\mathscr{L}}_{v}^{F_p,*}$ has essentially the same structure as the operator $\widehat{\mathscr{L}}^{F_p}$, except that the operators $e^i\wedge,\widehat{e}^i\wedge,{i}_{e_i},{i}_{\widehat{e}_i}$ are changed into ${i}_{e_i},{i}_{\widehat{e}_i},e^i\wedge,\widehat{e}^i\wedge$ respectively. We have
\begin{equation}\label{fe28}
\begin{split}
\big\langle s, \pa^\alpha\blk\lambda-\widehat{\mathscr{L}}^{F_p,*}_v\brk^{-\ell}s'\big\rangle_{0,p,1}=(-1)^{\vert\alpha\vert}\big\langle \blk\lambda-\widehat{\mathscr{L}}^{F_p}_v\brk^{-\ell}\pa^\alpha s,s'\big\rangle_{0,p,1}.
\end{split}
\end{equation}
 Note that
\begin{equation}\label{fe.30}
\langle s,s'\rangle_{0,p,1}=\langle \widehat{\theta}_\mu s,\widehat{\theta}_{\mu^{-1}}s'\rangle_{0,p,1}, \ \ \ \lV s\rV_{0,p,\mu}=\bV \widehat{\theta}_{\mu}s\bV_{0,p,1},\ \ \ \lV s'\rV_{0,p,\mu^{-1}}=\bV \widehat{\theta}_{\mu^{-1}}s'\bV_{0,p,1}.
\end{equation}
Now if we consider the weighted product $\lV\cdot\rV_{0,p,\mu^{-1}}$ for $\widehat{\mathscr{L}}_{v}^{F_p,*}$, an obvious analogue of \eqref{4.4.23} still holds: $\Vert(\lambda-\widehat{\mathscr{L}}^{F_p,*}_v)^{-\ell}\Vert^{0,0}_{p,\mu^{-1}}\leqslant C(1+\lv\lambda\rv^2)^{\ell+1}$. This together with \eqref{fe28} and \eqref{fe.30} yields
\begin{equation}\label{fe29}
	\begin{split}
\big\vert\big\langle s, \pa^\alpha\blk\lambda-\widehat{\mathscr{L}}^{F_p,*}_v\brk^{-\ell}s'\big\rangle_{0,p,1}\big\vert&\leqslant\Vert s\Vert_{0,p,\mu} \bV \blk\lambda-\widehat{\mathscr{L}}^{F_p,*}_v\brk^{-\ell} s'\bV_{0,p,\mu^{-1}} \\
&\leqslant C\blk1+\lv\lambda\rv^2\brk^{\ell+1} \Vert s\Vert_{0,p,\mu}\Vert s'\Vert_{0,p,\mu^{-1}}.
	\end{split}
\end{equation}
	By \eqref{fe28}, \eqref{fe.30} and \eqref{fe29}, we see that
	\begin{equation}\label{fe30}
\begin{split}
\bV\blk\lambda-\widehat{\mathscr{L}}^{F_p}_v\brk^{-\ell}\pa^\alpha s\bV_{0,p,\mu}&=\sup_{\Vert s'\Vert_{0,p,\mu^{-1}}\leqslant 1}\Bv\big\langle \blk\lambda-\widehat{\mathscr{L}}^{F_p}_v\brk^{-\ell}\pa^\alpha s,s'\big\rangle_{0,p,1}\Bv\\
&\leqslant C\blk1+\lv\lambda\rv^2\brk^{\ell+1}\bV s\bV_{0,p,\mu}.
\end{split}
	\end{equation} 
	By \eqref{4.4.23} and \eqref{fe30}, we prove \eqref{fe26} for $r=0$.
	
For $r>0$, by \eqref{fc37} and \eqref{fd9}, we can rewrite $\frac{\pa^{r}}{\pa v^{r}}\widehat{\mathscr{L}}_v^{F_p}$ in \eqref{fe24} as
\begin{equation}\label{fe31}
\sum_{\vert\beta\vert\leqslant r}C_\beta \Gamma_{\beta,{vZ}}\pa_i\pa_jZ^{\beta}+C''_\beta \Gamma'_{\beta,{vZ}}\pa_iZ^{\beta}+C''_\beta\Gamma''_{\beta,{vZ}}Z^{\beta},
\end{equation}
where $\Gamma_{\beta,Z},\Gamma_{\beta',Z},\Gamma_{\beta'',Z}\in\mathscr{C}_0^\infty(T_{x_0}X,\widehat{\mathbb{E}}_{p,x_0})$ are smooth families of Toeplitz operators (see \cref{Dc}). Let $\mathscr{R}'_{v,p}$ be the family of operators of the form
\begin{equation}
	\mathscr{R}_{v,p}'=\left\{\big[\Gamma_{1}Q_{1},\big[\Gamma_{2}Q_{2},\cdots\big[\Gamma_{k}Q_{k},\widehat{\mathscr{L}}^{F_p}_v\big]\cdots\big]\big]\right\}
\end{equation}
where $\Gamma_{i}\in\mathscr{C}_0^\infty(T_{x_0}X,\widehat{\mathbb{E}}_{p,x_0})$ are bounded with respect to $\lv\cdot\rv_{\mathscr{C}^0(T_{x_0}X,\widehat{\mathbb{E}}_{p,x_0})}$ uniformly for $p\in\mathbb{N}^*$ with all derivatives and $Q_1,\cdots Q_k\in\{\pa_i,Z_i\}_{i=1}^{m}$.

We handle now the operator $A^{\bm{k}}_{\bm{r}}\lk\lambda,v,p\rk\pa^{\alpha'}$. For each $\frac{\pa^{r_i}\widehat{\mathscr{L}}^{F_p}_v}{\pa v^{r_i}}$, we  write it in the form of \eqref{fe31} and move all the terms to the right-hand side as in \cite[Theorem 4.1.13]{mm07}. Then $\pa^\alpha A^{\bm{k}}_{\bm{r}}\lk\lambda,v,p\rk\pa^{\alpha'}$ is as the form $\sum_{\vert\beta\vert\leqslant r}\mathscr{L}_{v,\beta}Z^\beta$ where $\mathscr{L}_{v,\beta}$ is a linear combination of operators that can be split into the product of two parts:	\begin{equation}
		\begin{split}
			&\pa^{\alpha}\blk\lambda-\widehat{\mathscr{L}}^{F_p}_v\brk^{-k_0'}R_1\blk\lambda-\widehat{\mathscr{L}}^{F_p}_v\brk^{-k_1'}\cdots R_i\blk\lambda-\widehat{\mathscr{L}}^{F_p}_v\brk^{-k_i''}\\
			&\blk\lambda-\widehat{\mathscr{L}}^{F_p}_v\brk^{-(k_i'-k_i'')}\cdot R_{i+1}\cdots R_{j'}\blk\lambda-\widehat{\mathscr{L}}^{F_p}_v\brk^{-k_{j'}'}\pa^{\beta'},
		\end{split}
	\end{equation}
where $R_i\in \mathscr{R}_{v,p}'$, $\sum_{j=0}^{i-1}(k'_j-1)+k''_i-\ell>0$ and $\sum_{j=i}^{j'}(k'_j-1)-k''_i>2r+\ell+1$, then the $\lV\cdot\rV^{0,0}_{p,\mu}$ norm of each part is bounded by $C(1+{\vert}\lambda{\vert})^N$. This finishes the proof.\qed
	\end{pro}

\begin{theo}\label{4.16}
	For any $r\geqslant 0$ and $k > 0$, there exist $C > 0$ and $\ell\in\mathbb{N}^*$ such that
	\begin{equation}\label{4.4.31}
		\begin{split}
&\bigg\Vert\Big(\frac{\pa^{r}\widehat{\mathscr{L}}^{F_p}_v}{\pa v^{r}}-\frac{\pa^{r}\widehat{\mathscr{L}}^{F_p}_v}{\pa v^{r}}\Big|_{v=0}\Big) s\bigg\Vert_{-1,p,\mu}\leqslant Cv\sum_{|\alpha|\leqslant r+1}\lV Z^\alpha s\rV_{1,p,\mu},\\
&\bigg\Vert \Big(\frac{\pa^{r}}{\pa v^{r}}\blk\lambda-\widehat{\mathscr{L}}^{F_p}_v\brk^{-k}-\sum_{\lk{\bm{k}},{\bm{r}}\rk\in I_{k,r}}a^{\bm{k}}_{\bm{r}}A^{\bm{k}}_{\bm{r}}\lk\lambda,0,p\rk \Big)s\bigg\Vert_{0,p,\mu}\leqslant Cv\lk1+\vert\lambda\vert\rk^\ell\sum_{|\alpha|\leqslant 4r+1}\lV Z^\alpha s\rV_{0,p,\mu}.
		\end{split}
	\end{equation}
\end{theo}

\begin{pro}
	An application of Taylor expansion for \eqref{fc35} and \eqref{fc37} implies 
	\begin{equation}
		\Bli\Blk\frac{\pa^{r}\widehat{\mathscr{L}}^{F_p}_v}{\pa v^{r}}-\frac{\pa^{r}\widehat{\mathscr{L}}^{F_p}_v}{\pa v^{r}}\Big|_{v=0}\Brk s,s'\Bri_{0,p,\mu}\leqslant {Cv}\lV s'\rV_{1,p,\mu}\sum_{|\alpha|\leqslant r+1}\lV Z^\alpha s\rV_{1,p,\mu},
	\end{equation}
from which we get the first inequality of \eqref{4.4.31}.  The second inequality of \eqref{4.4.31} follows from Theorems \ref{4.13}, \ref{4.14} and the first inequality of \eqref{4.4.31} as \cite[Theorem 4.1.14]{mm07} follows from  \cite[Theorems 4.1.10, 4.1.12]{mm07}.\qed
\end{pro}

\subsection{Uniform estimation on the heat kernel}\label{Ff}

This section is analogous to \cite[\S\,4.1.5]{mm07}, with the necessary changes made.

\begin{theo}\label{4.17}
	For any $k,\ell,r\in\mathbb{N}$, there is $C>0$ such that for $t\geqslant1, Z,Z'\in T_{x_0}X$ with $\lv Z\rv,\lv Z'\rv\leqslant1$, we have
	\begin{equation}\label{4.4.38}
		\sup_{\lv\alpha\rv,\lv\alpha'\rv\leqslant k}\bigg\vert\frac{\pa^{{\vert}\alpha{\vert}+{\vert}\alpha'{\vert}}}{\pa Z^\alpha\pa {Z'}^{\alpha'}}\frac{\pa^{r}}{\pa v^{r}}\exp(-t\widehat{\mathscr{L}}^{F_p}_v)\lk Z,Z'\rk\bigg\vert_{\mathscr{C}^\ell(M)}\leqslant Ce^{\gamma t/2}.
	\end{equation}
\end{theo}

\begin{pro}
Recall $V_{d,\gamma/2}$ given in \eqref{fe1}, by the Cauchy integral formula, for $k\in \mathbb{N}^*$, we have
	\begin{equation}\label{4.4.39}
		\exp(-t\widehat{\mathscr{L}}^{F_p}_v)=\frac{(k-1)!}{2\pi {i}\lk -t\rk^{k-1}}\int_{\pa V_{d,\gamma/2}}e^{-t\lambda}\blk\lambda-\widehat{\mathscr{L}}^{F_p}_v\brk^{-k}{d} \lambda.
	\end{equation}
From Theorem \ref{4.15} and \eqref{4.4.39}, we obtain that for $\lv\alpha_i\rv\leqslant \ell$ and $t\geqslant 1$, we have
	\begin{equation}\label{4.4.40}
		\big\Vert\pa^{\alpha_1}\exp(-t\widehat{\mathscr{L}}^{F_p}_v)\pa^{\alpha'}\big\Vert_{p,\mu}^{0,0}\leqslant Ce^{\gamma t/2}.
	\end{equation}
	Now by \eqref{4.4.40} and the Sobolev inequality, for $Z,Z'\in T_{x_0}X$ that $\lv Z\rv,\lv Z'\rv\leqslant1$, we have
	\begin{equation}
		\sup_{\lv\alpha\rv,\lv\alpha'\rv\leqslant l}\bigg\vert\frac{\pa^{{\vert}\alpha{\vert}+{\vert}\alpha'{\vert}}}{\pa Z^\alpha\pa {Z'}^{\alpha'}}\exp(-t\widehat{\mathscr{L}}^{F_p}_v)\lk Z,Z'\rk\bigg\vert\leqslant Ce^{\gamma t/2}.
	\end{equation}
This implies \eqref{4.4.38} for $r=\ell=0$. To obtain \eqref{4.4.38} for $r\geqslant1$, we see from \eqref{4.4.39} that for $k\in\mathbb{N}^*$,
	\begin{equation}\label{4.4.42}
		\frac{\pa^{r}}{\pa v^{r}}\exp(-t\widehat{\mathscr{L}}^{F_p}_v)=\frac{(k-1)!}{2\pi {i}\lk-t\rk^{k-1}}\int_{\pa V_{d,\gamma/2}}e^{-t\lambda}\frac{\pa^{r}}{\pa v^{r}}\blk\lambda-\widehat{\mathscr{L}}^{F_p}_v\brk^{-k}d \lambda.
	\end{equation}
	By \eqref{4.4.21}, \eqref{4.4.23} and \eqref{4.4.42}, we get
	\begin{equation}
		\BV\frac{\pa^{r}}{\pa v^{r}}\pa^{\alpha_1}\exp(-t\widehat{\mathscr{L}}^{F_p}_v)\pa^{\alpha'}s\BV_{0,p,\mu}\leqslant Ce^{\gamma t/2}\sum_{\vert\alpha\vert\leqslant r}\Vert Z^\alpha s\Vert_{0,p,\mu}.
	\end{equation} 
	Taking $\mathrm{supp}(s)\in B^{T_{x_0}X}(0,2)$, by Sobolev inequality again, we get \eqref{4.4.38} for $\ell=0$.
	
	Finally, for any vector $U$ on $M$,
	\begin{equation}
\nabla^{\pi^*\widehat{\mathbb{E}}_p}_U\exp(-t\widehat{\mathscr{L}}^{F_p}_v)=\frac{(k-1)!}{2\pi {i}\lk-t\rk^{k-1}}\int_{\pa V_{d,\gamma/2}}e^{-t\lambda}\nabla^{\pi^*\widehat{\mathbb{E}}_p}_U\blk\lambda-\widehat{\mathscr{L}}^{F_p}_v\brk^{-k}d \lambda.
	\end{equation}
	Now we use a similar formula as \eqref{4.4.21} for $\nabla^{\pi^*\widehat{\mathbb{E}}_p}_U\blk\lambda-\widehat{\mathscr{L}}^{F_p}_v\brk^{-k}$ by replacing $\frac{\pa^{r_i}}{\pa v^{r_i}}\widehat{\mathscr{L}}^{F_p}_v$ with $\nabla^{\pi^*\widehat{\mathbb{E}}_p}_U\widehat{\mathscr{L}}^{F_p}_v$. Remark that $\nabla^{\pi^*\widehat{\mathbb{E}}_p}_U\widehat{\mathscr{L}}^{F_p}_v$ is a differential operator on $T_{x_0}X$ with the same structure as $\widehat{\mathscr{L}}^{F_p}_v$, it has the same type as \eqref{fc35}. Then using the above argument, we conclude that \eqref{4.4.21} also holds for $\ell\geqslant1$. We complete the proof.\qed
\end{pro}

For $k\in\mathbb{N}^*$ large enough, set
\begin{equation}\label{4.4.44}
	\begin{split}
		\mathscr{S}_{t,p}^{(r)}&=\frac{(k-1)!}{2\pi {i}(-t)^{k-1}r!}\int_{\pa V_{d,\gamma/2}}e^{-t\lambda}\sum_{({\bm{k}},{\bm{r}})\in I_{k,r}}a^{\bm{k}}_{\bm{r}}A^{\bm{k}}_{\bm{r}}\lk\lambda,0,p\rk{d} \lambda,\\
		\mathscr{S}^{(r)}_{v,t,p}&=\frac{1}{r!}\frac{\pa^{r}}{\pa v^{r}}\exp(-t\widehat{\mathscr{L}}^{F_p}_v)-\mathscr{S}_{t,p}^{(r)}(t),
	\end{split}
\end{equation}
then $\mathscr{S}^{(r)}_{t,p}$ and $\mathscr{S}^{(r)}_{v,t,p}$ do \emph{not} depend on the choice of $k$.  We denote by $\mathscr{S}_{t,p}^{(r)}(Z,Z')$ (resp. $\mathscr{S}_{v,t,p}^{(r)}(Z,Z')$) the smooth kernel of $\mathscr{S}^{(r)}_{p}(t)$ (resp. $\mathscr{S}^{(r)}_{v,t,p}$) with respect to ${d} v_{TX}$.

\begin{theo}\label{4.18}	
	For any $r\in\mathbb{N}$, there is $C>0$ such that for $t\geqslant1,Z,Z'\in T_{x_0}X$ with $\lv Z\rv,\lv Z'\rv\leqslant1$, we have
	\begin{equation}\label{4.4.45}
		\lv\mathscr{S}^{(r)}_{v,t,p}(Z,Z')\rv\leqslant Cv^{1/(m+1)}e^{\gamma t/2}.
	\end{equation}
\end{theo}

\begin{pro}
	By \eqref{4.4.44} and the Cauchy integral formula, we get
	\begin{equation}\label{ff10}
		\begin{split}
			\mathscr{S}^{(r)}_{v,t,p}=\frac{(k-1)!}{2\pi {i}\lk -t\rk^{k-1}r!}\int_{\pa V_{d,\gamma/2}}e^{-t\lambda}\cdot\Big(\frac{\pa^{r}}{\pa v^{r}}\blk\lambda&-\widehat{\mathscr{L}}^{F_p}_v\brk^{-k}\\
			&-\sum_{({\bm{k}},{\bm{r}})\in I_{k,r}}a^{\bm{k}}_{\bm{r}}A^{\bm{k}}_{\bm{r}}\lk\lambda,0,p\rk\Big) d\lambda.
		\end{split}
	\end{equation}
By \eqref{4.4.31} and \eqref{ff10}, when $t\geqslant 1$, we have
	\begin{equation}\label{4.4.47}
		\big\Vert \mathscr{S}^{(r)}_{v,t,p}(t)s\big\Vert_{0,p,\mu}\leqslant Cve^{\gamma t/2}\sum_{|\alpha|\leqslant 4r+1}\lV Z^\alpha s\rV_{0,p,\mu}.
	\end{equation}
By Theorem \ref{4.17} and \eqref{4.4.47}, using a similar argument as \cite[Theorem 4.17]{mm07} follows from \cite[Theorem 4.16, (4.1.67)]{mm07}, we get \eqref{4.4.45}.\qed
\end{pro}

\begin{theo}\label{4.19}
	For any $r,\ell\in\mathbb{N}$, there is $C > 0$ such that for $t\geqslant 1$, we have 
	\begin{equation}\label{4.4.50}
		\Big\vert\exp(-t\widehat{\mathscr{L}}^{F_p}_v)(0,0)-\sum_{i=0}^r\mathscr{S}^{(i)}_{t,p}(0,0)v^i\Big\vert_{\mathscr{C}^\ell(M)}\leqslant Cv^{r+1}e^{\gamma t/2}.
	\end{equation}
\end{theo}

\begin{pro}
	By Theorem \ref{4.18} and \eqref{4.4.44}, we have
	\begin{equation}\label{4.135}
		\frac{1}{r!}\frac{\pa^{r}}{\pa v^{r}}\exp(-t\widehat{\mathscr{L}}^{F_p}_v)\Big{\vert}_{v=0}=\mathscr{S}^{(r)}_{t,p}.
	\end{equation}
	Now by Theorem \ref{4.17} and \eqref{4.4.44}, $\mathscr{S}^{(r)}_{t,p}$ has the same estimation as $\frac{\pa^{r}}{\pa v^{r}}\exp(-t\widehat{\mathscr{L}}^{F_p}_v)$ in \eqref{4.4.38}. From \eqref{4.4.38} and the Taylor expansion
	\begin{equation}
		f(v)-\sum_{r=0}^k\frac{1}{r!}\frac{\pa^{r}f}{\pa v^{r}}(0)v^r=\frac{1}{k!}\int_{0}^v(v-y)^k\frac{\pa^{k+1}f}{\pa v^{k+1}}(y)\de y,
	\end{equation}
	we get \eqref{4.4.50}.\qed	
\end{pro}

\subsection{The asymptotics of $\mathscr{S}^{(i)}_{t,p}(0,0)$}\label{Fg}

Let $\mathscr{P}^{F_p}_t(Z,Z')$ be the kernel of $\exp(-t\widehat{\mathscr{L}}^{F_p}_{0})$ with respect to $dv_{TX}$. Put
\begin{equation}\label{fg6}
	\begin{split}
		\sigma^{F_p}=&-\frac{1}{4}\big\langle R^{TX}e_i,e_j \big\rangle_{x_0}\widehat{e}^i\wedge\widehat{e}^j-\frac{1}{4p}\omega^{F_p,2}_{x_0}+\frac{1}{4p^2}\big\vert\widehat{\omega}^{F_p}\big\vert^2_{x_0}+\frac{1}{4p}\widehat{\omega}^{F_p,2}_{x_0}\\
		&-\frac{1}{2p}\big(\nabla^{\widehat{TX}\otimes F_p,u}\widehat{\omega}^{F_p}\big)_{x_0}-\frac{1}{2p}z\omega^{F_p}_{x_0}.
	\end{split}
\end{equation}
By \eqref{fc28}, we have
\begin{equation}\label{fg7}
	\mathscr{P}^{F_p}_t(Z,Z')=\frac{1}{(4\pi t)^{\frac{m}{2}}}e^{-\frac{\lv Z-Z'\rv^2}{4t}-t\sigma^{F_p}}.
\end{equation}

For $k\in\mathbb{N}$, the $k$-simplex $\triangle_k$ is given by $\{(t_1,\cdots,t_k)\mid 0\leqslant t_1\leqslant t_2 \cdots\leqslant t_k\leqslant 1 \}$.
For $t>0$, we will write $t\triangle_k$ for the rescaled simplex $\{(t_1,\cdots,t_k)\mid 0\leqslant t_1\leqslant t_2 \cdots\leqslant t_k\leqslant t \}$.

For multi-index ${\bm{r}}=(r_1,\cdots,r_k)_{k\in\mathbb{N}}$ where $r_1,\cdots,r_k\in\mathbb{N}^*$, set
\begin{equation}\label{fg8}
	\begin{split}
	\mathscr{S}_{t,{\bm{r}},p}&=\int_{t\triangle_k}\bigg(\mathscr{P}^{F_p}_{t-t_k}\frac{\pa^{r_k}\widehat{\mathscr{L}}^{F_p}_v}{\pa v^{r_k}}\Big|_{v=0}\mathscr{P}^{F_p}_{t_k-t_{k-1}}\cdots\frac{\pa^{r_1}\widehat{\mathscr{L}}^{F_p}_v}{\pa v^{r_1}}\Big|_{v=0}\mathscr{P}^{F_p}_{t_1}\bigg)(0,0)\prod_{j=1}^kdt_j\\
	&=\int_{t\triangle_k}\int_{(TX)^k}\mathscr{P}^{F_p}_{t-t_k}(0,Z^{{(k)}})\bigg(\frac{\pa^{r_k}\widehat{\mathscr{L}}^{F_p}_{v,Z^{{(k)}}}}{\pa v^{r_k}}\Big|_{v=0}\mathscr{P}^{F_p}_{t_k-t_{k-1}}\bigg)(Z^{{(k)}},Z^{{(k-1)}})\\
	&\ \ \ \ \ \ \ \ \ \ \ \ \ \cdots \bigg(\frac{\pa^{r_1}\widehat{\mathscr{L}}^{F_p}_{v,Z^{{(1)}}}}{\pa v^{r_1}}\Big|_{v=0}\mathscr{P}^{F_p}_{t_1}\bigg)(Z^{{(1)}},0)\prod_{j=1}^kdv_{TX}\big(Z^{{(j)}}\big)\prod_{j=1}^kdt_j
	\end{split}
\end{equation}
where $(TX)^k$ means integral with respect to $k$-copies of $T_{x_0}X$: $\big\{Z^{{(i)}}=\big(Z^{(i)}_1,\cdots,Z^{(i)}_m\big)\mid 1\leqslant i\leqslant k\big\}$, and the subscript $Z^{(j)}$ in $\frac{\pa^{r_j}\widehat{\mathscr{L}}^{F_p}_{v,Z^{{(j)}}}}{\pa v^{r_j}}$ means acting on the coordinate $Z^{(j)}$. By the Duhamel's principle \cite[\S\,2.7]{bgv}, Theorem \ref{Fc13}, \eqref{fe24}, \eqref{4.4.44} and \eqref{4.135}, we have
\begin{equation}\label{fg9}
	\mathscr{S}^{(i)}_{t,p}(0,0)=\sum_{\substack{ r_1+\cdots +r_k=i}}\frac{(-1)^k}{\prod_{j=1}^k r_j!}\mathscr{S}_{t,{\bm{r}},p}.
\end{equation}
In particular,
\begin{equation}\label{fg..5}
	\begin{split}
\mathscr{S}^{(0)}_{t,p}(0,0)=\frac{e^{-t\sigma^{F_p}}}{(4\pi t)^{\frac{m}{2}}}, \  \mathscr{S}^{(1)}_{t,p}(0,0)=-\int_{0}^t\int_{T_{x_0}X}\mathscr{P}^{F_p}_{t-t_1}(0,Z^{})\bigg(\frac{\pa\widehat{\mathscr{L}}^{F_p}_{v,Z}}{\pa v}\Big|_{v=0}\mathscr{P}^{F_p}_{t_1}\bigg)(Z,0)dZdt_1.
	\end{split}
\end{equation}

By \eqref{fg7}, we see that there are singularities inside the integral of \eqref{fg8}, and we will prove that there are no singularities after taking the integral. Let us give an auxiliary result first. Put $P_t(x,y)$ the classical heat kernel on $\mathbb{R}$:
$P_t(x,y)={(4\pi t)}^{-1/2}e^{-(x-y)^2/4t}$. For $t>0$, $(t_1,\cdots,t_k)\in t\triangle_k$ and $\{\alpha_i^{(j)},\beta_i^{(j)}\in\mathbb{N}\}_{1\leqslant i\leqslant m,1\leqslant j\leqslant k}$, set
\begin{equation}\label{fg.5}
\begin{split}&f_t(t_k,\cdots,t_1)=
\prod_{i=1}^m \Bigg\{P_{t-t_k}(0,Z^{(k)}_{i})\cdot Z_{i}^{(k),\alpha^{(k)}_{i}}\frac{\pa^{\beta^{(k)}_{i}}}{\pa Z_{i}^{(k),\beta^{(k)}_{i}}}P_{t_k-t_{k-1}}(Z^{(k)}_{i},Z^{(k-1)}_{i})\\
 &\cdots Z_{i}^{(2),\alpha^{(2)}_{i}}\frac{\pa^{\beta^{(2)}_{i}}}{\pa Z_{i}^{(2),\beta^{(2)}_{i}}}P_{t_2-t_1}(Z^{(2)}_{i},Z^{(1)}_{i})\cdot Z_{i}^{(1),\alpha_{i}}\frac{\pa^{\beta_{i}}}{\pa Z_{i}^{(1),\beta_{i}}}P_{t_1}(Z^{(1)}_{i},0)\prod_{j=1}^kdZ^{(j)}_i\Bigg\}.
	\end{split}
\end{equation}
\begin{lemma}
For $k\in\mathbb{N}^*$, there are $\ell\in\mathbb{N}$ and $C>0$ such that for $t\geqslant1$ and $(t_1,\cdots,t_k)\in t\triangle_k$, we have
\begin{equation}\label{fg.6}
\lv f_t(t_k,\cdots,t_1)\rv\leqslant	C\big(1+t^{\ell}\big).
\end{equation}
\end{lemma}
\begin{pro}
By \eqref{fg.5}, we only need to prove for the special case $m=1$. We prove by induction on $k\in\mathbb{N}^*$. Let us consider the integral of the following form:
\begin{equation}\label{fg10}
	\int_{-\infty}^{\infty}P_{t_1}(Z_2,Z_1)Z_1^{{\ell}'}\frac{\pa^\ell}{\pa Z_1^\ell}P_{t_0}(Z_1,Z_0)dZ_1, \ \ \text{for\ }Z_0,Z_1,Z_2\in\mathbb{R}, t_0,t_1>0.
\end{equation}
We introduce the generating function of the integral in \eqref{fg10}. Put
\begin{equation}\label{fg11}
	f_1(Z_2,Z_0)=\sum_{\ell,\ell'\in\mathbb{N}}\int_{-\infty}^{\infty}P_{t_1}(Z_2,Z_1)\frac{(\mu_1 Z_1)^{{\ell}'}}{({\ell}')!}\frac{\lambda_1^\ell}{\ell!}\frac{\pa^\ell}{\pa Z_1^\ell}P_{t_0}(Z_1,Z_0)dZ_1,
\end{equation}
For any analytic function $g(y)$, we have $\sum_{\ell\in\mathbb{N}}\frac{\lambda^\ell}{\ell!}\frac{\pa^\ell}{\pa y^\ell}g(y)=g(y+\lambda)$, this gives
\begin{equation}\label{fg13}
	\begin{split}
		&f_1(Z_2,Z_0)=\int_{-\infty}^{\infty}P_{t_1}(Z_2,Z_1)\exp(\mu_1Z_1)P_{t_0}(Z_1+\lambda_1,Z_0)dZ_1\\
		&=\int_{-\infty}^{\infty}\frac{1}{4\pi (t_0t_1)^{\frac{1}{2}}}\exp\Big({-\frac{t_1+t_0}{4t_1t_0}\Big(Z_1-\frac{t_0Z_2}{t_1+t_0}-\frac{t_1(Z_0-\lambda_1)}{t_1+t_0}-\frac{2t_1t_0\mu_1}{t_1+t_0}\Big)^2}\Big)dZ_1\\
		&\ \ \ \cdot\exp\Big({-\frac{(Z_2-Z_0+\lambda_1)^2}{4(t_0+t_1)}+\mu_1\Big(\frac{t_0Z_2}{t_1+t_0}+\frac{t_1(Z_0-\lambda_1)}{t_1+t_0}\Big)+\frac{\mu_1^2t_0t_1}{(t_1+t_0)}}\Big)\\
		&=P_{t_1+t_0}(Z_2,Z_0-\lambda_1)\exp\Big({\mu_1\Big(\frac{t_0Z_2}{t_1+t_0}+\frac{t_1(Z_0-\lambda_1)}{t_1+t_0}\Big)+\frac{\mu_1^2t_0t_1}{(t_1+t_0)}}\Big).
	\end{split}
\end{equation}
Similar to \eqref{fg11} and \eqref{fg13}, we define
\begin{equation}\label{fg14}
	\begin{split}
		&f_j(Z_{j+1},Z_0)=\int_{-\infty}^{\infty}P_{t_j}(Z_{j+1},Z_j)\exp(\mu_jZ_j)f_{j-1}(Z_j+\lambda_j,Z_0)dZ_j,\\
		&w_0=0,\ \ w_j=\Big(\sum_{i=0}^jt_i\Big)^{-1}{\sum_{0\leqslant i<\ell\leqslant j}\mu_\ell t_i}\ \ \text{for\ }j\in\mathbb{N}^*.
	\end{split}
\end{equation}
Then we can prove inductively that
\begin{equation}\label{fg17..}
	\begin{split}
		f_j(Z_{j+1},Z_0)=P_{\sum_{i=0}^jt_i}\Big(Z_{j+1},Z_0-\sum_{i=1}^j\lambda_j\Big)\exp\Big(w_jZ_{j+1}+\sum_{i=1}^{j}(\mu_j+w_{j-1})w_it_i&\\
		+\sum_{i=1}^{j}(\mu_{i}+w_{i-1}-w_{i})(Z_0-\sum_{\ell=1}^{i}\lambda_\ell)\Big)&.
	\end{split}
\end{equation}
Plugging $Z_0=Z_{j+1}=0$ in \eqref{fg17..}, we obtain
\begin{equation}\label{fg17}
	\begin{split}
		f_j(0,0)=\Big(4\pi\sum_{i=0}^jt_i\Big)^{-\frac{1}{2}}\exp\Big(-\Big(4\sum_{i=0}^j\lambda_i\Big)^{-1}\Big(\sum_{i=1}^{j}\lambda_i\Big)^2+\sum_{i=1}^{j}(\mu_j+w_{j-1})w_it_i\\
		-\sum_{i=1}^{j}(\mu_{i}+w_{i-1}-w_{i})\Big(\sum_{\ell=1}^{i}\lambda_\ell\Big)\Big).
	\end{split}
\end{equation}
By \eqref{fg15} and \eqref{fg17}, for $k\in\mathbb{N}$, there is $C>0$ such that for any multi-index $\alpha,\beta\in\mathbb{N}^{j}$ with $\lv\alpha\rv,\lv\beta\rv\leqslant k$ and $\sum_{i=0}^jt_j=t\geqslant1$, the coefficient of $\lambda^\alpha\mu^\beta$ in $f_j(0,0)$ is dominated by $C(1+t^{k})$, which implies \eqref{fg.6}.\qed
\end{pro}

Recall the asymptotic trace symbol $\tro_{[i]}[T_p]$ in \eqref{cb20}.
\begin{theo}\label{Fg26}
For $i,j,\ell\in \mathbb{N}$ and ${\delta}>0$, there exists $C>0$ such that for $t\geqslant1$,
\begin{equation}\label{fg24}
	\begin{split}
			&\Bv\tro_{[j]}\big[\mathscr{S}^{(i)}_{t,p}(0,0)\big]\Bv_{\mathscr{C}^\ell(M)}\leqslant Ce^{{\delta} t},\\
				&\bbv p^{-n}\tro^{F_p}[\mathscr{S}^{(i)}_{t,p}(0,0)]-\sum_{k=0}^jp^{-k}\tro_{[k]}[\mathscr{S}^{(i)}_{t,p}(0,0)]\bbv_{\mathscr{C}^\ell(M)}\leqslant Ce^{{\delta} t}p^{-j-1}.
	\end{split}
\end{equation}
\end{theo}
\begin{pro}
By Theoreom \ref{Fc13}, \eqref{fg7}, \eqref{fg8} \eqref{fg.6}, $\mathscr{S}_{{\bm{r}},p}$ is a sum of integrals of the form
\begin{equation}\label{fg20}
	\int_{t\triangle_k}\bm{f}_t(t_k,\cdots,t_1){e}^{-(t-t_k)\sigma^{F_p}}\tau_k{e}^{-(t_k-t_{k-1})\sigma^{F_p}}\cdots\tau_1{e}^{-(t_1-t_{0})\sigma^{F_p}},
\end{equation}
where $\tau_{j}$ is a smooth family of Toeplitz operators (see \cref{Dc}) with respect to the parameter $x_0\in M$ for $0\leqslant j\leqslant k$, and $\bm{f}_t(t_k,\cdots,t_1)$ is the sum of the product of functions in the form of \eqref{fg.5}. By \eqref{fg.6}, we get
\begin{equation}\label{fg21}
	\bv \bm{f}_t(t_k,\cdots,t_1)\bv\leqslant C(1+t^\ell)
\end{equation}
for some $C>0,\ell\in\mathbb{N}$. By \eqref{fg6}, it is clear that
\begin{equation}\label{fg22}
	\mathrm{Spec}\big(\sigma^{F_p}\big)=\mathrm{Spec}\Big(\frac{1}{4p^2}\big\vert\widehat{\omega}^{F_p}\big\vert^2_{x_0}\Big).
\end{equation}
By Theorem \ref{Cc10}, \eqref{fg9}, \eqref{fg20}, \eqref{fg21}, \eqref{fg22}, $\mathscr{S}^{(i)}_{t,p}(0,0)$ verifies estimations the same as \eqref{cc41} and \eqref{cc42}, from which we get \eqref{fg24} for $\ell=0$. Since $\mathscr{S}^{(i)}_{t,p}(0,0)$ is a smooth family of Toeplitz operators with respect to $x_0\in M$, by Remark \ref{dc6}, we get \eqref{fg24} for the general norm $\lv\cdot\rv_{\mathscr{C}^\ell(M)}$.\qed
\end{pro}

\subsection{The convergence when $0<t\leqslant 1$}\label{Fi}

As $t\rightarrow 0$, there is singularity in \eqref{4.4.50}, we consider a different rescaling following the proof of \cite[Theorem 9.31]{bmz17}. In this subsection, we always assume that $0<t\leqslant 1$. Set
\begin{equation}
	{\mathscr{N}}^{F_p}_{t,x_0}=\frac{t}{p^2}\theta_{\sqrt{\frac{p}{t}}}K_{\frac{\sqrt{t}}{p}}{\mathscr{L}}_{x_0}^{F_p}K_{\frac{p}{\sqrt{t}}}\theta_{\sqrt{\frac{t}{p}}}.
\end{equation}
We denote by $\widehat{\mathscr{L}}^{F_p}_{t,x_0}$ the operator obtained from ${\mathscr{N}}^{F_p}_{t,x_0}$ by replacing $c(e_i), \widehat{c}(e_i)$ with $c_{\frac{t}{p}}(e_i)$ and $\widehat{c}_{\frac{1}{p}}(\widehat{e}_i)$ respectively. Like \eqref{4.2.8}, we have
\begin{equation}\label{fi.2}
	\begin{split}
\theta_{\frac{1}{\sqrt{t}}}\tro_s^{\Lambda(T^*X)\otimes F_p}\big[\exp(-t{\mathscr{M}}_{x_0}^{F_p})\lk0,0\rk\big]dv_X=(4\pi)^{\frac{m}{2}}\tro^{F_p}\int^{\widehat{B}}\Big[\exp(-\widehat{\mathscr{L}}^{F_p}_{t,x_0})(0,0)\Big]^{\mathrm{max}}.
	\end{split}
\end{equation}

Let ${N_2}$ and $N_3$ be the number operators of the exterior algebras $\mathbb{R}[z]\hti\pi^*\Lambda\lk T^*S\rk\hti\Lambda\lk T^*X\rk$ and $\Lambda(\widehat{T^*X})$ respectively  acting by multiplication of the degree. For any $a>0$, set
\begin{equation}\label{fi.3..}
	\overline{\theta}_a=a^{N_2}, \ \ \ \ \ \ \widetilde{\theta}_a=a^{N_3}.
\end{equation} 
By evaluating the tensors on the right hand side at ${v\sqrt{t}Z}$, we take
\begin{equation}\label{fi4}
	\begin{split}
		\widehat{\mathscr{L}}^{F_p}_{v,t}=-g^{ij}\Big\{&\big(\pa_i+{v\sqrt{t}}\overline{\theta}_{\sqrt{{vt}}}^{-1}\widetilde{\theta}_{\sqrt{v}}^{-1}\widehat{\Gamma}^{0}\lk\pa_i\rk\widetilde{\theta}_{\sqrt{v}}\overline{\theta}_{\sqrt{{vt}}}+\frac{\sqrt{t}}{p}\widehat{\Gamma}^{F_p,u}\lk\pa_i\rk\big)\\
		&\cdot\big(\pa_j+{v\sqrt{t}}\overline{\theta}_{\sqrt{{vt}}}^{-1}\widetilde{\theta}_{\sqrt{v}}^{-1}\widehat{\Gamma}^{0}\lk\pa_i\rk\widetilde{\theta}_{\sqrt{v}}\overline{\theta}_{\sqrt{{vt}}}+\frac{\sqrt{t}}{p}\widehat{\Gamma}^{F_p,u}\lk\pa_j\rk\big)\\
		+v\sqrt{t}&\big(\nabla^{TX_0}_{\pa_i}\pa_j+{v\sqrt{t}}\overline{\theta}_{{v\sqrt{t}}}^{-1}\widetilde{\theta}_{v}^{-1}\widehat{\Gamma}^{0}(\nabla^{TX_0}_{\pa_i}\pa_j)\widetilde{\theta}_{\sqrt{v}}\overline{\theta}_{\sqrt{{vt}}}+\frac{\sqrt{t}}{p}\widehat{\Gamma}^{F_p,u}(\nabla^{TX_0}_{\pa_i}\pa_j)\big)\Big\}\\
		+f_\varepsilon\Big\{&\frac{v^2t}{4}r^X-\frac{v^2t}{8}\big\langle R^{TX}(e_i,e_j)e_k,e_\ell \big\rangle c_{{vt}}(e_i)c_{{vt}}(e_j)\widehat{c}_{{v}}\lk \widehat{e}_k\rk\widehat{c}_{{v}}\lk \widehat{e}_\ell\rk\\
		&-\frac{v\sqrt{t}}{8}\big\langle R^{TX}(f_\alpha^H,f_\beta^H\big)e_k,e_\ell \big\rangle f^\alpha f^\beta\widehat{c}_{{v}}\lk \widehat{e}_k\rk\widehat{c}_{{v}}\lk \widehat{e}_\ell\rk\\
		&-\frac{v\sqrt{vt}}{4}\big\langle R^{TX}\big(e_i,f_\alpha^H\big)e_k,e_\ell \big\rangle c_{{vt}}(e_i)f^\alpha\widehat{c}_{{v}}\lk \widehat{e}_k\rk\widehat{c}_{{v}}\lk \widehat{e}_\ell\rk \\
		&-\frac{vt}{2}c_{{vt}}(e_i)c_{{vt}}(e_j)\frac{1}{4p}\omega^{F_p,2}(e_i,e_j)-\frac{1}{2}f^\alpha f^\beta\frac{1}{4p}\omega^{F_p,2}\big(f_\alpha^H,f_\beta^H\big)\\
		&-{\sqrt{vt}}c_{{vt}}(e_i)f^\alpha\frac{1}{4p}\omega^{F_p,2}\big(e_i,f_\alpha^H\big)+\frac{t}{4p^2}\bv\widehat{\omega}^{F_p}\bv^2\\
		&+\frac{vt}{2}\widehat{c}_{v}(\widehat{e}_i)\widehat{c}_{v}(\widehat{e}_j)\frac{1}{4p}\widehat{\omega}^{F_p,2}(e_i,e_j)-\sqrt{vt}f^\alpha\widehat{c}_{v}(\widehat{e}_i)\frac{1}{2p}\nabla^{\widehat{TX}\otimes F_p,u}_{f^H_\alpha}\widehat{\omega}^{F_p}(e_i)\\
		&-{vt}{c}_{{vt}}(e_i)\widehat{c}_{{v}}(\widehat{e}_j)\frac{1}{2p}\nabla^{\widehat{TX}\otimes F_p,u}_{e_i}\widehat{\omega}^{F_p}(e_j)\\
		&-{\sqrt{vt}}zc_{{vt}}(e_i)\frac{1}{2p}{\omega}^{F_p}(e_i)-zf^\alpha\frac{1}{2p}{\omega}^{F_p}\big(f^H_\alpha\big)\Big\},
	\end{split}
\end{equation}
then we get $\widehat{\mathscr{L}}^{F_p}_{1/p,t}=\widehat{\mathscr{L}}^{F_p}_{t,x_0}$. Similar to \eqref{fc31}, we have
\begin{equation}\label{fi5}
	\begin{split}
		{v\sqrt{t}}\overline{\theta}_{\sqrt{{vt}}}^{-1}\widetilde{\theta}_{\sqrt{v}}^{-1}\widehat{\Gamma}^{0}_{{v\sqrt{t}Z}}\lk\pa_i\rk\widetilde{\theta}_{\sqrt{v}}\overline{\theta}_{\sqrt{{vt}}}&={vZ_j}\Gamma'_{ij,v\sqrt{t}Z}-{vtZ_j}\Gamma''_{ij,{v\sqrt{t}Z}},
	\end{split}
\end{equation}
where the support of  $\Gamma'_{i,Z},\Gamma''_{i,Z}$ are contained in $B^{T_{x_0}X}(0,2{\varepsilon})$, and 
\begin{equation}
	\begin{split}
		\Gamma'_{ij,0}=&\frac{1}{4}\Big(\sum_{k<\ell}\big\langle R^{TX}(e_k,e_\ell)\pa_i,\pa_j\big\rangle e^k\wedge e^\ell+\sum\big\langle R^{TX}(f_\alpha,e_k)\pa_i,\pa_j\big\rangle f^\alpha\wedge e^k\\
		&+\sum_{\alpha<\beta}\big\langle R^{TX}(f_\alpha,f_\beta)\pa_i,\pa_j\big\rangle f^\alpha\wedge f^\beta\Big),\\
\Gamma''_{ij,0}=&-\frac{1}{4}\sum_{k<\ell}\big\langle R^{TX}(e_k,e_\ell)\pa_i,\pa_j\big\rangle e^k\wedge e^\ell.
	\end{split}
\end{equation}

Since the first term on the right-hand side of \eqref{fi5} is not bounded for $Z\in T_{x_0}X$, we set a new norm to make the following operators uniformly bounded for $0\leqslant v,t\leqslant 1$:
\begin{equation}
	\begin{split}
		\mathbbm{1}_{{v\sqrt{t}\lv Z\rv}\leqslant{\varepsilon}}\cdot{vZ_j}(e^i\wedge-{vt}{i}_{{e}_i}),\ \ \ \mathbbm{1}_{{v\sqrt{t}\lv Z\rv}\leqslant{\varepsilon}}\cdot(e^i\wedge-{vt}{i}_{{e}_i}),
	\end{split}
\end{equation}
where $\mathbbm{1}_{A}$ for $A\subseteq T_{x_0}X$ represents the characteristic function of $A$.
\begin{defi}
	For $k,\ell\in\mathbb{N}$, if $N_1s=\ell s$ (see \eqref{fi.3..}), we set
	\begin{equation}\label{fi8}
		\begin{split}
	\Vert s\Vert_{0,p}^2=\int_{\mathbb{R}^{n}}\bV s(Z)\bV^2(1+\vert Z\vert^2)^{m+1-\ell}dZ,\ \ 
\Vert s\Vert_{k,p}^2=\sum_{\alpha\in\mathbb{N}^{m},\lv\alpha\rv\leqslant k}\Vert\pa^\alpha s\Vert_{0,p}^2.
		\end{split}
	\end{equation}
	Let $H_{k,p}$ be the completion Sobolev space with respect to $\lV \cdot\rV_{k,p}$.
	Let $H_{-k,p}$ be the Sobolev space of negative order with the norm given by
	\begin{equation}
		\lV s\rV_{-k,p}=\sup_{0\neq s'\in H_{k,p}}\frac{\bv\li s,s'\ri_{0,p}\bv}{\lV s'\rV_{k,p}}
	\end{equation} 
\end{defi}

\begin{theo}\label{Fi27}
	For $r\in\mathbb{N}, 0\leqslant v,t\leqslant1, p\in\mathbb{N}^*$, we have
	\begin{equation}\label{fi10}
	\frac{\pa^r}{\pa v^r}\widehat{\mathscr{L}}^{F_p}_{v,t}=\mathscr{X}^{(r)}_{v,t,ij}(Z)\pa_i\pa_j+\mathscr{Y}^{(r)}_{v,t,i}(Z)\pa_i+\mathscr{Z}^{(r)}_{v,t}(Z),\ \ \text{for\ } \mathscr{X}^{(r)}_{v,t,ij}(Z)=-\frac{\pa^rg^{ij}_{v\sqrt{t}Z}}{\pa v^r},
	\end{equation}
where $\mathscr{Y}^{F_p,(r)}_{v,t,i}, \mathscr{Z}^{F_p,(r)}_{v,t}$ are in $\mathscr{C}^\infty(T_{x_0}X,\widehat{\mathbb{E}}_{p,x_0})$ with the following properties:
	\begin{enumerate}
		\item	There is $C>0$ such that the $\lV\cdot\rV_{0,p}$ norm of $\mathscr{X}^{F_p,(r)}_{v,t,ij}$  $\mathscr{Y}^{F_p,(r)}_{v,t,i}, \mathscr{Z}^{F_p,(r)}_{v,t}$ is dominated by $C\bV(1+\vert Z\vert^r)s\bV_{0,p}$ uniformly for $0\leqslant v,t\leqslant1, p\in\mathbb{N}^*$.
		
		\item Operators $\mathscr{X}^{F_p,(r)}_{0,t,ij}$ $\mathscr{Y}^{F_p,(r)}_{0,t,i}, \mathscr{Z}^{F_p,(r)}_{0,t}$ are polynomials in $Z$ and $\sqrt{t}$. Moreover, $\mathscr{X}^{F_p,(r)}_{0,t,ij}(Z)$ is a homogeneous polynomial in $Z$ of degree $r$, the degrees in $Z$ of $\mathscr{Y}^{F_p,(r)}_{0,t,i}(Z)$ and $\mathscr{Z}^{F_p,(r)}_{0,t}(Z)$ are no more that $r$.

		\item		In particular,
		\begin{equation}\label{fi12}
			\begin{split}
				{\widehat{\mathscr{L}}}^{F_p}_{0,t}=&-{\Delta}^{T_{x_0}X}-\frac{1}{4}\big\langle R^{TX}e_i,e_j \big\rangle_{x_0}\widehat{e}^i\wedge\widehat{e}^j-\frac{1}{4p}\omega^{F_p,2}_{x_0}+\frac{t}{4p^2}\big\vert\widehat{\omega}^{F_p}\big\vert^2_{x_0}\\
				&+\frac{t}{4p}\widehat{\omega}^{F_p,2}_{x_0}-\frac{\sqrt{t}}{2p}\nabla^{\widehat{TX}\otimes F_p,u}\widehat{\omega}^{F_p}_{x_0}-\frac{1}{2p}z\omega^{F_p}_{x_0},\\
			\frac{\pa}{\pa v}\widehat{\mathscr{L}}^{F_p}_{v,t}\Big|_{v=0}=&\mathscr{Y}^{F_p,(1)}_{0,t,i}(Z)\pa_i+\mathscr{Z}^{F_p,(1)}_{0,t}(Z)=-\frac{Z_j}{4}\Big\{\big\langle R^{TX}_{x_0}\pa_i,\pa_j\big\rangle\\
			&-t\big\langle \widehat{R}^{TX}_{x_0}\pa_i,\pa_j\big\rangle-\frac{t}{2p}\omega^{F_p,2}_{x_0}\big(\pa_i,\pa_j\big)\Big\}\pa_i+\mathscr{Q}^{}_{t}(Z)+\mathscr{Q}'_t,
			\end{split}
		\end{equation}
	where $\mathscr{Q}^{F_p}_{t}(Z)$ and $\mathscr{Q}_t^{F_p}{'}$ are homogeneous polynomials in $Z$ of degree $1$ and $0$ respectively with coefficients in $\widehat{\mathbb{E}}_{p,x_0}$, and
		\begin{equation}\label{fi13}
			\begin{split}
				\mathscr{Q}_t^{F_p}{'}=&\frac{t}{2}R^{\widehat{TX}}-\frac{1}{2}\widehat{R}^{TX}-\frac{t}{4}\big\langle {R}^{TX}(f_\alpha^H,e_i){e}_k,{e}_l\big\rangle \widehat{e}^k\wedge\widehat{e}^l\wedge f^\alpha\wedge{i}_{e_i}\\
				&+e^i\wedge{i}_{e_j}\frac{t}{4p}\omega^{F_p,2}(e_i,e_j)+f^\alpha\wedge{i}_{e_i}\frac{t}{4p}\omega^{F_p,2}\big(f_\alpha^H,e_i\big)+\widehat{e}^i\wedge{i}_{\widehat{e}_j}\frac{t}{4p}\widehat{\omega}^{F_p,2}(\widehat{e}_i,\widehat{e}_j)\\
				&-f^\alpha\wedge{i}_{\widehat{e}_i}\frac{\sqrt{t}}{2p}\nabla^{\widehat{TX}\otimes F_p,u}_{f^H_\alpha}\widehat{\omega}^{F_p}(e_i)-{e}^i\wedge{i}_{\widehat{e}_j}\frac{\sqrt{t}}{2p}\nabla^{\widehat{TX}\otimes F_p,u}_{e_i}\widehat{\omega}^{F_p}(e_j)\\
				&-\widehat{e}^j\wedge{i}_{e_i}\frac{t^{3/2}}{2p}\nabla^{\widehat{TX}\otimes F_p,u}_{e_i}\widehat{\omega}^{F_p}(e_j)+z{i}_{e_i}\frac{t}{2p}\omega^{F_p}(e_i).
			\end{split}
		\end{equation}
	\end{enumerate}
\end{theo}
\begin{pro}
The first statement follows from \eqref{fi4}, \eqref{fi5} and the weight in the norm \eqref{fi8}. The rest part follows from \eqref{fi4} as we obtain Theorem \ref{Fc13} from \eqref{4.2.9}.\qed
\end{pro}
By Theorem \ref{Fi27}, using the same way as we get Theorems \ref{4.11}-\ref{Fg26}, we could prove the following Theorems \ref{Fi29}-\ref{Fi38} that all hold uniformly for $0\leqslant v,t\leqslant1,p\in\mathbb{N}^*$.
\begin{theo}\label{Fi29}
	There exists $C>0$ such that
	\begin{equation}
		\begin{split}
			&\mathrm{Re}\big\langle \widehat{\mathscr{L}}^{F_p}_{v,t}s,s \big\rangle_{0,p}\geqslant C^{-1}\lV\nabla s\rV_{0,p}^2-C\Vert s\Vert_{0,p}^2,\\
			&\lv\mathrm{Im}\big\langle \widehat{\mathscr{L}}^{F_p}_{v,t}s,s \big\rangle_{0,p}\rv\leqslant C\Vert s\Vert_{1,p}\Vert s\Vert_{0,p},\\
			&\left\vert\big\langle \widehat{\mathscr{L}}^{F_p}_{v,t}s,s' \big\rangle_{0,p}\right\vert\leqslant C\Vert s\Vert_{1,p}\cdot\Vert s'\Vert_{1,p}.
		\end{split}
	\end{equation}
\end{theo}

\begin{prop}
For $k\in\mathbb{N}$, there is $C>0$ such that for any $Q_1,\cdots Q_k\in\{{\pa_i},Z_i\}_{i=1}^{m}$,
	\begin{equation}
	\Bv\big\langle\big[Q_{1},[Q_{2},\cdots[Q_{k},\widehat{\mathscr{L}}^{F_p}_{v,t}]\cdots]\big]s,s'\big\rangle_{0,p}\Bv\leqslant C\Vert s\Vert_{1,p}\Vert s'\Vert_{1,p}.
	\end{equation}
\end{prop}

\begin{theo}
	There exist $d_1,d_2> 0$ such that if $\lambda\in V_{d_1,d_2}$ where $V_{d_1,d_2}$ is given in \eqref{fe1}, the resolvent  $\blk\lambda-\widehat{\mathscr{L}}^{F_p}_{v,t}\brk^{-1}$ exists. Moreover, there is $C>0$ such that
	\begin{align}
\bV\blk\lambda-\widehat{\mathscr{L}}^{F_p}_{v,t}\brk^{-1}\bV^{0,0}_p\leqslant C,\ \ \ \ \bV\blk\lambda-\widehat{\mathscr{L}}^{F_p}_{v,t}\brk^{-1}\bV^{-1,1}_p\leqslant C\blk1+\lv\lambda\rv^2\brk.
	\end{align} 
\end{theo}

\begin{prop}
	For any $\lambda\in V_{d_1,d_2}$, $k\in\mathbb{N}$, the resolvent $\blk\lambda-\widehat{\mathscr{L}}^{F_p}_{v,t}\brk^{-1}$ maps $H_{k,p}$ to $H_{k+1,p}$. Moreover, for any multi-index $\alpha\in\mathbb{N}^m$, there exists $C_{}>0$ such that
	\begin{equation}
		\bV Z^\alpha\blk\lambda-\widehat{\mathscr{L}}^{F_p}_{v,t}\brk^{-1}s\bV_{k+1,p}\leqslant C_{}\blk1+\lv\lambda\rv^2\brk^{k+\lv\alpha\rv+1}\sum_{\alpha'\leqslant\alpha}\bV Z^{\alpha'}s\bV_{k,p}.
	\end{equation}
\end{prop}

For $({\bm{k}},{\bm{r}})\in I_{k,r}, \lambda\in V_{d_1,d_2}$, we set
\begin{equation}
	B^{\bm{k}}_{\bm{r}}\lk v,\lambda,p\rk=\blk\lambda-\widehat{\mathscr{L}}^{F_p}_{v,t}\brk^{-k_0}\frac{\pa^{r_1}\widehat{\mathscr{L}}^{F_p}_{v,t}}{\pa v^{r_1}}\blk\lambda-\widehat{\mathscr{L}}^{F_p}_{v,t}\brk^{-k_1}\cdots \frac{\pa^{r_j}\widehat{\mathscr{L}}^{F_p}_{v,t}}{\pa v^{r_j}}\blk\lambda-\widehat{\mathscr{L}}^{F_p}_{v,t}\brk^{-k_j}.
\end{equation}
Then there exist $a^{\bm{k}}_{\bm{r}}\in\mathbb{R}$ such that
\begin{equation}
	\frac{\pa^{r}}{\pa v^{r}}\blk\lambda-\widehat{\mathscr{L}}^{F_p}_{v,t}\brk^{-k}=\sum_{({\bm{k}},{\bm{r}})\in I_{k,r}}a^{\bm{k}}_{\bm{r}}B^{\bm{k}}_{\bm{r}}\lk v,\lambda,p\rk
\end{equation}

\begin{theo}
Let $\ell\in\mathbb{N}$, for any multi-index $\alpha,\alpha'\in\mathbb{N}^m$ with $\lv\alpha_i\rv\leqslant \ell$ and $k> 2(\ell+ r+1)$, if $({\bm{k}},{\bm{r}})\in I_{k,r}$, there are $C> 0$ and $j\in\mathbb{N}$ such that for any $\lambda\in V_{d_1,d_2}$,
	\begin{equation}
		\bV\pa^{\alpha} B^{\bm{k}}_{\bm{r}}\lk v,\lambda,p\rk\pa^{\alpha'}s\bV_{0,p}\leqslant C(1+{\vert}\lambda{\vert})^j \sum_{{\vert}\alpha{\vert}\leqslant r}\lV Z^\alpha s\rV_{0,p}.
	\end{equation}
\end{theo}

\begin{theo}
	For any $r\in\mathbb{N}$ and $k\in\mathbb{N}^*$, there exist $C > 0$ and $\ell\in\mathbb{N}$ such that for $\lambda\in V_{d_1,d_2}$, we have
	\begin{equation}
\begin{split}
&\bigg\Vert\Big(\frac{\pa^{r}\widehat{\mathscr{L}}^{F_p}_{v,t}}{\pa v^{r}}-\frac{\pa^{r}\widehat{\mathscr{L}}^{F_p}_{v,t}}{\pa v^{r}}\Big|_{v=0}\Big) s\bigg\Vert_{-1,p}\leqslant Cv\sum_{|\alpha|\leqslant r+1}\lV Z^\alpha s\rV_{1,p},\\
	&\bigg\Vert \Big(\frac{\pa^{r}}{\pa v^{r}}\blk\lambda-\widehat{\mathscr{L}}^{F_p}_{v,t}\brk^{-k}-\sum_{\lk{\bm{k}},{\bm{r}}\rk\in I_{k,r}}a^{\bm{k}}_{\bm{r}}B^{\bm{k}}_{\bm{r}}\lk 0,\lambda,p\rk \Big)s\bigg\Vert_{0,p}\\
	&\ \ \ \ \ \ \ \ \ \ \ \ \ \ \ \ \ \ \ \ \ \ \ \ \ \ \ \ \ \ \ \ \ \ \ \ \ \ \leqslant Cv\lk1+\vert\lambda\vert\rk^\ell\sum_{|\alpha|\leqslant 4r+1}\lV Z^\alpha s\rV_{0,p}.
\end{split}
	\end{equation}
\end{theo}

\begin{theo}
	For $k,\ell,r\in\mathbb{N}$, there is $C>0$ such that for $Z,Z'\in T_{x_0}X$ with $\lv Z\rv,\lv Z'\rv\leqslant1$, we have
	\begin{equation}
		\sup_{\lv\alpha\rv,\lv\alpha'\rv\leqslant k}\bigg\vert\frac{\pa^{{\vert}\alpha{\vert}+{\vert}\alpha'{\vert}}}{\pa Z^\alpha\pa {Z'}^{\alpha'}}\frac{\pa^{r}}{\pa v^{r}}\exp(-\widehat{\mathscr{L}}^{F_p}_{v,t})\lk Z,Z'\rk\bigg\vert_{\mathscr{C}^\ell(M)}\leqslant C.
	\end{equation}
\end{theo}

For $k$ large enough, set
\begin{equation}
	\begin{split}
		\mathscr{T}^{(r)}_{t,p}&=\frac{(k-1)!}{2\pi {i}(-1)^{k-1}r!}\int_{\pa V_{d_1,d_2}}e^{-\lambda}\sum_{({\bm{k}},{\bm{r}})\in I_{k,r}}a^{\bm{k}}_{\bm{r}}B^{\bm{k}}_{\bm{r}}\lk 0,\lambda,p\rk\de \lambda,\\
		\mathscr{T}^{(r)}_{v,t,p}&=\frac{1}{r!}\frac{\pa^{r}}{\pa v^{r}}\exp(-t\widehat{\mathscr{L}}^{F_p}_{v,t})-\mathscr{T}^{(r)}_{t,p}.
	\end{split}
\end{equation}
Then $\mathscr{T}^{(r)}_{t,p}$ and $\mathscr{T}^{(r)}_{v,t,p}$ do not depend on the choice on $k$.  We denote by $\mathscr{T}^{(r)}_{t,p}(Z,Z')$ (resp. $\mathscr{T}^{(r)}_{v,t,p}(Z,Z')$) the smooth kernel of $\mathscr{T}^{(r)}_{t,p}$ (resp. $\mathscr{T}^{(r)}_{v,t,p}$) with respect to $d v_{TX}$.

\begin{theo}
	For $r\in\mathbb{N}$, there is $C>0$ such that for $Z,Z'\in T_{x_0}X$ with $\lv Z\rv,\lv Z'\rv\leqslant1$,
	\begin{equation}
		\lv\mathscr{T}^{(r)}_{v,t,p}(Z,Z')\rv\leqslant Cv^{1/(m+1)}.
	\end{equation}
\end{theo}

\begin{theo}
	For any $r,\ell\in\mathbb{N}$, there is $C > 0$ such that
	\begin{equation}\label{fi.25}
		\Big\vert\exp(-\widehat{\mathscr{L}}^{F_p}_{v,t})(0,0)-\sum_{i=0}^r\mathscr{T}^{(i)}_{t,p}(0,0)v^i\Big\vert_{\mathscr{C}^\ell(M)}\leqslant Cv^{r+1}.
	\end{equation}
\end{theo}

For $s>0$, let $\mathcal{P}^{F_p}_s(x,y)$ be the kernel of $\exp(-s\widehat{\mathscr{L}}^{F_p}_{0,t})$ with respect to $dv_{TX}$. For multi-index ${\bm{r}}=(r_1,\cdots,r_k)_{k\in\mathbb{N}}$ where $r_1,\cdots r_k\in\mathbb{N}^*$, like \eqref{fg8}, we set
\begin{equation}\label{fi27}
	\begin{split}
	\mathscr{T}_{{\bm{r}},t,p}=
	\int_{\triangle_k}\Big(\mathcal{P}^{F_p}_{1-t_k}\frac{\pa^{r_k}\widehat{\mathscr{L}}^{F_p}_{v,t}}{\pa v^{r_k}}\Big|_{v=0}\mathcal{P}_{t_k-t_{k-1}}\cdots\frac{\pa^{r_1}\widehat{\mathscr{L}}^{F_p}_{v,t}}{\pa v^{r_1}}\Big|_{v=0}\mathcal{P}^{F_p}_{t_1}\Big)(0,0)\prod_{j=1}^kdt_j,
	\end{split}
\end{equation}
then by the Duhamel's principle \cite[\S\,2.7]{bgv}, we have
\begin{equation}\label{fi28}
	\mathscr{T}^{(i)}_{t,p}(0,0)=\sum_{\substack{ r_1+\cdots +r_k=i}}\frac{(-1)^k}{\prod_{j=1}^k r_j!}\mathscr{T}_{{\bm{r}},t,p},\ \ \mathscr{T}^{(0)}_{t,p}(0,0)=\frac{{e}^{-\sigma^{F_p}_t}}{(4\pi)^{\frac{m}{2}}}.
\end{equation}
\begin{theo}\label{Fi38}
	For $i,j,\ell\in \mathbb{N}$, there exists $C>0$ such that
	\begin{equation}\label{fi30}
		\begin{split}
			&\Bv\tro_{[j]}\big[\mathscr{T}^{(i)}_{t,p}(0,0)\big]\Bv_{\mathscr{C}^\ell(M)}\leqslant C,\\
			&\bbv p^{-n}\tro^{F_p}[\mathscr{T}^{(i)}_{t,p}(0,0)]-\sum_{k=0}^jp^{-k}\tro_{[k]}[\mathscr{T}^{(i)}_{t,p}(0,0)]\bbv_{\mathscr{C}^\ell(M)}\leqslant Cp^{-j-1}.
		\end{split}
	\end{equation}
\end{theo}
Note that there is singularity inside the integral of \eqref{fi27}, so similar to the proof of Theorem \ref{Fg26}, here we use \eqref{fg.6} when $t=1$ to get Theorem \ref{Fi38}.

\subsection{The proof of Theorem \ref{F1}}\label{Fh}

Now we begin to prove the main theorem of this section. First, we deal with the case of $t\geqslant1$. Recall the symbols $[\cdot]^z$ in \eqref{bj26}, $[\cdot]^{\text{max}}$ in \eqref{fc57}  and $\int^{\widehat{B}}$ in \eqref{ab1}, then we set
\begin{equation}\label{fg3}
\mathscr{I}_{x_0}=(4\pi)^{\frac{m}{2}}p^{-n}\tro^{F_p}\Big[\theta_{\sqrt{t}}^{-1}\int^{\widehat{B}}\big[\exp(-t\widehat{\mathscr{L}}^{F_p}_{x_0})(0,0)\big]^{\text{max}}\Big]^{z}.
\end{equation}
By \eqref{fa5}, \eqref{fb.82} and \eqref{4.2.8}, we have
\begin{equation}\label{fg4}
\begin{split}
	\bbv p^{-n}\frac{1}{\sqrt{p}}\psi_{1/\sqrt{p}}h\big( A',g^{\Omega^\bullet(X,F_p)}_{4t/p^2}\big)-\int_X\mathscr{I}\bbv_{\mathscr{C}^\ell(S)}\leqslant C{e}^{\gamma t-\sqrt{c_1c_2}c\varepsilon p}.
\end{split}
\end{equation}
Taking $v=1/p$ in \eqref{4.4.50}, we get
\begin{equation}\label{fg5}
\Big\vert\exp(-t\widehat{\mathscr{L}}^{F_p}_{x_0})(0,0)-\sum_{i=0}^rp^{-i}\mathscr{S}^{(i)}_{t,p}(0,0)\Big\vert_{\mathscr{C}^\ell(M)}\leqslant Cp^{-r-1}{e}^{\gamma t}.
\end{equation}
Put
\begin{equation}\label{fh6}
\bm{c}_i(t)=\sqrt{{2\pi i}}\varphi\bigg[(4\pi)^{\frac{m}{2}}\theta_{\sqrt{t}}^{-1}\int^{\widehat{B}}\sum_{0\leqslant j\leqslant i}\tro_{[i-j]}\big[\mathscr{S}^{(j)}_{t,p}(0,0)\big]^{\text{max}}\bigg]^z.
\end{equation}
By \eqref{fg24}, \eqref{fg3}, \eqref{fg4} and \eqref{fg5}, we get \eqref{fa1} and \eqref{fa3} for $t\geqslant1$.

Likewise, if we take
\begin{equation}\label{fh7}
c_i(t)=\sqrt{{2\pi i}}\varphi\bigg[(4\pi)^{\frac{m}{2}}\int^{\widehat{B}}\sum_{0\leqslant j\leqslant i}\tro_{[i-j]}\big[\mathscr{T}^{(j)}_{t,p}(0,0)\big]^{\text{max}}\bigg]^z,
\end{equation}
then by \eqref{fa5}, \eqref{fb.82}, \eqref{fi.2}, \eqref{fi.25} and \eqref{fi30}, we see that \eqref{fa1} and \eqref{fa3} hold for $0<t\leqslant 1$. Note that $c_i(t)\neq \bm{c}_i(t)$ in general, while when separating them with respect to \eqref{ba4}, they contain the same component in top degree of $\Lambda(T^*X)$, this gives
\begin{equation}\label{fh8}
\int_X\bm{c}_i(t)= \int_X{c}_i(t)
\end{equation}
for $i\in\mathbb{N}$ and $t>0$. Hence for $t>0$, we can always define $c_i(t)$ by \eqref{fh7}. This completes the proof.

\section{The full asymptotics of the analytic torsion forms}\label{sG}

In this section, when $\widehat{\vartheta}^L$ is nondegenerate, we obtain the full asymptotics of the analytic torsion forms $\mathcal{T}\blk T^HM,g^{TX},\nabla^{F_p},g^{F_p}\brk$ as $p\rightarrow+\infty$.

This section is organized as follows. In \cref{Ga}, we give the existence of the full asymptotic expansion of torsion forms. The proof is divided into two key steps, involving large and small values of the parameter $t>0$. In \cref{Gc}, we prove that the $\Gamma$-torsion forms share the same asymptotic behavior as the torsion forms. In \cref{Gb}, we prove that $\mathscr{T}_{1,t,p}(0,0)$ does not contribute to the first two terms in the expansion of torsion forms.

\subsection{The main result: a precise version of Theorem \ref{1.1}}\label{Ga}

For the enlarged fibration $\widetilde{M}\rightarrow \widetilde{S}$ as in \cref{Bi}, by \eqref{bh13} we have
\begin{equation}\label{ga1}
	\begin{split}
h\big(\widetilde{A}',g^{\Omega^\bullet(X,F_p)}_{4t/p^2}\big)\big|_{s=1}=h\big({A}',g^{\Omega^\bullet(X,F_p)}_{4t/p^2}\big)+ds\wedge h^{\wedge}\big({A}',g^{\Omega^\bullet(X,F_p)}_{4t/p^2}\big).
	\end{split}
\end{equation}
By Theorem \ref{F1}, there exist smooth forms $\{\widetilde{c}_i(t)\}_{i\in\mathbb{N}}$ as in \eqref{fa1} for $h\big(\widetilde{A}',g^{\Omega^\bullet(X,F_p)}_{4t/p^2}\big)$, then for some smooth forms $\{d_i(t)\}_{i\in\mathbb{N}}$, for $t>0$ we have
\begin{equation}\label{ga2}
\widetilde{c}_i(t)=c_i(t)+ds\wedge d_i(t).
\end{equation}
According to \eqref{fa1}, \eqref{ga1} and \eqref{ga2}, for  $k,\ell\in\mathbb{N}$, there is $C>0$ such that for any $t>0$,
\begin{equation}\label{ga3}
	\begin{split}
		\bbv p^{-n-1}\psi_{1/\sqrt{p}}h^{\wedge}\big( A',g^{\Omega^\bullet(X,F_p)}_{4t/p^2}\big)-\sum_{i=0}^k\int_Xd_i(t)p^{-i}\bbv_{\mathscr{C}^\ell(S)}\leqslant C{e}^{-(a-{\gamma})t}p^{-k-1},
	\end{split}
\end{equation}
moreover, 
\begin{equation}\label{ga4}
	\bbv\int_Xd_k(t)\bbv_{\mathscr{C}^\ell(S)}\leqslant C{e}^{-(a-{\gamma})t}.
\end{equation}

By \eqref{bf4}, the restriction on $M\times \{1\}$ of the operator $\widehat{\mathscr{L}}^{F_p}_{v,t}$ on $\widetilde{M}$ is exactly the operator in \eqref{fi4} by counting also the direction $\frac{\pa}{\pa s}$ in the term $-\sqrt{vt}f^\alpha\widehat{c}_{v}(\widehat{e}_i)\frac{1}{2p}\nabla^{\widehat{TX}\otimes F_p,u}_{f^H_\alpha}\widehat{\omega}^{F_p}(e_i)$, so we only need to add the following extra term in \eqref{fi4}:
\begin{equation}\label{ga5}
	\begin{split}
\frac{1}{2}\sqrt{vt}ds\wedge\widehat{c}_{{v}}(\widehat{e}_i)\frac{1}{2p}\widehat{\omega}^{F_p}(e_i).
	\end{split}
\end{equation}
By \eqref{ga5}, the degrees of $\sqrt{t}$ in all the terms of \eqref{fi4} that contain $ds$ are no less than $1$. Therefore, when $0<t\leqslant 1$, \eqref{ga3} and \eqref{ga4} can be replaced by 
\begin{equation}\label{ga.5}
	\begin{split}
		&\bbv {p^{-n-1}}{\psi_{1/\sqrt{p}}}h^{\wedge}\big( A',g^{\Omega^\bullet(X,F_p)}_{4t/p^2}\big)-\sum_{i=0}^k\int_Xd_i(t)p^{-i}\bbv_{\mathscr{C}^\ell(S)}\leqslant{C\sqrt{t}}{p^{-k-1}},\\ &\bbv\int_Xd_k(t)\bbv_{\mathscr{C}^\ell(S)}\leqslant C\sqrt{t}.
\end{split}
\end{equation}

Now we state and prove the main result of this article.
\begin{theo}\label{4.21}
	If $\widehat{\vartheta}^L$ is nondegenerate, then for any $i\in\mathbb{N}$, the integral
	\begin{equation}\label{ga6}
		W_i^{L,\xi}=-\sqrt{{2\pi i}}\varphi\int_{0}^{+\infty}{d_i(t)}\frac{d t}{t}\in \Omega^\bullet(M,o(TX))
	\end{equation}
	is well defined.
	Moreover, for $k,\ell\in\mathbb{N}$, there is $C>0$ such that as $p\rightarrow+\infty$, we have
	\begin{equation}\label{ga7}
	\bbv p^{-n-1} \psi_{1/\sqrt{p}}\mathcal{T}\blk T^HM,g^{TX},\nabla^{F_p},g^{F_p}\brk-\sum_{i=0}^k\int_XW^{L,\xi}_ip^{-i}\bbv_{\mathscr{C}^\ell(S)}\leqslant Cp^{-k-1}.
	\end{equation}
Also, for differential forms $\gamma
	_i$ defined in \eqref{dd1}, we have
	\begin{equation}\label{ga8}
	d W_i^{L,\xi}=\int_X\big[e\big(TX,\nabla^{TX}\big)\gamma_i\big].
	\end{equation}
\end{theo}

\begin{pro}
	 By Theorem \ref{De8}, we have $H^\bullet(X,F_p)=0$ for $p\in\mathbb{N}^*$ large enough, then by \eqref{bi1},
	\begin{equation}\label{ga.10}
		\mathcal{T}\blk T^HM,g^{TX},\nabla^{F_p},g^{F_p}\brk=-\int_0^{+\infty}h^{\wedge}\big( A',g^{\Omega^\bullet(X,F_p)}_{t}\big)\frac{\de t}{t}.
	\end{equation}
	Using the change of variable $t\rightarrow \frac{4t}{p^2}$ on the right hand side of \eqref{ga.10}, the integral becomes
	\begin{equation}\label{ga10}
		\begin{split}
			\mathcal{T}\blk T^HM,g^{TX},\nabla^{F_p},g^{F_p}\brk=-\int_0^{+\infty}h^{\wedge}\big( A',g^{\Omega^\bullet(X,F_p)}_{4t/p^2}\big)\frac{\de t}{t}.
		\end{split}
	\end{equation}

By \eqref{ga4} and \eqref{ga.5}, the form $W_i^{L,\xi}$ in \eqref{ga6} is well-defined, and we have
	\begin{equation}
		\begin{split}
			p^{k+1}\bbv p^{-n-1} \psi_{1/\sqrt{p}}\mathcal{T}\blk T^HM,g^{TX},\nabla^{F_p},g^{F_p}\brk-\sum_{i=0}^k\int_XW_i^{L,\xi}p^{-i}\bbv_{\mathscr{C}^\ell(S)}&\\
			\leqslant \int_0^1C\sqrt{t}\frac{\de t}{t}+\int_1^{+\infty}C{e}^{-at}\frac{\de t}{t}\leqslant C&,
		\end{split}
	\end{equation}
from which we get \eqref{ga7}. And \eqref{ga8} follows directly from \eqref{bi2}, \eqref{dd1} and \eqref{ga7}.\qed
\end{pro}

\subsection{The asymptotics of the $\Gamma$-torsion forms}\label{Gc}

Let $\widehat{M}\to S$ be a smooth fibre bundle with fibre $\widehat{X}$, and let $\Gamma$ be a discrete group acting fibrewise freely and properly discontinuously on $\widehat{M}$ such that the $\Gamma\backslash \widehat{M}=M$. Let $\widehat{\pi}\colon\widehat{M}\to M$ be the obvious projection.

Let $F$ be a Hermitian vector bundle on $M$. If $Q(\widehat{x},\widehat{x}')$ be a continuous fibrewise kernel acting on $\Omega^\bullet(\widehat{M},\widehat{\pi}^*F)$ that commute with $\Gamma$, we can define its Von Neumann $\Gamma$-supertrace. Note that $\tro_s[Q(x,x)]$ is $\Gamma$-invariant that it descends to $M$, we define $\tro_s^\Gamma[Q]$ as the fibrewise integral of $\tro_s[Q(x,x)]$ on $M$: for each $b\in S$, let $\widehat{X}_b$ and $X_b$ be the corresponding fibres in $\widehat{M}$ and $M$, then we have an isomorphism $\Gamma\backslash\widehat{X}_b\cong X_b$, let $\mathcal{X}_b\subset \widehat{X}_b$ be an associated fundamental domain, we have

Let $F$ be a Hermitian vector bundle on $M$. If $Q$ is an operator acting on $\mathscr{C}^\infty(\widehat{M},\widehat{\pi}^*F)$ with a smooth fibrewise kernel $Q(\widehat{x},\widehat{x}')$ for $(\widehat{x},\widehat{x}')\in \widehat{M}\times_S\widehat{M}$. If $Q$ commutes with $\Gamma$, then $\tro[Q(x,x)]$ is $\Gamma$-invariant and it descends to a function on $M$. For each $b\in S$, let $\widehat{X}_b$ and $X_b$ be the corresponding fibres in $\widehat{M}$ and $M$, then we have an isomorphism $\Gamma\backslash\widehat{X}_b\cong X_b$, and the Von Neumann $\Gamma$-trace of $Q$ is given by
\begin{equation}
	\tro^\Gamma[Q_b]=\int_{{X}_b}\tro[Q(x,x)]dv_{{X}_b}(x).
\end{equation}

We will now apply the formalism of the previous sections. Note that we can lift geometric data $(T^HM,g^{TX},F_p,g^{F_p})_{{p\in{\mathbb{N}^*}}}$ on $M$ to $(T^H\widehat{M},g^{\widehat{\pi}^*TX},\widehat{\pi}^*F_p,g^{\widehat{\pi}^*F_p})_{{p\in{\mathbb{N}^*}}}$ on $\widehat{M}$. 

For $t>0$, we can define $h^{\wedge,\Gamma}(A',g_t^{\Omega(X,\widehat{\pi}^*F_p)})\in\Omega^\bullet(S)$ by the same formula as in \eqref{bh13}, by replacing the supertrace by the corresponding $\Gamma$-supertrace $\tro_s^\Gamma$. We assume that the nondegeneracy condition \eqref{de1} holds. As observed by Bismut-Ma-Zhang \cite[\S\,6.6]{bmz17}, when the nondegeneracy condition \eqref{de1} holds, the spectral gap property \eqref{de2} is still valid over $\widehat{M}$, therefore, for ${p\in{\mathbb{N}^*}}$ large enough, as $t\to+\infty$, $\vert h^{\wedge,\Gamma}( A',g^{\Omega^\bullet(X,\widehat{\pi}^*F_p)}_{t})\vert=O({e}^{-p^2t})$, then we can define the following  $\Gamma$-torsion form similar to \eqref{bi1} and \eqref{ga.10}:
\begin{equation}\label{gc5}
	\mathcal{T}^\Gamma\blk T^HM,g^{TX},\nabla^{F_p},g^{F_p}\brk=-\int_0^{+\infty}h^{\wedge,\Gamma}\big( A',g^{\Omega^\bullet(X,\widehat{\pi}^*F_p)}_{t}\big)\frac{\de t}{t}\in\Omega^\text{even}(S).
\end{equation}
In particular, if $\widehat{X}$ is the universal covering of $X$, the $\Gamma$-torsion is also called the $L^2$-torsion, denoted by $\mathcal{T}_{L^2}\blk T^HM,g^{TX},\nabla^{F_p},g^{F_p}\brk$. For the definition of $\Gamma$-torsion forms of general bundles whose Novikov-Shubin invariant is positive, see \cite{gammat} for more details.

\begin{theo}\label{Gc5}
	Under the nondegeneracy condition \eqref{de1}, there is $a>0$ such that for $\ell\in\mathbb{N}$, there is $C>0$ such that as $p\rightarrow+\infty$, 
	\begin{equation}\label{gc.15.}
		\begin{split}
			\Bv\mathcal{T}^\Gamma\blk T^HM,g^{TX},\nabla^{F_p},g^{F_p}\brk-\mathcal{T}\blk T^HM,g^{TX},\nabla^{F_p},g^{F_p}\brk\Bv_{\mathscr{C}^\ell(S)}\leqslant C{e}^{-a{p}}.
		\end{split}
	\end{equation}
	In particular, the $\Gamma$-torsion forms verify the same asymptotic expansion as in \eqref{ga7}:
	\begin{equation}\label{gc4}
		\bbv p^{-n-1} \psi_{1/\sqrt{p}}\mathcal{T}^\Gamma\blk T^HM,g^{TX},\nabla^{F_p},g^{F_p}\brk-\sum_{i=0}^k\int_XW_i^{L,\xi}p^{-i}\bbv_{\mathscr{C}^\ell(S)}\leqslant Cp^{-k-1}.
	\end{equation}
\end{theo}
\begin{pro}
	For $\widehat{x}_0\in\widehat{M}$, set $\mathscr{M}^{F_p}_{\widehat{x}_0}=\widehat{\pi}^*{\mathscr{M}}^{F_p}_{\widehat{\pi}(\widehat{x}_0)}$ where ${\mathscr{M}}^{F_p}_{\widehat{\pi}(\widehat{x}_0)}$ is defined in \eqref{fc10}, we clearly have
	\begin{equation}\label{gc6}
		\exp(-t{\mathscr{M}}^{\widehat{\pi}^*F_p}_{\widehat{x}_0})(0,0)=\widehat{\pi}^*\big[\exp(-t{\mathscr{M}}^{F_p}_{\widehat{\pi}(\widehat{x}_0)})(0,0)\big].
	\end{equation}
	Since the spectral estimate \eqref{de2} holds over $\widehat{M}$, by \eqref{fa5} and \eqref{fb.82}, we see that for $\ell\in\mathbb{N}$, there is $C>0$ such that for ${p\in{\mathbb{N}^*}}$, we have 
	\begin{equation}\label{gc7}
		\begin{split}
			\Bv\trs\big[\exp(-t\mathcal{M}^{\widehat{\pi}^*F_p})(\widehat{x}_0,\widehat{x}_0)-\exp(-t{\mathscr{M}}^{\widehat{\pi}^*F_p}_{\widehat{x}_0})(0,0)\big]\Bv_{\mathscr{C}^{\ell}(\widehat{M})}\leqslant C{e}^{-(a-\gamma) t-\frac{c_2}{2t}-\sqrt{\gamma c_2}p}.
		\end{split}
	\end{equation}
	By \eqref{fb.82}, \eqref{gc6} and \eqref{gc7}, we obtain
	\begin{equation}\label{gc8}
		\begin{split}
			\Bv\trs\big[\exp(-t\mathcal{M}^{\widehat{\pi}^*F_p})(\widehat{x}_0,\widehat{x}_0)\big]-\trs\big[\exp(-t\mathcal{M}^{F_p})(\widehat{\pi}(\widehat{x}_0),\widehat{\pi}(\widehat{x}_0))\big]\Bv_{\mathscr{C}^{\ell}(\widehat{M})}\\	\leqslant C{e}^{-(a-\gamma) t-\frac{c_2}{2t}-\sqrt{\gamma c_2}p}.
		\end{split}
	\end{equation}
	Similar to \eqref{fa5}, we set
	\begin{equation}\label{gc9}
		\frac{1}{\sqrt{p}}\psi_{1/\sqrt{p}}h^\Gamma\big( A',g^{\Omega^\bullet(X,F_p)}_{4t/p^2}\big)=\sqrt{{2\pi i}}\varphi\bigg[\theta_{\sqrt{t}}^{-1}\tro_s^{\Gamma}\big[\exp(-t\mathcal{M}^{\widehat{\pi}^*F_p})\big]\bigg]^z.
	\end{equation}
	By \eqref{fa5}, \eqref{gc8} and \eqref{gc9}, we have
	\begin{equation}\label{gc10}
		\begin{split}
			\frac{1}{\sqrt{p}}\psi_{1/\sqrt{p}}\lv h^\Gamma\big( A',g^{\Omega^\bullet(X,F_p)}_{4t/p^2}\big)-h\big( A',g^{\Omega^\bullet(X,F_p)}_{4t/p^2}\big)\rv_{\mathscr{C}^{\ell}(S)}\leqslant C{e}^{-(a-\gamma) t-\frac{c_2}{2t}-\sqrt{\gamma c_2}p}.
		\end{split}
	\end{equation}
	As in the argument before Theorem \ref{Bh17}, we apply \eqref{gc10} to the enlarged manifolds $\widetilde{\widehat{M}}\to\widetilde{\widehat{S}}$ and $\widetilde{{M}}\to\widetilde{{S}}$, when restricting to $s=1$, like \eqref{ga1}, we get
	\begin{equation}\label{gc11}
		\begin{split}
			p^{-1}\psi_{1/\sqrt{p}}\lv h^{\wedge,\Gamma}\big( A',g^{\Omega^\bullet(X,F_p)}_{4t/p^2}\big)-h^{\wedge}\big( A',g^{\Omega^\bullet(X,F_p)}_{4t/p^2}\big)\rv_{\mathscr{C}^{\ell}(S)}\leqslant C{e}^{-(a-\gamma) t-\frac{c_2}{2t}-\sqrt{\gamma c_2}p}.
		\end{split}
	\end{equation}
	From \eqref{ga10} and \eqref{gc5}, we see that for some $c>0$,
	\begin{equation}\label{gc12}
		\begin{split}
			\Bv\psi_{1/\sqrt{p}}\mathcal{T}^\Gamma\blk T^HM,g^{TX},\nabla^{F_p},g^{F_p}\brk-\psi_{1/\sqrt{p}}\mathcal{T}\blk T^HM,g^{TX},\nabla^{F_p},g^{F_p}\brk\Bv_{\mathscr{C}^\ell(S)}\\
			\leqslant\int_0^{+\infty} Cp{e}^{-(a-\gamma) t-\frac{c_2}{2t}-\sqrt{\gamma c_2}p}\frac{dt}{t}\leqslant C{e}^{-c{p}},
		\end{split}
	\end{equation}
which gives \eqref{gc.15.}. By \eqref{ga7} and \eqref{gc.15.}, we get \eqref{gc4}.\qed
\end{pro}

\subsection{A reduction for $W_1$}\label{Gb}

In this subsection, we give a reduction formula for $W_0,W_1$, under the assumption that $\dim_\mathbb{C}\xi=1$ and $F=\mathbb{C}$ (mainly used in \eqref{gb.10}). By \eqref{ga2} and \eqref{ga6}, let us compute $c_0(t)$ and $c_1(t)$ first.

By \eqref{fh7}, we get
\begin{equation}\label{gb1}
	\begin{split}
\int_X\big(c_0(t)+p^{-1}c_1(t)\big)=(4\pi)^{\frac{m}{2}}\int_X\bigg[\int^{\widehat{B}}\big(\tro_{[0]}+p^{-1}\tro_{[1]}\big)\big[\mathscr{T}^{(0)}_{t,p}(0,0)\big]^{\text{max}}&\\
+p^{-1}\tro_{[0]}\big[\mathscr{T}^{(1)}_{t,p}(0,0)\big]^{\text{max}}\bigg]^{z}&.
	\end{split}
\end{equation}

\begin{lemma}
We have
\begin{equation}\label{gb.2}
\int_X\bigg[\int^{\widehat{B}}\tro_{[0]}\big[\mathscr{T}^{(1)}_{t,p}(0,0)\big]^{\mathrm{max}}\bigg]^z=0.
\end{equation}
\end{lemma}

\begin{pro}
Put
\begin{equation}\label{gb2}
	\begin{split}
		\sigma_t^{F_p}=&-\frac{1}{4}\big\langle R^{TX}e_i,e_j \big\rangle\widehat{e}^i\wedge\widehat{e}^j-\frac{1}{4p}\omega^{F_p,2}+\frac{t}{4p^2}\big\vert\widehat{\omega}^{F_p}\big\vert^2+\frac{t}{4p}\widehat{\omega}^{F_p,2}\\
		&-\frac{\sqrt{t}}{2p}\nabla^{\widehat{TX}\otimes F_p,u}\widehat{\omega}^{F_p}-\frac{1}{2p}z\omega^{F_p}.
	\end{split}
\end{equation}
By \eqref{fi27} and \eqref{fi28}, we get
\begin{equation}\label{gb3}
	\begin{split}
		&\tro_{[0]}\big[\mathscr{T}^{(1)}_{t,p}(0,0)\big]=-\tro_{[0]}\bigg[\int_{0}^1\int_{TX}\mathcal{P}^{F_p}_{1-{t_1}}(0,Z)\frac{\pa {\widehat{\mathscr{L}}}^{F_p}_{v,{t},Z}}{\pa v}\bigg|_{v=0}\mathcal{P}^{F_p}_{{t_1}}(Z,0)dZd{t_1}\bigg]\\
		&=-\tro_{[0]}\bigg[\int_{0}^1\int_{TX}(16\pi^2 (1-{t_1})t_1)^{-\frac{m}{2}}{e}^{-\frac{\lv Z\rv^2}{4(1-{t_1})}-(1-t_1)\sigma_{t}^{F_p}}\frac{\pa {\widehat{\mathscr{L}}}^{F_p}_{v,{t},Z}}{\pa v}\bigg|_{v=0}{e}^{-\frac{\lv Z\rv^2}{4{t_1}}-t_1\sigma_{t}^{F_p}}dZd{t_1}\bigg].
	\end{split}
\end{equation}
For $1\leqslant i,j\leqslant m$, we set 
\begin{equation}\label{gb4}
	a_{t,ij}=Z_iZ_j\Big\{\big\langle R^{TX}_{x_0}\pa_i,\pa_j\big\rangle-t\big\langle \widehat{R}^{TX}_{x_0}\pa_i,\pa_j\big\rangle-\frac{t}{2p}\omega^{F_p,2}_{x_0}\big(\pa_i,\pa_j\big)\Big\},
\end{equation}
which is anti-symmetric for $i,j$. By \eqref{fi10}, \eqref{fi12} and \eqref{gb4} , we get 
\begin{equation}
	\sum_i{e}^{-\frac{\lv Z\rv^2}{4(1-t_1)}}\mathscr{Y}^{(1)}_{0,t,i}\pa_i{e}^{-\frac{\lv Z\rv^2}{4t_1}}=-\frac{1}{2t_1}\sum_{i,j}{e}^{-\frac{\lv Z\rv^2}{4(1-t_1)}}a_{t,ij}{e}^{-\frac{\lv Z\rv^2}{4t_1}}=0.
\end{equation}
By \eqref{fi12}, $\mathscr{Q}^{F_p}_t(Z)$ is a homogeneous polynomial of $Z$ with degree $1$, so we have
\begin{equation}
	\int_{T_{x_0}X}{e}^{-\frac{\lv Z\rv^2}{4(1-t_1)}}\mathscr{Q}^{F_p}_t(Z){e}^{-\frac{\lv Z\rv^2}{4t_1}} dZ=0.
\end{equation}
Now we consider the term $\mathscr{Q}_t^{F_p}{'}$, we claim that for $0\leqslant t_1\leqslant 1$, we have
\begin{equation}\label{gb.8}
	\begin{split}
\int_{X}\bigg[\int^{\widehat{B}}\tro_{[0]}\Big[{e}^{-(1-t_1)\sigma_{t}^{F_p}}\mathscr{Q}_t^{F_p}{'}{e}^{-t_1\sigma_{t}^{F_p}}\Big]^{\text{max}}\bigg]^z=0.
\end{split}
\end{equation}
Recall $\mathscr{Q}^{F_p}{'}$, $\sigma^{F_p}$ and $\overline{\theta}_a$ given in \eqref{fc29}, \eqref{fg6}, and \eqref{fi.3..} respectively, then
\begin{equation}
t\overline{\theta}_{\sqrt{t}}^{-1}\sigma^{F_p}\overline{\theta}_{\sqrt{t}}=\sigma_{t}^{F_p},\ \  t\overline{\theta}_{\sqrt{t}}^{-1}\mathscr{Q}^{F_p}{'}\overline{\theta}_{\sqrt{t}}=\mathscr{Q}_t^{F_p}{'},
\end{equation}
which gives
\begin{equation}
\begin{split}
t^{-\frac{m}{2}}\int_X\bigg[\theta^{-1}_{\sqrt{t}}\int^{\widehat{B}}	\tro_{[0]}&\Big[{e}^{-(t-t_1t)\sigma^{F_p}}t\mathscr{Q}^{F_p}{'}{e}^{-t_1t\sigma^{F_p}}\Big]^{\text{max}}\bigg]^z\\
&=\int_X\bigg[\int^{\widehat{B}}	\tro_{[0]}\Big[{e}^{-(1-t_1)\sigma_{t}^{F_p}}\mathscr{Q}_t^{F_p}{'}{e}^{-t_1\sigma_{t}^{F_p}}\Big]^{\text{max}}\bigg]^z.
\end{split}
\end{equation}
Therefore, instead of \eqref{gb.8}, we only need to prove that for $0\leqslant t'\leqslant t$, we have
\begin{equation}\label{gb.9}
	\begin{split}
\int_X\bigg[\int^{\widehat{B}}	\tro_{[0]}\Big[{e}^{-(t-t')\sigma^{F_p}}t\mathscr{Q}^{F_p}{'}{e}^{-t'\sigma^{F_p}}\Big]^{\text{max}}\bigg]^z=0.
	\end{split}
\end{equation}
The algebra of Toeplitz operators is \emph{non-commutative}, while in \eqref{gb.9}, we only need to get the first term of asymptotic traces, by Theorem \ref{Cb3}, when $\dim_\mathbb{C}\xi=1$, we can formally treat them as \emph{commutative}: for $k\in\mathbb{N}^*$ and Toeplitz operators $T_1,\cdots,T_k$, if $i_1,\cdots,i_k$ is any permutation of  $\{1,\cdots,k\}$, we have $\tro_{[0]}[T_1\cdots T_k]=\tro_{[0]}[T_{i_1}\cdots T_{i_k}]$. Then, by \eqref{fc24} and \eqref{fc29},  $\mathscr{Q}^{F_p}{'}$ works like an odd derivation: for Toeplitz operator values differential forms $\{\omega_\ell\}_{\ell=0,1}$ on $M$, we have
\begin{equation}
\tro_{[0]}\big[\mathscr{Q}^{F_p}{'}(\omega_0\wedge\omega_1)\big]=\tro_{[0]}\big[(\mathscr{Q}^{F_p}{'}\omega_0)\wedge\omega_1\big]+(-1)^{\deg\omega_0}\tro_{[0]}\big[\omega_0\wedge(\mathscr{Q}^{F_p}{'}\omega_1)\big],
\end{equation}
therefore, we get
\begin{equation}\label{gb.10}
	\begin{split}
		&\tro_{[0]}\big[{e}^{-(t-t')\sigma^{F_p}}t\mathscr{Q}^{F_p}{'}{e}^{-t'\sigma^{F_p}}\big]=\tro_{[0]}\big[{e}^{-(t-t')\sigma^{F_p}}\cdot t\mathscr{Q}^{F_p}{'}(-t'\sigma^{F_p})\cdot{e}^{-t'\sigma^{F_p}}\big]\\
		&=\tro_{[0]}\big[t\mathscr{Q}^{F_p}{'}(-t'\sigma^{F_p})\cdot{e}^{-t\sigma^{F_p}}\big]={t'}\tro_{[0]}\Big[\mathscr{Q}^{F_p}{'}\big({e}^{-t\sigma^{F_p}}\big)\Big].
	\end{split}
\end{equation}
In \eqref{gb.9}, we only need terms of top degree in $\Lambda(T^*X)\hti\Lambda(\widehat{T^*X})$. However, by \eqref{fc29}, $\mathscr{Q}'{e}^{-t\sigma^{F_p}}$ can never be of top degree in $\Lambda(T^*X)\hti\Lambda(\widehat{T^*X})$. We finish the proof.\qed
\end{pro}

For any $a\in \mathbb{R}[z]\hti\mathbb{R}[ds]\hti\Lambda(T^*M)$ with the form $a=\alpha_1+z\alpha_2+ds\alpha_3+z\wedge ds\wedge\alpha_4$ where $\alpha_i\in \Lambda(T^*M)$, we set $[a]^{z\wedge ds}=\alpha_4$. Put
\begin{equation}\label{gb11}
\begin{split}
{\eta}_t^{F_p}=&{\sigma}_t^{F_p}+\frac{\sqrt{t}}{2}ds\frac{1}{2p}\widehat{\omega}^{F_p}=-\frac{1}{4}\big\langle R^{TX}e_i,e_j \big\rangle\widehat{e}^i\wedge\widehat{e}^j-\frac{1}{4p}\omega^{F_p,2}+\frac{t}{4p^2}\widehat{\omega}^{F_p}(e_i)^2\\
&+\frac{t}{4p}\widehat{\omega}^{F_p,2}-\frac{\sqrt{t}}{2p}\nabla^{\widehat{TX}\otimes F_p,u}\widehat{\omega}^{F_p}-\frac{1}{2p}z\omega^{F_p}+\frac{\sqrt{t}}{2}ds\frac{1}{2p}\widehat{\omega}^{F_p}.
\end{split}
\end{equation}

\begin{prop}
For $i=0,1$, we have
\begin{equation}\label{gb.15}
W_i^{L,\xi}=\int_0^{+\infty}\Big[\int^{\widehat{B}}\tro_{[i]}\big[\exp(-{\eta}_t^{F_p})\big]^{\mathrm{max}}\Big]^{z\wedge ds}\frac{dt}{t}.
\end{equation}
\end{prop}

\begin{pro}
By \eqref{gb.2}, $\mathscr{T}^{(1)}_{t,p}(0,0)$ does \emph{not} contribute to the right hand side of \eqref{gb1}. By \eqref{fi12}, \eqref{fi27}, \eqref{fi28}, \eqref{gb1}  and \eqref{gb2}, we get
	\begin{equation}\label{gb.11}
		\begin{split}
c_0(t)+p^{-1}c_1(t)&=(4\pi)^{\frac{m}{2}}\bigg[\int^{\widehat{B}}\big(\tro_{[0]}+p^{-1}\tro_{[1]}\big)\big[\mathscr{T}^{(0)}_{t,p}(0,0)\big]^{\text{max}}\bigg]^{z}\\
&=\bigg[\int^{\widehat{B}}\big(\tro_{[0]}+p^{-1}\tro_{[1]}\big)\big[\exp(-\sigma_t^{F_p})\big]^{\text{max}}\bigg]^{z}.
		\end{split}
	\end{equation}
	By \eqref{ga2}, \eqref{ga5}, \eqref{gb11} and \eqref{gb.11},  we have
	\begin{equation}\label{gb12}
	d_0(t)+p^{-1}d_1(t)=\bigg[\int^{\widehat{B}}\big(\tro_{[0]}+p^{-1}\tro_{[1]}\big)\big[\exp(-{\eta}_t^{F_p})\big]^{\text{max}}\bigg]^{z\wedge ds}.
	\end{equation}
By \eqref{ga6} and \eqref{gb12}, we obtain \eqref{gb.15}.\qed
\end{pro}

\begin{remark}
By Theorem \ref{Cb7}, \eqref{gb.15} for $i=0$ is just the formula for $W_0$ given by Bismut-Ma-Zhang \cite[(9.89)]{bmz17}. By \eqref{gb.15}, we can get an explicit formula for $W_1$ for $i=0$ and the product formula of Toeplitz operators \eqref{cb24} for higher order terms (see \cite{mm07}). While direct expansion is very complicated and less illuminating. In the next section, when $G$ is a reductive Lie group, we will rewrite \eqref{gb12} in a concise and clear form. 
\end{remark}

\section{The case when the structure group of $P_G$ is reductive}\label{sH}

In this section, we give an explicit formula for $W^{L,\xi}_1$ in \eqref{ga7} when the structure group $G$ of $P_G$ is a reductive linear Lie group, $\dim_\mathbb{C}\xi=1$ and $F=\mathbb{C}$ trivial. 

This section is organized as follows. In \cref{Ha}, we introduce the reductive Lie group $G$, its compact form $U$ and complexification $G_\mathbb{C}$. In \cref{Hb}, we introduce the reduction of the flat principle bundle $P_G\rightarrow M$. In \cref{Hc}, we express certain Lie derivative operators as Toeplitz operators (see Definition \ref{Cb2}). In \cref{Hd}, for $A\in\mathfrak{u}$, we get the asymptotics trace of ${e}^{\frac{A}{p}}$ on $H^{(0,0)}(N,L^p\otimes {{{\xi}}})$ using the Kirillov formula. In \cref{He}, we obtain an explicit formula for $W_1^{L,\xi}$. In \cref{H.e}, we discuss a special case of $W_1^{L,\xi}$ when $N$ is an adjoint orbit. In \cref{Hf}, we compute $W_1^{L,\xi}$ for some concrete examples.

\subsection{Reductive Lie groups}\label{Ha}

Let $G$ be a connected real reductive Lie group with Lie algebra $\mathfrak{g}$, and let $\Theta\in\mathrm{Aut}(G)$ be a Cartan involution of $G$. Let $K\subset G$ be the fixed point set of $\Theta$ in $G$, which is a maximal compact subgroup of $G$, and let $\mathfrak{k}$ be its Lie algebra. Let $\mathfrak{p}\subset\mathfrak{g}$ be the eigenspace of $\Theta$ associated with the eigenvalue $-1$. Then we have the Cartan decomposition
\begin{equation}\label{ha1}
	\mathfrak{g}=\mathfrak{p}\oplus\mathfrak{k}.
\end{equation}

Let $U$ be the compact form of $G$ with Lie algebra
\begin{equation}
	\mathfrak{u}=\sqrt{-1}\mathfrak{p}\oplus\mathfrak{k}.
\end{equation}
Let $G_\mathbb{C}$ be the complexification of $G$ with Lie algebra
\begin{equation}
	\mathfrak{g}_\mathbb{C}=\mathfrak{g}\otimes_\mathbb{R}\mathbb{C}=\mathfrak{u}\otimes_\mathbb{R}\mathbb{C}=\mathfrak{u}_\mathbb{C}
\end{equation}
Then $G_\mathbb{C}$ is also the complexification of $U$, and $U$ is a maximal compact subgroup of $G_\mathbb{C}$. 

\subsection{Principle bundles with reductive structure groups}\label{Hb}

Now we use the same assumptions as in \cref{Da}. In particular, $p\colon P_G\rightarrow M$ is a principal flat $G$-bundle where $G$ is a reductive Lie group given in \cref{Ha}. Since $G/K$ is contractible, there are smooth sections of the fibre bundle $P_G \times_G G/K$ and each corresponds to a reduction of the principal $G$-bundle $p\colon P_G\rightarrow M$ to a principal $K$-bundle $p\colon P_K\rightarrow M$.

Set $\mathfrak{g}_r=P_G\times_G \mathfrak{g}\cong P_K\times_K\mathfrak{g}$. Then \eqref{ha1} gives
\begin{equation}\label{hb3}
\mathfrak{g}_r=\mathfrak{p}_r\oplus \mathfrak{k}_r.
\end{equation}

We denote the flat connection on $P_G$ by a $\mathfrak{g}$-valued flat connection form $\theta^\mathfrak{g}$. By projection of $\theta^\mathfrak{g}$ on $\mathfrak{p}$ and $\mathfrak{k}$ with respect to the decomposition \eqref{ha1}, we get
\begin{equation}\label{hb4}
\theta^\mathfrak{g}=\theta^\mathfrak{p}+\theta^\mathfrak{k},
\end{equation}
then $\theta^\mathfrak{k}$ can be viewed as a connection form on $P_K$ and $\theta^\mathfrak{p}$ is a section of $T^*X\otimes\mathfrak{p}_r$.

Let $\Theta^\mathfrak{k}$ be the curvature of $\theta^\mathfrak{k}$. Since $\theta^\mathfrak{g}$ is flat, we get
\begin{equation}
\Theta^\mathfrak{k}=-\frac{1}{2}\big[\theta^\mathfrak{p},\theta^\mathfrak{p}\big]\ (\text{or\ }\Theta^\mathfrak{k}=-\theta^{\mathfrak{p},2}),\ \ \  \big[d+\theta^\mathfrak{k},\theta^\mathfrak{p}\big]=0.
\end{equation}
Let $\nabla^{\mathfrak{g}_r,u}$ be the connection on $\mathfrak{g}_r$ induced by $\theta^\mathfrak{k}$, which preserves the splitting \eqref{hb3}.

\subsection{Moment map and a class of Toeplitz operators}\label{Hc}

We use the assumptions and notation of \cref{Ca} and \cref{Da}. In particular, $N$ denotes a compact complex manifold of complex dimension $n$ and $(L,g^{L})$ is a holomorphic Hermitian line bundle on $N$. Let $\nabla^L$ be the Chern connection of $L$ with $c_1(L,g^L)=\omega$. 

We assume that the group $U$ acts holomorphically on $N$. If $A\in\mathfrak{u}$, let $A^{N}$ be the corresponding holomorphic vector field on $N$ with its $(1,0)$ part $A^{N,(1,0)}$ and $(0,1)$ part $A^{N,(0,1)}$), hence $A\in\mathfrak{u}\rightarrow -A^{N}$ is a morphism of Lie algebras.

We assume further that the action of $U$ on $N$ lifts to a holomorphic unitary action on $L$. Then $\omega$ is a $U$-invariant form. If $A\in\mathfrak{u}$, let ${L}_A$ denote the Lie derivative of $A$ on the smooth sections of $L$. If $A\in\mathfrak{u}$, the Kostant formula \cite[Definition 7.5]{bgv} gives
\begin{equation}\label{hc1}
{L}_A=\nabla^L_A-{2\pi i}\langle\mu_L,A\rangle,
\end{equation}
where $\mu_L$ is a moment map $\mu\colon N\rightarrow\mathfrak{u}^*$ such that, if $u\in U,x\in N$,
\begin{equation}
	\mu_L(ux)={}^t\mathrm{Ad}^{-1}(u)\mu(x),\ \ \ d\langle\mu_L,A\rangle-{i}_{A^{N}}\omega=0.
\end{equation}
Let $({{\xi}},g^{{\xi}})$ be another holomorphic  Hermitian line bundle on $N$ and $U$ acts holomorphically and unitarily on ${\xi}$ with moment map $\mu_{{\xi}}$.

Recall that $\mathscr{C}^\infty(N,L^p\otimes\xi)$ is equipped with the $L^2$-product induced by $(g^L,g^\xi)$, the Kähler metric is given by $g^{T_\mathbb{R}N}=\omega(\cdot,J\cdot)$ and $P_p$ is the orthogonal projection operator from $\mathscr{C}^\infty(N,L^p\otimes\xi)$ to $H^{(0,0)}(N,L^p\otimes\xi)$ as in \cref{Db}. Then $U$ acts on $H^{(0,0)}(N,L^p\otimes\xi)$ unitarily. This action can be extended to a holomorphic $G_\mathbb{C}$ action.

Let $g^{TN}$ be the metric on $TN$ induced by $g^{T_\mathbb{R}N}$ with the corresponding Chern connection $\nabla^{TN}$ and its curvature $R^{TN}$. Denote the $(1,0)$ part of $\nabla^{TN}$ by $\nabla^{TN}{'}$. If $A\in\mathfrak{u}$, $\nabla^{TN}{'} A^{N,(1,0)}$ is a skew-adjoint endomorphism of $TN$. 

The metric $g^{TN}$ induces a Hermitian metric $g^{\det TN}$ on the line bundle $\det TN$. Let $\nabla^{\det TN}$ denote the corresponding Chern connection on $\det TN$, then
\begin{equation}\label{hc9.}
c_1(\det TN,g^{\det TN})=-\frac{1}{{2\pi i}}\tro^{TN}[R^{TN}].
\end{equation}
The group $U$ acts holomorphically and unitarily on $\det TN$. The associated moment map $\mu_{\det TN}\colon N\rightarrow \mathfrak{u}^*$ is given by
\begin{equation}\label{hc4}
{2\pi i}\langle\mu_{\det TN},A\rangle=\tro^{TN}[\nabla^{TN}{'} A^{N,(1,0)}].
\end{equation}
\begin{theo}[{\cite[Theorem 3.1]{bmz17}}]
If $A\in\mathfrak{u}$, the following identity holds:
\begin{equation}\label{hc6}
\mathcal{L}^{}_A\big|_{H^{(0,0)}(N,L^p\otimes\xi)}=-{2\pi i} P_p\langle p\mu_L+\mu_\xi+\mu_{\det TN},A\rangle P_p.
\end{equation}
\end{theo}
\begin{pro}
By \eqref{da9} and \eqref{hc4}, we see that \eqref{hc6} is a special case of \eqref{db16}.\qed
\end{pro}

\subsection{The Kirillov formula}\label{Hd}

For a matrix $B$, if $\mathrm{sup}_{\lambda\in\mathrm{Spec}(B)}\vert\lambda\vert<2\pi$, we set
\begin{equation}\label{hd3}
	\mathrm{Td}(B)=\det\Big[\frac{B}{1-{e}^{-B}}\Big].
\end{equation}
If $B$ is self-adjoint, no condition on $B$ is necessary to define $\mathrm{Td}(B)$. From the above, for $A\in\mathfrak{g}_\mathbb{C}$ such that $\vert A\vert$ small, the following form is well defined:
\begin{equation}\label{hd4}
\mathrm{Td}_A(N)=\mathrm{Td}\Big[-\frac{R^{TN}}{{2\pi i}}+\nabla^{TN}{'} A^{N,(1,0)}\Big].
\end{equation}

The following form of the Kirillov formula is given in \cite[Theorem 3.5]{bmz17}.
\begin{theo}
For $p$ large enough, if $A\in\mathfrak{g}_\mathbb{C}$ and $\vert  A\vert$ small, 
\begin{equation}\label{hd6}
\tro^{H^{(0,0)}(N,L^p\otimes\xi)}({e}^A\big)=\int_N\mathrm{Td}_A(N)\exp\big({2\pi i} \langle p\mu_L+\mu_{\xi},A\rangle+pc_1(L,g^L)+c_1({{\xi}},g^{{{\xi}}})\big).
\end{equation}
\end{theo}

For $A\in\mathfrak{g}_\mathbb{C}$, let $R_{L}(A)$ be the integral of Duistermaat-Heckman given by
\begin{equation}\label{hd7}
R_{L}(A)=\int_N\exp\big({2\pi i}\langle\mu_L,A\rangle+c_1(L,g^L)\big).
\end{equation}

Formally, for $p\in\mathbb{N}^*$, let ${{{\xi}}^{\frac{1}{p}}}$ be the $p$-th root of $L$ and $(\det TN)^{\frac{1}{2p}}$ be the $2p$-th root of $\det TN$ at the level of cohomology and moment map, then
\begin{equation}\label{hd9}
	\begin{split}
		c_1\big({{{\xi}}^{\frac{1}{p}}}\big)&=\frac{1}{p}c_1({{\xi}},h^{{{\xi}}}),\ \ \ c_1\big((\det TN)^{\frac{1}{2p}}\big)=\frac{1}{2p}c_1(\det TN,g^{\det TN}),\\
		\mu_{{{{\xi}}^{\frac{1}{p}}}}&=\frac{1}{p}\mu_{{{\xi}}},\ \ \ \mu_{(\det TN)^{\frac{1}{2p}}}=\frac{1}{2p}\mu_{\det TN}.
	\end{split}
\end{equation}
\begin{prop}
If $A\in\mathfrak{g}_\mathbb{C}$, for $p$ large enough, we have
\begin{equation}\label{hd8}
p^{-n}\tro^{H^{(0,0)}(N,L^p\otimes\xi)}\big[{e}^{\frac{A}{p}}\big]=R_{L\otimes {{{\xi}}^{\frac{1}{p}}}\otimes (\det TN)^{\frac{1}{2p}}}(A)+\mathcal{O}\big(p^{-2}\big).
\end{equation}
\end{prop}

\begin{pro}
	Let $Q(z)$ be a holomorphic function near $0$ given by
	\begin{equation}\label{hd.10}
		Q(z)=\ln\frac{z}{1-{e}^{-z}}=\frac{z}{2}-\frac{z^2}{24}+\cdots.
	\end{equation}
By \eqref{hd3} and \eqref{hd.10}, we get $\mathrm{Td}(B)=\exp\big(\tro\big[Q(B)\big]\big)$. From \eqref{hc9.} and \eqref{hd4}, we have
\begin{equation}\label{hd10}
	\begin{split}
&\mathrm{Td}_{\frac{A}{p}}\big(N\big)=\mathrm{Td}_0(N)\Big(1+p^{-1}\tro\big[Q'({i}R^{TN}/2\pi)\nabla^{TN}{'}A^{N,(1,0)}\big]\Big)+\mathcal{O}\big(p^{-2}\big)\\
&=\Big(1+\frac{1}{2}c_1(\det TN)+\cdots\Big)\Big(1+p^{-1}\tro\big[Q'({i}R^{TN}/2\pi)\nabla^{TN}{'}A^{N,(1,0)}\big]\Big)+\mathcal{O}\big(p^{-2}\big).
	\end{split}
\end{equation}
Note that $Q'(0)=1/2$. By \eqref{cb9}, \eqref{hc4}, \eqref{hd6} and \eqref{hd10}, we get
\begin{equation}\label{hd12}
	\begin{split}
		p^{-n}\tro^{H^{(0,0)}(N,L^p\otimes\xi)}&\big[{e}^{\frac{A}{p}}\big]=\int_N\exp\big({2\pi i} \langle\mu_L,A\rangle\big)\frac{c_1^n(L)}{n!}\\
		+&\frac{1}{2p}\big({2\pi i} \langle\mu_{\det TN},A\rangle+4\pi{i}\langle\mu_{\xi},A\rangle\big)\exp\big({2\pi i} \langle\mu_L,A\rangle\big)\frac{c_1^n(L)}{n!}\\
		+&\frac{1}{2p} \big(c_1(\det TN)+2c_1({{\xi}})\big)\exp\big({2\pi i} \langle\mu,A\rangle\big)\frac{c^{n-1}_1(L)}{(n-1)!}+\mathcal{O}\big(p^{-2}\big).
	\end{split}
\end{equation}
By \eqref{hd9}, we could rewrite the right hand side of \eqref{hd12} as
\begin{equation}
\int_N\exp\Big({2\pi i} \big\langle\mu_{L\otimes {{{\xi}}^{\frac{1}{p}}}\otimes (\det TN)^{\frac{1}{2p}}},A\big\rangle+c_1\big(L\otimes {{{\xi}}^{\frac{1}{p}}}\otimes (\det TN)^{\frac{1}{2p}}\big)\Big)+\mathcal{O}\big(p^{-2}\big),
\end{equation}
which is exactly \eqref{hd8}.\qed
\end{pro}

\subsection{An explicit formula for $W_1^{L,\xi}$}\label{He}

For $p\in\mathbb{N}^*$, we denote by $\rho$ the holomorphic representation $G_\mathbb{C}\rightarrow \mathrm{End}\big(H^{(0,0)}(N,L^p\otimes{{\xi}})\big)$. Similar to \eqref{0.4...}, we have
\begin{equation}\label{hc7}
	{F}_p=P_K\times_KH^{(0,0)}(N,L^p\otimes {{\xi}}).
\end{equation}
Let $g^{{F}_p}$ be the Hermitian metric on $F_p$ induced by the $L^2$-metric on $H^{(0,0)}(N,L^p\otimes {{\xi}})$. The connection form $\theta^\mathfrak{g}$ induces a flat connection $\nabla^{F_p}$ on $F_p$, and $\theta^\mathfrak{k}$ gives a unitary connection $\nabla^{F_p,u}$ on $F_p$. Recall that $\theta^{\mathfrak{p}}$ is a section of $T^*M\otimes\mathfrak{p}_r$. By \eqref{hb4}, we get
\begin{equation}\label{hc8}
	\nabla^{F_p,u}=\nabla^{F_p}-\rho\theta^{\mathfrak{p}}
\end{equation}
where $\rho\theta^{\mathfrak{p}}\in\Omega^1(M,\mathrm{End}(F_p))$. By \eqref{bb3} and \eqref{hc8}, we have
\begin{equation}\label{hc9}
	\omega\big(\nabla^{F_p,u},g^{F_p}\big)=-2\rho\theta^{\mathfrak{p}},\ \ \ R^{F_p,u}=-\rho\theta^{\mathfrak{p},2}.
\end{equation}
Note that \eqref{hc8} and \eqref{hc9} are reductive version of Theorem \ref{Db4}. Also, by \eqref{hc6} and \eqref{hc9}, for a local orthonormal frame $\{e_i\}_{i=1}^m$ of $TX$, the condition \eqref{de1} becomes
\begin{equation}\label{hc11}
	\begin{split}
		\sum_{i=1}^m\bv\langle2\pi\mu_L,{\sqrt{-1}\widehat{\theta}^\mathfrak{p}(e_i)}\rangle\bv^2>0.
	\end{split}
\end{equation}

Let $S\mathfrak{g}$ be the symmetric algebra of $\mathfrak{g}$, which can be viewed as the algebra of real differential operators with constant coefficients on $\mathfrak{g}$, and we denote by $\overline{S}\mathfrak{g}$ its formal completion. Let $\widehat{\theta}^\mathfrak{p}$ be the restriction of ${\theta}^\mathfrak{p}$ on $\widehat{TX}$. For $t\geqslant0$, recall that $\sigma_t$ is a section of $\Lambda(T^*M)\hti\Lambda(\widehat{T^*X})\otimes S\mathfrak{g}_r$ given in \eqref{hd..14}.

Now we state and prove the main result of this section.
\begin{theo}\label{Hd4}
The first two terms of the asymptotic torsion \eqref{ga7} are given by
\begin{equation}\label{hd15}
	\begin{split}
p^{-n-1}\psi_{1/\sqrt{p}}\mathcal{T}(T^HM&,g^{TX},\nabla^{F_p},g^{F_p})=W^{L,{{\xi}}}_0+p^{-1}W_1^{L,{{\xi}}}+\mathcal{O}\big(p^{-2}\big)\\
&=\sqrt{{2\pi i}}\varphi\int_{0}^{+\infty}\int^{\widehat{B}}\frac{\theta^{\mathfrak{p}}\wedge\widehat{\theta^{\mathfrak{p}}}}{2}\exp(-\sigma_t) R_{L\otimes{{{\xi}}}^{\frac{1}{p}}\otimes (\det TN)^{\frac{1}{2p}}}(0)\frac{dt}{\sqrt{t}},
	\end{split}
\end{equation}
which also could be formally written in the following form:
	\begin{equation}\label{hd16}
		W_0^{L,{{\xi}}}+p^{-1}W_1^{L,{{\xi}}}=W_0^{L\otimes {{{\xi}}}^{\frac{1}{p}}\otimes (\det TN)^{\frac{1}{2p}}}+\mathcal{O}\big(p^{-2}\big).
	\end{equation}
\end{theo}

\begin{pro}
	Let us first explain \eqref{hd15} in more detail. Note that $\frac{\theta^{\mathfrak{p}}\wedge\widehat{\theta^{\mathfrak{p}}}}{2}\exp(-\sigma_t)$ is a smooth section of $\Lambda(\widehat{T^*X})\hti \Lambda(T^*M)\otimes \overline{S}\mathfrak{g}_r$, and it acts naturally on $R_{L\otimes{{{\xi}}}^{\frac{1}{p}}\otimes (\det TN)^{\frac{1}{2p}}}(\cdot)$, which is a function on a neighborhood of $0\in \mathfrak{g}_\mathbb{C}$ (see \cite[\S\,1.4]{bmz17} for more details). When evaluating at $0\in\mathfrak{g}_\mathbb{C}$, we see that $\frac{\theta^{\mathfrak{p}}\wedge\widehat{\theta^{\mathfrak{p}}}}{2}\exp(-\sigma_t)R_{L\otimes{{{\xi}}}^{\frac{1}{p}}\otimes (\det TN)^{\frac{1}{2p}}}(0)$ is a smooth section of $\Lambda(\widehat{T^*X})\hti \Lambda(T^*M)$.
	
Now we come back to the proof of \eqref{hd15}. We plug \eqref{hc8} and \eqref{hc9} into \eqref{gb11}, then by \eqref{gb12} and \eqref{hd8},  we get
\begin{equation}\label{hd17}
	\begin{split}
		d_0(t)+p^{-1}d_1(t)=\int^{\widehat{B}}\sqrt{t}\frac{\theta^{\mathfrak{p}}\wedge\widehat{\theta^{\mathfrak{p}}}}{2}\exp(-\sigma_t) R_{L\otimes {{{\xi}}^{\frac{1}{p}}}\otimes (\det TN)^{\frac{1}{2p}}}(0)+\mathcal{O}\big(p^{-2}\big),
	\end{split}
\end{equation}
which gives \eqref{hd15} together with \eqref{gb.15}.\qed
\end{pro}

\subsection{A special case for coadjoint orbit}\label{H.e}

Now we discuss a special case of Theorem \ref{Hd4} when $N$ is given by a coadjoint orbit. 

For the compact connected Lie group $U$ given in \cref{Ha} with Lie algebra $\mathfrak{u}$. For $\lambda\in\mathfrak{u}^*$, denote by $\mathcal{O}_\lambda=U\cdot\lambda$ the orbit of coadjoint action with the natural symplectic form $\omega_\lambda$: for $A,B\in\mathfrak{u}$ and $f\in\mathcal{O}_\lambda$,
\begin{equation}
\omega_\lambda(A^{\mathcal{O}_\mu},B^{\mathcal{O}_\mu})_{f}=\frac{1}{2\pi}\langle f,[A,B]\rangle,
\end{equation}
then the associated moment map $\mu_\lambda$ verifies that $2\pi\mu_\lambda$ is the inclusion 
\begin{equation}\label{he..1}
2\pi\mu_\lambda\colon\mathcal{O}_\lambda\hookrightarrow\mathfrak{u}^*.
\end{equation}
Similar to \eqref{hd7}, for $A\in\mathfrak{u}$, set
\begin{equation}\label{he.18}
	\begin{split}
		R_\lambda(A)=\int_{\mathcal{O}_\lambda}\exp\big({2\pi i}\langle \mu_\lambda,A\rangle+\omega_\lambda\big).
	\end{split}
\end{equation}

We fix a maximal torus $T$ of $U$ with Lie algebra $\mathfrak{t}$. We note that if $\lambda\in\mathfrak{u}^*$ is regular, we have $\mathcal{O}_\lambda\cong U/T$. The integral lattice $\Lambda\subset \mathfrak{t}$ is defined as the kernel of the exponential map $\exp\colon\mathfrak{t}\to T$, and the real weight lattice $\Lambda^*\subset \mathfrak{t}^*$ is defined by $\Lambda^*=\mathrm{Hom}(\Lambda,2\pi\mathbb{Z})$. We fix a set of positive roots $\Phi^+\subset\Lambda^*$, a positive open Weyl Chamber $\mathcal{C}^+$ and its closure $\overline{\mathcal{C}^+}$. By the Weyl character formula, finite-dimensional irreducible representations of $U$ are parameterized by $\Lambda^*\cap\overline{\mathcal{C}^+}$.

Set $\mathfrak{r}=[\mathfrak{t},\mathfrak{u}]$, then we have $\mathfrak{u}=\mathfrak{t}\oplus \mathfrak{r}$, and $\mathfrak{u}^*=\mathfrak{t}^*\oplus \mathfrak{r}^*$. Hence we identify $\Lambda^*\cap\overline{\mathcal{C}^+}$ with a subset of $\mathfrak{u}^*$. For $\lambda\in \Lambda^*\cap\overline{\mathcal{C}^+}$, let $V_\lambda$ be the unique finite dimensional irreducible representation of $U$ with highest weight $\lambda$.

If $\lambda\in \Lambda^*\cap\overline{\mathcal{C}^+}$, let $(L_\lambda,g^{L_\lambda})$ be the canonical prequantum line bundle on $\mathcal{O}_\lambda$ with $c_1(L_\lambda,g^{L_\lambda})=\omega_\lambda$. By the Borel-Weil-Bott theorem \cite[Theorem 8.8]{bgv}, we have an isomorphism
$H^{(0,0)}(\mathcal{O}_\lambda,L_\lambda)\cong V_\lambda$.

For $\lambda\in\Lambda^*\cap\overline{\mathcal{C}^+}$ and $\tau\in\Lambda^*$, in Theorem \ref{Hd4} we take
\begin{equation}
(N,L,{{\xi}},H^{(0,0)}(N,L^p\otimes {{\xi}}),F_p)=(U/T,L_\lambda,L_\tau,V_{p\lambda+\tau},P_G\times_GV_{p\lambda+\tau}).
\end{equation}

Set $\varrho_\mathfrak{u}=\frac{1}{2}\sum_{\alpha\in \Phi^+} \alpha$, then $\det T(U/T)\cong L_{2\varrho_\mathfrak{u}}$. Hence, by \eqref{hd8} and \eqref{he.18} we have
\begin{equation}\label{he23}
	p^{-n}\tro^{V_{p\lambda+\tau}}\big[{e}^{\frac{A}{p}}\big]=R_{\lambda+p^{-1}(\varrho_{\mathfrak{u}}+\tau)}(A)+\mathcal{O}\big(p^{-2}\big).
\end{equation}
Note that we could also derive \eqref{he23} from the Kirillov character formula \cite[Theorem 8.4]{bgv}: for any $A\in\mathfrak{u}$, let $\mathrm{ad}(A)$ be the morphism of $\mathfrak{u}$ given by $\mathrm{ad}(A)B=[A,B]$, then
\begin{equation}\label{he24}
	\tro^{V_\lambda}\big[{e}^A\big]=j^{-\frac{1}{2}}(A)R_{\lambda+\varrho_{\mathfrak{u}}}(A),\ \ \text{for\ }j(A)=\text{det}\Big(\frac{\sinh(\mathrm{ad}(A)/2)}{\mathrm{ad}(A)/2}\Big).
\end{equation}
Since $\mathrm{ad}(A)$ is anti-symmetric, $j^{-\frac{1}{2}}(\cdot)$ is an even function of $A\in \mathfrak{u}$ and $j^{-\frac{1}{2}}(0)=1$, we have $j^{-\frac{1}{2}}(A/p)=1+\mathcal{O}(p^{-2})$. Then by \eqref{he.18} and \eqref{he24}, we recover \eqref{he23}:
\begin{equation}
	\begin{split}
		p^{-n}\tro^{V_{p\lambda+\tau}}\big[{e}^\frac{A}{p}\big]=p^{-n}\big(1+\mathcal{O}(p^{-2})\big)R_{p\lambda+\tau+\varrho_{\mathfrak{u}}}\Big(\frac{A}{p}\Big)=R_{\lambda+p^{-1}(\varrho_{\mathfrak{u}}+\tau)}(A)+\mathcal{O}\big(p^{-2}\big).
	\end{split}
\end{equation}

By \eqref{hc11} and \eqref{he..1}, $\theta^\mathfrak{g}$ is nondegenerated with respect to $\lambda$ if and only if for $u\in U$,
\begin{equation}\label{he.19}
	\sum_{i=1}^m\bv\langle u\cdot2\pi\lambda,{\sqrt{-1}\widehat{\theta}^\mathfrak{p}(e_i)}\rangle\bv^2>0,
\end{equation}
and when \eqref{he.19} holds, we set
\begin{equation}
W^\lambda=\sqrt{{2\pi i}}\varphi\int_{0}^{+\infty}\int^{\widehat{B}}\frac{\theta^{\mathfrak{p}}\wedge\widehat{\theta^{\mathfrak{p}}}}{2}\exp(-\sigma_t) R_{\lambda}(0)\frac{dt}{\sqrt{t}},
\end{equation}
then we have the following corollary from \eqref{he23} exactly as \eqref{hd16} follows from \eqref{hd8}:
\begin{coro}\label{He6}
If  $\theta^\mathfrak{g}$ is nondegenerated with respect to $\lambda\in\Lambda^*\cap\overline{\mathcal{C}^+}$ and $\tau\in\Lambda^*$, we have
	\begin{equation}\label{he21}
		\begin{split}
p^{-n-1}\psi_{1/\sqrt{p}}\mathcal{T}(T^HM,g^{TX},\nabla^{F_p},g^{F_p})=\int_XW^{\lambda+p^{-1}(\varrho_\mathfrak{u}+\tau)}+\mathcal{O}\big(p^{-2}\big).
		\end{split}
	\end{equation}
\end{coro}

\subsection{The $W_1$ for $\mathrm{SL}(2,\mathbb{C})$}\label{Hf}

In this subsection, we mainly follow \cite[\S\,8.7]{bmz17}. Now, we assume that $G=\mathrm{SL}(2,\mathbb{C})$. The Cartan involution is $\Theta\colon g\rightarrow (g^*)^{-1}$, then $K=\mathrm{SU}(2)$. Note that $U=\mathrm{SU}(2)\times\mathrm{SU}(2)$ is a compact form of $G=\mathrm{SL}(2,\mathbb{C})$, in which $K=\mathrm{SU}(2)$ embeds by the diagonal embedding.

A maximal torus $T=S^1$ in $\mathrm{SU}(2)$ is the $1$-parameter group $\exp(2\pi t\bm{k}),t\in S^1=\mathbb{R}/\mathbb{Z}$ where $\bm{k}$ is a generator of the Lie algebra of $\mathfrak{t}$, and $2\pi \bm{k}$ is a coroot in $T$. Then
\begin{equation}\label{hf9}
 \Phi^+=\big\{\frac{\bm{k}}{\pi}\big\}, \ \ \ \ \varrho_\mathfrak{u}=\frac{\bm{k}}{2\pi},
\end{equation}
and the set $\Lambda^*\cap\overline{C^+}$ of dominant weights is just given by $\mathbb{N}/2\pi\cdot\bm{k}$.

For $a\in \mathbb{N}$, the coadjoint orbit $N_a$ of $a\bm{k}/2\pi\in\mathfrak{t}^*$ for $\mathrm{SU}(2)$ can be identified with a point for $a=0$, with $\mathbb{CP}^1$ for $a>0$. The orbit $\mathcal{O}_a$ carries a canonical line bundle $L_a$. 

For $a,b\neq 0$, put $\mathcal{O}_{a,b}=\mathcal{O}_a\times \overline{\mathcal{O}}_b$. Let $q_1,q_2$ be the obvious projections from $\mathcal{O}_{a,b}$ on $\mathcal{O}_{a}$ and $\overline{\mathcal{O}}_b$. For $a,b\in\mathbb{N}$, set $L_{a,b}=q_1^*L_a\times q_2^*\overline{L}_b$. Then $g\in\mathrm{SL}(2,\mathbb{C})$ acts on $\mathcal{O}_{a,b}$ by the map $(z,z')\rightarrow(gz,(\Theta g)z')$, and the action of $\mathrm{SL}(2,\mathbb{C})$ lifts to $L_{a,b}$. Moreover,
\begin{equation}
H^{(0,0)}(\mathcal{O}_{a,b},L_{a,b})\cong \text{Sym}^{a}\mathbb{C}^2\otimes \text{Sym}^{b}\overline{\mathbb{C}}^2,
\end{equation}
where $\mathbb{C}^2$ is the tautological representation of $\mathrm{SL}(2,\mathbb{C})$, $\overline{\mathbb{C}}^2$ is its complex conjugate representation and $\text{Sym}^k$ means the $k$-th symmetric power.

Now we give the change of notation concerning the previous sections. Note that $G/K=\mathbb{H}^3$, the hyperbolic space with constant sectional curvature $-4$. Let $\Gamma$ be a torsion-free discrete cocompact subgroup of $\mathrm{SL}(2,\mathbb{C})$. Then $\Gamma\backslash\mathbb{H}^3$ is a compact $3$-dimensional hyperbolic manifold. For $a,b\in\mathbb{N}^*,a',b'\in\mathbb{Z}$ with $a\neq b$, we set
\begin{equation}
(S,X,P_G,P_K,N,L,{{\xi}})=(\{\mathrm{pt}\},\Gamma\backslash\mathbb{H}^3,\Gamma\backslash(\mathbb{H}^3\times G),\Gamma\backslash G/K,\mathcal{O}_{a,b},L_{a,b},L_{a',b'}).
\end{equation}
As in \eqref{hc7}, we have
\begin{equation}
	\begin{split}
		F_p=\Gamma\backslash G\times_KH^{(0,0)}(N,L^p\otimes {{\xi}})\cong \Gamma\backslash G\times_K \big(\text{Sym}^{ap+a'}\mathbb{C}^2\otimes \text{Sym}^{bp+b'}\overline{\mathbb{C}}^2\big),
	\end{split}
\end{equation}
then we could work under the settings of Corollary \ref{He6}.

\begin{prop}
The degeneracy condition \eqref{he.19} is equivalent to $a\in\mathbb{N}^*,b=0$ or $a,b\in\mathbb{N}^*,a\neq b$, and we have
\begin{equation}\label{hf14}
	\begin{split}
W^{L_{a,0}}_0=\frac{2}{\pi}a^2,\ \ \ \		W^{L_{a,b}}_0=\begin{cases}
			\frac{2}{3\pi}(3a^2b-b^3), \ a>b,\\
			\frac{2}{3\pi}(3ab^2-a^3), \ a<b.
		\end{cases}
	\end{split}
\end{equation}
\end{prop}

The first equality in \eqref{hf14} was the main result obtained by Müller \cite[Theorem 1.1]{mu12} and Bismut-Ma-Zhang \cite[Theorem 8.1]{bmz17} gave both equalities in \eqref{hf14}. We note that in \cite{mu12}, the curvature of $\mathbb{H}^3$ is $-1$ instead of $-4$.

For $a\in\mathbb{N}^*,a'\in\mathbb{Z},b=b'=0$ or $a,b\in\mathbb{N}^*,a\neq b,a',b'\in\mathbb{Z}$, set
\begin{equation}\label{hf16}
	\begin{split}
	W^{L_{a,0},L_{a',0}}_1&=\frac{4}{\pi}a(a'+1),\\
	W^{{L_{a,b},L_{a',b'}}}_1&=\begin{cases}
			\frac{2}{\pi}\big((a^2-b^2)(b'+1)+2ab(a'+1\big),\  a>b,\\
			\frac{2}{\pi}\big((b^2-a^2)(a'+1)+2ab(b'+1\big),\  a<b.
		\end{cases}
	\end{split}
\end{equation}
\begin{theo}\label{Hf9}
For $a,b\in\mathbb{N}^*,a',b'\in\mathbb{Z}, a\neq b$, as $p\rightarrow +\infty$, we have
\begin{equation}\label{hf17}
	p^{-3}\mathcal{T}(\Gamma\backslash\mathbb{H}^3,{F_p})=\big(W^{{L_{a,b},L_{a',b'}}}_0+p^{-1}W^{{L_{a,b},L_{a',b'}}}_1\big)\mathrm{Vol}(\Gamma\backslash\mathbb{H}^3)+\mathcal{O}\big(p^{-2}\big).
\end{equation}
For $a\in\mathbb{N}^*,a'\in\mathbb{Z},b,b'=0$, as $p\rightarrow +\infty$,
\begin{equation}\label{hf18}
	p^{-2}\mathcal{T}(\Gamma\backslash\mathbb{H}^3,{F_p})=\big(W^{L_{a,0},L_{a',0}}_0+p^{-1}W^{L_{a,0},L_{a',0}}_1\big)\mathrm{Vol}(\Gamma\backslash\mathbb{H}^3)+\mathcal{O}\big(p^{-2}\big).
\end{equation}
\end{theo}

\begin{pro}
By \eqref{hf14}, we get the first order expansion in \eqref{hf16}. By Corollary \ref{He6} and \eqref{hf9}, we replace $a,b$ with $a+p^{-1}(1+a'),b+p^{-1}(1+b')$ in $W_0^{L_{a,0}}$ and $W_0^{L_{a,b}}$ and evaluate the coefficients of $p^{-1}$, then we get the second orders in the expansions \eqref{hf16}.\qed
\end{pro}

\begin{remark}
By Theorem \ref{Gc5}, we see that $\mathcal{T}(\Gamma\backslash\mathbb{H}^3,{F_p})-\mathcal{T}_{L^2}(\Gamma\backslash\mathbb{H}^3,{F_p})=\mathcal{O}({e}^{-c{p}})$. By \cite[Page.423, Example (3)]{bergeron_venkatesh_2013}, if $F=\text{Sym}^{a}\mathbb{C}^2\otimes \text{Sym}^{b}\overline{\mathbb{C}}^2$, we have
\begin{equation}\label{hf.20}
	\begin{split}
\frac{\mathcal{T}_{L^2}(\Gamma\backslash\mathbb{H}^3,{F})}{\mathrm{Vol}(\Gamma\backslash\mathbb{H}^3)}=\frac{1}{6\pi}\big\{&(a+b+2)^3-\lv a-b\rv^3\\
&+3(a+b+2)\lv a-b\rv(a+b+2-\lv a-b\rv)\big\}.
	\end{split}
\end{equation}
Replacing $(a,b)$ with $(ap+a',0)$ or $(ap+a',bp+b')$ in \eqref{hf.20}, we recover Theorem \ref{Hf9}.
\end{remark}

\providecommand{\bysame}{\leavevmode\hbox to3em{\hrulefill}\thinspace}
\providecommand{\MR}{\relax\ifhmode\unskip\space\fi MR }
\providecommand{\MRhref}[2]{%
	\href{http://www.ams.org/mathscinet-getitem?mr=#1}{#2}
}
\providecommand{\href}[2]{#2}


\end{document}